\pdfoutput=1
\documentclass[11pt]{article}
\usepackage[english]{babel}
\usepackage{amsmath,amssymb,graphicx,latexsym,theorem,bbm}%,showkeys,showlabels} %\showlabels{bibitem}
\usepackage{mathabx}
\usepackage{xcolor}
\definecolor{blue}{rgb}{0,0,1}
\definecolor{red}{rgb}{1,0,0}
\def\red{\begin{color}{red}}
\def\ered{\end{color}}
\usepackage{times}
\usepackage[textsize=footnotesize,color=blue!40, bordercolor=white]{todonotes}
\usepackage[pagebackref,bookmarks]{hyperref} 
\hypersetup{pagebackref=true,bookmarks=true,
       hyperfigures=true,
	unicode=true,
	colorlinks=true,
	citecolor=blue,
	linkcolor=blue,
	anchorcolor=red
}

\newcommand{\assign}{:=}
\newcommand{\nobracket}{}
\newcommand{\nocomma}{}
\newcommand{\noplus}{}
\newcommand{\nosymbol}{}
\newcommand{\tmcolor}[2]{{\color{#1}{#2}}}
\newcommand{\tmem}[1]{{\em #1\/}}
\newcommand{\tmmathbf}[1]{\ensuremath{\boldsymbol{#1}}}
\newcommand{\tmop}[1]{\ensuremath{\operatorname{#1}}}

\newcommand{\tmstrong}[1]{\textbf{#1}}
\newcommand{\tmtextbf}[1]{\text{{\bfseries{#1}}}}
\newcommand{\tmtextit}[1]{\text{{\itshape{#1}}}}

\newtheorem{lemma}{\sc{Lemma}}
\numberwithin{lemma}{section}
\newtheorem{corollary}[lemma]{\sc Corollary}
\newtheorem{definition}[lemma]{\sc Definition}
{\theorembodyfont{\rmfamily}}
{\theorembodyfont{\rmfamily\small}}
\numberwithin{exercise}{section}

{\theorembodyfont{\rmfamily}\newtheorem{note}[lemma]{\sc Note}}
{\theorembodyfont{\rmfamily}\newtheorem{remark}[lemma]{\sc Remark}}
\newtheorem{theorem}[lemma]{\sc Theorem}

\newcommand{\nonconverted}[1]{\mbox{}}
\def\qed{\hfill{Q.E.D.}}
\def\Lim{\mathop{Lim}}

\newcommand{\pI}[1]{\{{#1}\}^\mathsf{I}}
\def\mr{\mathsf{R}}
\newcommand{\ptwoe}[1]{\{{#1}\}^\mathsf{E}}
\def\etwoe{\ptwoe}
\usepackage{xcolor}
\definecolor{blue}{rgb}{0,0,1}

\def\ed{\end{document}}
\def\g{\gamma}
\def\d{\delta}
\def\r{\rho}
\def\a{\alpha}
\def\w{\omega}
\def\b{\beta}
\def\s{\sigma}
\def\S{\Sigma}
\def\l{\lambda}
\def\z{\zeta}

\def\me{{\mathsf{E}}}
\def\mr{\mbox{$\mathsf{R}$}}
\def\mi{\mathsf{I}}    %{{\mbox{$\mathsf{I}$}}}
\def\mj{{{\mathsf{J}}}}
\def\mf{{\mathsf{F}}}
\def\mg{{\mbox{$\mathsf{G}$}}}
\def\mk{{\mbox{$\mathsf{K}$}}}

\def\mh{\mbox{$\mathsf{H}$}}
\def\ed{\end{document}}
\renewcommand{\qed}{\hfill QED}

\def\asterisk{\ast}
\def\ptt{\mbox{\em{II}}}
\def\poo{\mbox{\em{I}}}
\def\po{ ${\poo}$ }
\def\pt{${\ptt}$ }

\def\calr{\mathcal{R}}

\newcommand{\call}{ {\mathcal{L}}}

\newcommand{\gott}{\ensuremath{\mathfrak{T}}}

\newcommand{\zx}[1]{\ensuremath{{^{#1}\zeta}}}

\newcommand{\sx}[1]{\ensuremath{{^{#1} \Sigma}}}

\newcommand{\phiabp}[2]{\ensuremath{\Phi_{\beta +1}^{(  \alpha )}}}

\newcommand{\ptwo}{{\mbox{{\em II\/}}}}
\newcommand{\pone}{{\mbox{{\em I\/}}}}

\newcommand{\nod}{{\noindent}}

\newcommand{\ul}{\ensuremath{\ulcorner}}
\newcommand{\ur}{\ensuremath{\urcorner}}

\newcommand{\bu}{\ensuremath{\bullet}}
\newcommand{\nat}{\ensuremath{\mathbbm{N}}}
\newcommand{\re}{\ensuremath{\mathbbm{R}}}
\newcommand{\baire}{\ensuremath{\mathbbm{N}}\textsuperscript{\ensuremath{\mathbbm{N}}}}
\newcommand{\bai}{\textsuperscript{\ensuremath{\omega}}\ensuremath{\omega}}

\newcommand{\cant}{2\textsuperscript{\ensuremath{\mathbbm{N}}}}

\newcommand{\rest}{{\upharpoonright}}
\newcommand{\emp}{{\varnothing}}
\newcommand{\pa}[1]{\ensuremath{\langle #1 \rangle}}

\newcommand{\game}{\ensuremath{\Game}}

\newcommand{\vp}{\ensuremath{\varphi}}

\newcommand{\da}{{\downarrow}}
\newcommand{\ua}{{\uparrow}}
\newcommand{\la}{{\langle}}
\newcommand{\ra}{{\rangle}}
\newcommand{\dom}{\text{dom}}
\newcommand{\back}{{\backslash}}
\newcommand{\power}{{\mathcal{P}}}

\newcommand{\ie}{{\itshape{i.e.}}{\hspace{0.25em}}}
\newcommand{\eg}{{\itshape{e.g.}}{\hspace{0.25em}}}
\newcommand{\etc}{{\itshape{etc.}}{\hspace{0.25em}}}
\newcommand{\via}{{\itshape{via }}{\hspace{0.25em}}}
\newcommand{\cf}{{\itshape{cf. }}{\hspace{0.1em}}}

\newcommand{\all}{{\forall}}
\newcommand{\Equi}{{\,\Longleftrightarrow\,}}
\newcommand{\equi}{{\,\longleftrightarrow\,}}
\newcommand{\ex}{{\exists}}
\newcommand{\Imp}{{\,\Rightarrow\,}}
\newcommand{\sset}{ { \,\subseteq\, } }
\newcommand{\lr}{{\leftrightarrow}}
\newcommand{\imp}{{\,\longrightarrow\,}}

\newcommand{\rem}{{\noindent}{\bfseries{Remark: }}}
\newcommand{\fin}[1]{\ensuremath{[  ]^{< \omega}}}

\newcommand{\pf}{{\noindent}{\textbf{Proof: }}}
\newcommand{\dfs}{\ensuremath{=_{\ensuremath{\operatorname{df}}}}}
\newcommand{\edfs}{\ensuremath{  \Equi_{\ensuremath{\operatorname{df}}}  }}

\newcommand{\oddpagetext}[1]{\newcommand{\pageoddheader}{{\small }}}
\newcommand{\evenpagetext}[1]{\newcommand{\pageevenheader}{{\small }}}

\title{Higher Type ITTM-recursion and Determinacy of Infinite games}

\author{P.D. Welch\\
University of Bristol}
\date{May 15, 2026}

\begin{document}

\maketitle

\abstract{We outline a theory of type 2 recursion for Infinite Time Turing Machines {\em \`a la Kleene}. We establish a connection between classical descriptive set theory and ittm theory, by calculating the complexity of its halting problem as exactly that of a complete $\Game \Sigma^0_3$ (or $G_{\d\s}$) set. This mirrors exactly what Kleene, Moschovakis {\em et al.} achieved for Kleene's type 2 recursion and $\Sigma^0_1$ (or Open) Determinacy.}
We ascertain the least ordinal which is not generalised recursive in this sense, and its characterisation {\via}a concept of {\em infinite nestings} in G\"odel's constructible hierarchy. The results do not require large cardinal axioms, and are all provable within analysis.\footnote{{\em Key Words and Phrases}: Higher type recursion, infinite time Turing machine, infinite games, determinacy, constructibility.
{\em MSC2020 Mathematics Subject Classification: Primary} 03D65, 03D10, 03D78, 03E60, 0345, {\em Secondary:} 03D75, 68Q10.
}

\setcounter{tocdepth}{3}
\tableofcontents
%%%%%%%%%%%%%%%%
%\ed

\section{Introduction}

%%%%%%%%%%%%%%%
The purpose of this paper is to outline a theory of higher type recursion
(actually just type 2) for infinite time Turing machines (ittm's) in the
manner of Kleene from the late 1950's and early 1960's, in particular where he
used (ordinary) Turing machines arranged on wellfounded trees to present
computations and their subcomputations (see {\cite{Kl62b}}, {\cite{Kl62a}}).
What is {\tmem{ not}} meant here, is the computation on elements
of Cantor space, with an oracle some $A \subseteq \cant$ that was also outlined in
{\cite{HL}}. Although the basic ittm model allows for computation
using reals (identified with elements of Cantor space), that simple recursive
model does not have the Kleenean feature of calling {\em subcomputations}, which
enabled him to build up a model of recursion using Turing machines that, for
example, showed that the Kleene-decidable sets coincided with the hyperarithmetic (or $\Delta^1_1$):
we want to also build a model that equally goes beyond the basic oracle
machine.

We outline a theory of this{\tmem{ generalised type-2 ittm recursion}} here. A
number of choices were made as to what features the model should have, and no
doubt there are variations. A typical oracle in this theory is then a total
type-2 functional $\mathsf{I} : \,^k \omega \times \hspace{0.17em}^l
(^{\omega} 2) \rightarrow \omega$ - just as for Kleene type-2 recursion, and
we develop the theory of recursive-in-$\mathsf{I}$ functions
$\{e\}^{{\mathsf{I}}} : \,^k \omega \times \hspace{0.17em}^l (^{\omega}
2) \rightarrow \omega$ as exemplified by a machine architecture
$P_e^{\mathsf{I}} (\mathbf{m}, \mathbf{x})$ with $e$ thought of as an
index of a program - as always a finite ordinary Turing programme, but
enhanced with a query call instruction initiating subcomputations, just as Kleene
did.

Once the model is defined, the question of its properties arises. For Kleene
recursion the existential quantifier $^2\mathsf{E}$ played a crucial role,
and the theory was developed often with the stipulation that `normal'
functionals should be candidates for oracles, {\ie} that an oracle
$\mathsf{I}$ was normal, if $\mathsf{^2 E}$ was itself Kleene recursive in
$\mathsf{I}$. \ \ We use $\mathsf{E} $= $\mathsf{^2 E}$ here also as a `base'
oracle for our arguments, although here its role is rather trivial.

One immediately wants to ask: what are the semi-decidable sets? What is $H$, the
`$\mathsf{E}$-halting problem'?
$$H = H^\me = \left\{ e \mid \{ e \}^{\mathsf{E}} (e) \da \right\}.$$

The analogue of $\omega^{\tmop{ck}}_1 = \omega_1^{\mathsf{E}, \tmop{ck}}$ -
the least ordinal not the order type of a (Kleene) recursive in $\mathsf{E}$
wellordering of $\omega$, is here $\alpha_0^{\mathsf{E}}$, the least ordinal
not the order type of a wellordering of $\omega$ which is generalised ittm
recursive in $\mathsf{E}$. It turns out that this is a rather large countable
ordinal.

Infinite Time Turing Machine theory has developed since its inception in
{\cite{HL}} and there have been many ingenious ways of expanding its role to
define `computable' concepts. It was clear from {\cite{HL}} that the model
could produce codes for some initial segment of the $L_{\alpha}$-hierarchy.
One generalisation was due to Peter Koepke and simply extends the tape to have
length $\tmop{On}$, \cf{\!\!\cite{K05}}. This produces codes for any $L_{\alpha}$ -
given a mark for the ordinal $\alpha$ on its tape. Thus a satisfaction
relation for $(L, \in)$ could be so computable, whilst the original ITTM's
were limited to doing this for an initial segment of true $\omega_1$. There
had been the hope that such devices could give, when not new proofs of fine
structural results of Jensen, such as say the $\Box$
principle, at least new insight into the fine-structure of $L$. \ \ However
this seems not really to be the case: the $\Sigma_2$-nature of the liminf rule
used for cell update at limit stages works against the very $\Sigma_1$, or
even $\Sigma_0$ method of skolem hulls used in fine-structural proofs. So
although $\alpha$-recursion theory can be somewhat reformulated as
$\alpha$-length tape ittm models, no new deeper fine-structural insights seemed
to be forthcoming. It seemed that direct connection, or elucidation, of contingent areas of `classical' constructibility theory or  descriptive set theory was lacking.

However the author has for some while had the thought that the strength of
ittm theory was either exactly, or at least close to, being lined up with that
of $\Sigma^0_3$-Determinacy in the area of Gale-Stewart games, that is, two person perfect information games played on integers. (\cf\!\!, \eg\!\!,\cite{Mosch4}. The reasons for this insight are a little
difficult to state just here.) \ If this were to be so, it would be a
characterisation of this classical descriptive set theoretic property by means
of this generalised machine theory: it would also be the occasion of a true
application of ittm-theory to a problem in classical descriptive set theory, and as far as we are aware, it
would be the first such.

This can be realised as follows:

\ 

\nod{\sc Theorem \ref{Psicomplete}}\label{Thm1.1}
 {\em  Let $G_3 $ be a complete $\game \Sigma^0_3$-set. \ Then $H \equiv_1 G_3$, that
  is, the halting set for generalised type-2 ittm recursion, is recursively isomorphic to a complete $\game
  \Sigma^0_3$-set of integers.}\\
%\end{theorem}

That is to say, if the $\Sigma^0_3$ sets of reals are (ordinarily) recursively
listed as say $A_0, A_1, \nocomma \ldots \nocomma, A_n, \ldots$, then there is
a pencil and paper algorithm $f$ given an element $k \in H$ to calculate an $f
(k)$ so that Player \tmtextit{I} has a winning strategy in $A_{f (k)}${\tmem{
and vice versa}}. That is, $f : \omega \imp \omega$ is an ordinary Turing
computable bijection and we thus have, given $H,$ a complete listing of those
games in which $\mathit{I}$ wins, and conversely from such a listing we can
retrieve $H$ using the inverse of $f$. Thus $H$ and a $\game \Sigma^0_3$-set
are (ordinarily){\tmem{ recursively isomorphic}}. A further analogy with
$\omega_1^{\tmop{ck}}$ and $\Sigma^0_1$-Determinacy emerges: just as
strategies for  Player $\mathit{I}$ in $\Sigma^0_1$-games appear all the way
up to stage $L_{\omega_1^{\tmop{ck}}}$ in the $L$-hierarchy, so strategies for
Player $\mathit{I}$ in $\Sigma^0_3$-games appear all the way up to stage
$L_{\alpha_0^{\mathsf{E}}}$, thus giving a second characterisation of what
$\alpha_0^{\mathsf{E}}$ ``is''. It is thus the analogue of $\omega_1^{ck}$, and, as mentioned above, is the first ordinal not generalized ittm recursive (Lemma \ref{leastnonrecursive}).\\

\nod{\sc Theorem \ref{Thm1.2}}
 {\em If $A$ is a $\Sigma^0_3 (x)$ set so that the game $G (A)$ is won by Player
  $\mathit{I}$, then there is a generalised-ittm-recursively computable (in
  $x$) winning strategy $\sigma$ for $\mathit{I}$. That is, for some index $e$
  dependent on the definition of $A$, but not $x$, $\{ e \}^{\mathsf{E}} (x)
  \da \sigma$.}\\

We can then extend Theorem \ref{Psicomplete} to generalise the classical results that:

\begin{theorem}[Kleene, et al.]\label{ThmKleene}
Let $H^K$ be the halting set for Kleene's type 2 recursion. Let $G_1$ be a complete $\Game \S^0_1$ set of integers, and $\Psi_0$ the $\S_1$-$Th(L_{\w_1^{ck}})$, \ie, the set of $\S_1$-sentences true in $L_{\w_1^{ck}}$.
Then:
$$\Psi_0 \equiv_1 H^K\equiv_1 G_1.$$
\end{theorem}
\nod{\sc Theorems \ref{H=Psi}, \ref{Psicomplete}}\label{Thm1.3}
{\em Let $\Psi= \S_1$-$Th(L_{\a_0^\me})$. Then:
$$\Psi \equiv_1 H\equiv_1 G_3.$$
}

From these ideas further classical style results can be developed: we prove that the ittm semi-recursive in a functional $\mi$ relations form a Spector class; we prove also a Stage Comparison Theorem  (\ref{OrdComp}) for  computations generalised ittm recursive for a functional $\mi$, and thence a Gandy Selection Theorem \ref{GandySelection}. This allows for the usual closure and regularity properties for sets semi-recursive in $\mi$, some of which we state without proofs, as they follow very closely those in \cite{Mosch4} or \cite{M}.

\

\nod{\sc Theorem 4.13 [Spector-Gandy Theorem]}\label{SpectorGandy} {\em The following are equivalent for an $A\sset
\omega$:\\
(i) $A$ is semi-recursive in $\mi$;\\
(ii) There exists a $\Sigma_1$ $\vp(v_0)\in {\call}_{\dot I}$ so that
 $$m\in A \equi L_{\alpha^\mi_0}[\mi]\models\vp[m,\mi];$$
(iii) There exists a $P$ recursive in $\mi$ so that 
$$m\in A \equi \ex y \mbox{ recursive in } \mi \, (P(m,y)).$$
}

The pages that follow are intended to be an introduction to this generalised
type-2 ittm recursion; this is outlined in Sections 1-5.  In Section 2 the regular infinite time turing machine  is introduced. This goes back to the original definition in \cite{HL}. An important example is the {\em Theory Machine} (TM - Section 2.2) that is used to compute codes for levels of the constructible hierarchy and their $\Sigma_2$-theories. 
Generalized type 2 ittm recursion in a functional $\mf :{\bai} \imp \omega$ is introduced in Section 4.1. In Section 6 we prove the characterisation of the halting problem for generalized recursions in $\me$, in terms of winning strategies for $\Sigma^0_3$ (or  $G_{\delta\sigma}$) games).  Some familiarity with the notions of Gale-Stewart perfect information games on
integers will be assumed here (\cf \cite{Mosch4}); for the argument (in 6.2) that infinite nestings (\cf Definition \ref{Def3.11}) imply $\Sigma^0_3$-Determinacy some familiarity with the arguments of \cite{W2011} will be helpful.\\

{\em Historical and Related Remarks} {\em Infinite Nestings} were first defined in \cite{W2011}. These were at the $\Sigma_2$ level,  and  then were generalised to $\Sigma_n$-nestings and both  extensively used in \cite{AgWe2026}. The $\Sigma^0_3$-determinacy argument with nestings was used as template for the proof of determinacy of $\Sigma^0_3$-long games 
played on reals in \cite{AgWe2021}. 
In \cite{We12} the result is explicitly stated that  all games on integers at this level have strategies $\Delta_2$-definable over the smallest $L_{\beta_0}$ that would support an infinite nesting.
It was  shown in  \cite{W2011} that the reals of $L_{\alpha^\me_0}$ were $\game\Sigma^0_3$.
In \cite{We12} the least non-recursive in $\me$ ordinal $\a_0^\me$ was identified by using an argument that a program searching for strategies for such games, needed to run for this length of time in order to locate them (thus establishing Theorem \ref{Thm1.2}).  Here we have identified $\a_0^\me$ by a different argument that searches directly for an infinite ascending sequence of ordinals that could be the lower ordinals in such a nesting. This is both conceptually simpler and avoids the dependency on notions not immediately connected to the nesting idea. 

In \cite{Ha18} Hachtman showed the equivalence between infinite $\Sigma_2$-nestings supported at those levels of the $L$ hierarchy where every set is countable, and those levels whose reals comprised  models of {\boldmath $\Pi^1_2$} monotone induction.

Although we do not touch on reverse mathematical themes here, these results can be couched in terms of subsystems of analysis (\cf \cite{Si99}). Using levels of restricted Comprehension $\mathsf{CA}$: that $\Pi^1_3$-$\mathsf{CA}\vdash\Pi^0_3$-$Det$ was also shown in \cite{We12}. The argument there shows that $\Pi^1_1$-$\mathsf{CA}+ Det(\Pi^0_3)$ proves the existence of $\beta$-models of $\Delta^1_2$-$\mathsf{CA}$ (namely that there are levels of the $L$-hierarchy which are $\Sigma_2$-admissible). This was generalized by Montalban-Shore \cite{MoSh12} to showing (amongst other results) that $\Pi^1_{n+2}$-$\mathsf{CA}$
$\vdash n$-$\Pi^0_3$-$Det$, where $ n$-$\Pi^0_3$ denotes the $n$-level of the difference hierarchy on $\Pi^0_3$. They then have:  $\Pi^1_1$-$\mathsf{CA}+ Det(n$-$\Pi^0_3)$ proves the existence of $\beta$-models of $\Delta^1_{n+2}$-$\mathsf{CA}$.
Section 7 gives more detail and additional recent results of this kind.

It is very reasonable to ask: if there are natural examples of complete $\Sigma^0_3$-sets of reals or of their complements, complete $\Pi^0_3$-sets. Example (1) $\{x\in \baire\mi \,\, \mathrm{Lim}_{n\rightarrow\infty}x(n)\neq \infty\}$. (2) In \cite{Hardy93} it is remarked that $G\dfs\{x\in \re\mid  \mathrm{Lim}_{n\rightarrow\infty} sin (n!\pi x) =0\}\supseteq \mathbbm{Q}$. $G$ can be seen to be a complete $\Pi^0_3$-set of reals.  (3) In \cite{Gold94} Goldstern shows that the set of real numbers that are {\em not} normal in base 10, is a complete $\Sigma^0_3$-sets of reals.
(4) The set of three quantifier $\ex\all\ex \ldots$ sentences true of arithmetic is a complete $\Sigma^0_3$ set of integers. Indeed questions involving finitude are often $\Sigma^0_3$. (5) The set of $(e,X)$ so that $P_e^X(m)$ converges for cofinitely many $m$ is a complete $\Sigma^0_3$ subset of $ \nat \times\baire$. (Here $P_e^X$ is the $e$'th Turing machine with oracle $X$.) (6) The set of $X\in \cant$ so that $X$ codes a countable graph with at most finitely many connected components.

\section{Infinite Time Turing Machines}

Any set theoretical or  computability theoretic notation here is standard. We define an  {\em admissible set} as a model of the Kripke-Platek axioms: Extensionality, Pairing, Foundation, $\Delta_0$-Collection, and  $\Delta_0$-Separation, and we shall include here the axiom of Infinity. We refer the reader to \cite{Bar} for a development of this set theory. This set of axioms we shall denote by $\mathsf{KP}$. By $\mathsf{KPI}$ we mean $\mathsf{KP}$ with the additional axiom ``$\all x\ex m( ({\mathsf{KP}})^m \wedge x\in m$." By ``$\S_n$-$\mathsf{KP}$'' (or ``{\em $\S_n$-admissibility}''), we shall mean
$\mathsf{KP}$ augmented by  $\Sigma_n$-Replacement (which in our relativised constructibility models, is sufficient to prove $\Delta_n$-Separation). $\Sigma_1$-Replacement is provable in $\mathsf{KP}$, thus $\mathsf{KP}$ is $\S_1$-$\mathsf{KP}$. We call an ordinal $\g\,\, \S_n$-{\em admissible}, if $(\S_n$-$\mathsf{KP})^{L_\g}$.

\subsection{Basic Construction}

\raisebox{-0.0016366612111292965\height}%{\includegraphics[width=15.77961432506887cm,height=5.129738947920766cm]{TestMoganTheoryMachine-1.pdf}}
{\includegraphics[width=10cm,height=3cm]{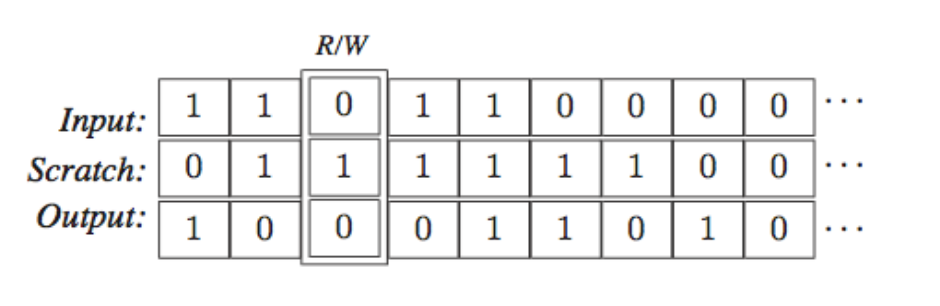}}
\

Let the (ordinary) Turing programmes be enumerated as $\pa{P_e \mid e \in \omega}$ with
`indices' $e$. We suppose the programme is listed as a finite list of 
instructions $I_0, \ldots, I_k$. \ The {\tmem{current instruction number}} at
time $\alpha + 1$, $I (\alpha + 1)$, and its R/W (``Read/Write") action {\etc}is determined by
the programme $P_e$ in the usual Turing way, given the instruction number $I(\alpha)$ and R/W head position at time $\alpha$. Of course we intend the ordinal $\alpha$ to be possibly transfinite.

At limit times $\lambda$ we use Liminf operations to decree behaviour. If the cells of the machine
are enumerated $\pa{C_i \mid i \in \omega}$ with values at time $\nu$ denoted
by $\pa{C_i (\nu) \mid i \in \omega}$ then we set at limit time $\lambda$:
$$C_i (\lambda) = 1 \Equi \ex \alpha <
\lambda \all \beta < \lambda (\alpha < \beta \imp C_i (\beta) = 1) $$
Thus, if the value of the $C_i (\alpha)$ alternates 0/1 cofinally often in
some limit ordinal $\lambda$ then \ $C_i (\lambda) = 0$. The {\tmem{current
instruction number}} $I (\lambda)$ at a limit time $\lambda$ is
$\tmop{Liminf}_{\alpha \imp \lambda} I (\alpha)$. ``OT'' abbreviates ``output
tape''. The {\em current read/write  head position} at limit time $\lambda$, $R(\l)$, we define as   $\tmop{Liminf}^\ast_{\alpha \imp \lambda} R (\alpha)$. $\tmop{Liminf}^\ast$ here has the same value as $\tmop{Liminf}$ except when the latter is $\w$ (when the $R/W$ head has wandered off to the end of the tape). In which case the head is positioned back on the first triad of cells $\pa{C_0,C_1,C_2}$ before continuing. 
This is slightly different to the architecture of \cite{HL}. There the machine entered a special ``limit state" $q_L$ at limit times, and the $R/W$ was always set back to read the first triad. One can show this makes no difference to the class of computable functions.  However, what it does do if we regard the transition table/program as written out logically, is that the instruction number, $I(\l)$, is then at the head of the outermost program loop, so subroutine, that the machine entered cofinally often before $\lambda$, and the $R/W$ head is placed at that entrypoint.
\begin{definition}
  We say $P_e (x)$ is {\tmem{convergent}} and write
  $$P_e (x) \da \mbox{ \tmem{ if $\ex \tau
  ($ the OT of $P_e (x) $ remains constant after time } $\tau$)}.$$
  
  If the contents of that OT is the real $y$ we write $P_e (x) \da y$. \ If
  we wish to indicate the time $\tau$ at which this convergence first occurs, we write $P_e (x)
  \da^{\tau} y$. It is {\tmem{divergent}} otherwise, and we write $P_e (x)
  \ua$.
\end{definition}

A synonym for convergence is being ``eventually settled''. The machine may not
formerly halt, but it makes no further changes to the OT. We treat formal
machine halting as a special kind of convergence.

\begin{definition}
  The {\tmem{snapshot at time or stage $\alpha$ is}} $s_{\alpha} = s_{\alpha}
  (e, x)$ {\tmem{of}} $P_e (x)$ is a sequence consisting of the current
  instruction number $I (\alpha)$ about to be performed, an integer $R
  (\alpha)$ representing the position of the R/W head, and an
  $\omega$-sequence of $0 / 1$ representing the cell contents at time $\alpha
  : \pa{C_k (\alpha) | k < \omega \nobracket}$.
\end{definition}

We note the following:

\

{\nod}{\tmstrong{Observation}}: a course of computation of some $P_e (x)$
with input $x \in \cant$ is absolute to $L [x]$ and can be defined by a
$\Delta_1$-recursion so that for $\Lim(\a)$  $\pa{s_{\beta} \mid \beta < \alpha} \in
\Delta^{L_{\alpha} [x \}}_1$ with the snapshot $s_{\alpha}$ at stage $\alpha$
being $\Sigma^{L_{\alpha} [x]}_2$.

\begin{definition}
  \label{Def3}An ordinal $\xi$ is called {\tmem{$\Sigma_2$-extendible}} if
  there is $\delta > \xi$ with $L_{\xi} \prec_{\Sigma_2} L_{\delta} .$ We call
  such a $(\xi, \delta)$ {\tmem{a ($\Sigma_2$-)extendible pair}}. \ The least
  such pair is usually written $(\zeta, \Sigma)$. 
\end{definition}

\begin{lemma} (i) A $\Sigma_2$-extendible ordinal is always $\Sigma_2$-admissible,
and is a limit of such. (ii) If $L_\xi$ is $\Sigma_2$-extendible to $L_\delta$ is also a limit of $\Sigma_2$-admissibles, but need not be
itself $(\Sigma_1 \nobracket$-)admissible. (iii) For the least pair $(\zeta, \Sigma \nobracket$), $L_{\Sigma}$ is the
unique such end extension. 
\end{lemma}
\pf
 (i) Argue that $\Sigma_2$-Collection
  holding in $L_{\xi}$ follows from its $\Sigma_2$-extendibility.
  (ii) Argue that there are arbitrarily large $\Sigma_2$-admissibles below $\xi$, and so the same must be true below $\delta$.(iii) To show
  that, {\eg}, $\nocomma \Sigma$ is not admissible, consider $T^2_{\alpha}
  \dfs \Sigma_2$-$\tmop{Th} (L_{\alpha}) \in L_{\alpha + 1}$. By
  extendibility, $T \dfs T^2_{\zeta} = T^2_{\Sigma}$. Let $\alpha_0 = \zeta$
  and $\alpha_n > \alpha_{n - 1}$ be least with $T \cap n = T^2_{\alpha_n}
  \cap n$ (where we think of sentences as recursively coded by integers).\qed\\

{\tmstrong{{\nod}Remark}}: The notion easily relativises: we say that $\xi$ is
$\tmmathbf{x}$-{\tmem{extendible}} if there is $\delta > \xi$ with $L_{\xi}
[\tmmathbf{x}] \prec_{\Sigma_2} L_{\delta} [\tmmathbf{x}]$ and so forth.

\

By the Observation we have:

\begin{lemma}
  $P_e (n)$ has identical snapshots at times $\zeta$, $\Sigma$: $s_{\zeta} =
  s_{\Sigma}$.
\end{lemma}

Having identical snapshots is close to being forever looping.

\begin{definition}
  \label{3.1} \ We shall say that a computation such as $P_e (x)$ `exhibits
  final looping behaviour'  (`at stage $\sigma$', or `by stage $\tau$'), if
  there are stages or times $\xi < \sigma (\leq \tau$) with  (a)  identical
  snapshots at $\xi$ and $\sigma$, and moreover  (b)  no cell that had a
  stable value at time $\xi$ changes that value in the interval $(\xi,
  \sigma)$. We say the computation `has entered a (final) repeating loop' at time or
  stage $\alpha$ if $\xi \leq \alpha$ for such a pair $(\xi, \sigma)$.
\end{definition}

\subsubsection{The  Jensen $J$-hierarchy}
Some of our results, and means of argumentation, rely on analysing the behaviours of various levels of $L$ and how that influences machine behaviour. In particular we shall want to consider taking elementary hulls, usually at the level of $\S_1$, or $\S_2$ elementarity. The discussion here, and in particular for the universal ``theory machine" (``TM") to come, is much facilitated by the use of these hulls and $\Sigma_2$ skolem functions. As is well known the G\"odel levels $L_\a$ are not friendly for this kind of reasoning: they are not closed under ordered pairs for example. Jensen's alternative
hierarchy $J_\a$ for $\a \in On$ provides a rearrangement of the sets of $L$ into a hierarchy of transitive models of a weak set theory, which are all rudimentarily closed. We shall define this hierarchy, and show that we can define {\em uniform $\S _2$-Skolem functions} for all the levels in the region of interest for our ittm theory.

The reader who is unfamiliar with or unwilling to study too closely the details of this hierarchy, can simply skip these sections and treat all talk of $J_\a$'s as being about $L_\a$'s with very little loss of understanding as to the effect on the ittm discussion. Then they are asked to take the existence of such skolem functions on trust.

\begin{definition}[The rudimentary functions]\label{rudimentaryfunctions}
  
  A function $f:V^n\imp V$ is {\em rudimentary} if it can be finitely generated
by the following schemata:\\ $
(i) f(\vec x) = x_i\\
(ii)  f(\vec x) = x_i\back x_j\\
(iii) f(\vec x) = \{x_i,x_j\}\\
(iv) f(\vec x)= h(g(\vec x))\\
(v) f(y,x)=\bigcup_{z\in y}g(z,x).
$
\end{definition}
These schemata are exactly those for the primitive recursive set functions, minus the recursion scheme itself.
We shall sometimes have recourse  to the {\em ``rudimentary-in-A'' } functions, where we have additionally:\\
$
(vi) f(\vec x) = x_0\cap A$ for $A\in V$ a set.

The rudimentary functions were invented independently by Gandy \cite{Ga74} and Jensen \cite{JeKa71}. Schemata (iv) and (v) essentially involve infinitely many rudimentary functions, but it is possible to show that there is a finite set of {\em basic rudimentary functions} from which all other rudimentary functions can be built up. We shall not need that result here, so we do not repeat this analysis (see for example \cite{De} or \cite{Je72}).

\begin{definition}(i) The {\em rudimentary closure} of a set $U$ is the closure of $U$ under
the rudimentary functions.
The {\em rudimentary closure in $A$} of a set $U$ is the closure of $U$ under
the rudimentary  functions.

(ii) For $U$ the set $rud(U)$ is defined as the rud.~closure of $U\cup\{U\}$; 
$rud_{A_0}(U)$ is defined as the rud.($A_0$)~closure of $U\cup\{U\}$.
\end{definition}

\begin{lemma} If $U$ is transitive, so are $rud(U)$, $rud_{A_0}(U)$.
\end{lemma}

\begin{definition}
[The Jensen $J$-hierarchy]{\em \cite{Je72}}\mbox{ }\\ $$J_0 = \emp; J_{\a +1}= rud(J_\a);  Lim(\l)\rightarrow J_\lambda = \bigcup_{\a < \l} J_\a;  L=\bigcup_{\a \in On} J_\a.$$\\ The $J[A]$ -hierarchy is defined in the same way using $rud_A$ closure.
\end{definition}

Each $L_\l$ is also rud. closed for $\l$ a limit. Then we also have that $J_\a=L_\a \Equi \omega\cdot\a =\a$. So a modest amount of closure, for example primitive set recursive closure, ensures that the levels are identical. The following are taken from \cite{Je72} once more.

\begin{lemma}
For any $\Delta_0$ $\varphi
  (v_0, \ldots v_n)$ there is a rudimentary function $F_{\varphi}$
  % (thus acombination of basis functions) 
  so that $$\varphi [x_0, \ldots x_n] \Equi
  F_{\varphi} (x_0, \ldots x_n) \neq 0.$$

\end{lemma}

The theory of rudimentary functions give us $\Delta_0$-Separation.
\begin{lemma}
  \label{rud=comp}For any $\Delta_{0}\,    \varphi (v_0, \ldots v_n)$ there is a rudimentary
  function $F_{}$ so that
 $$F (a,v_{0} , \ldots ,v_{k-1} ,v_{k+1} ,
  \ldots v_{n} ) =  \{v_{k} \in a \mid \varphi (v_{0} , \ldots v_{n} )\}.$$
\end{lemma}

We sketch a presentation of the ``Theory Machine'' of {\cite{FrWe08}}. We use  the $J_{\alpha}$-hierarchy  presented above. Those not so familiar
with this can pretty much read ``$L_{\alpha}$'' for ``$J_{\alpha}$''
throughout with little distortion of the truth; recall that in any case if
$\alpha = \omega \cdot \alpha$ then $L_{\alpha} = J_{\alpha}$. For example in the next definition if $\b >0$ then $J_\b$ will be admissible, and thus $J_\b=L_\b$. 

\begin{definition}
  Let $S^n_{\alpha} \dfs \{ \beta < \alpha \mid J_{\beta} \prec_{\Sigma_n}
  J_{\alpha} \} \cup \{ 0 \}$.
\end{definition}

\begin{note}
  (i)  For $0< \beta \in S^1_{\alpha}$ we shall say  ``$\beta$ is $\Sigma_1$-stable in
  $\alpha$''.
  
  (ii) \ If $S^1_{\alpha}$ is unbounded in $\alpha$ we say that $\alpha$ is
  {\em non-$\Sigma_1$-projectible}, and in fine structural terminology
  ``$\rho^1_{\alpha} = \alpha$''. One may reason that $S^1_{\alpha} \in
  \Pi_1^{J_{\alpha}}$. Further $\power(\omega)\cap L_\alpha$ in such a case is a model of $\Pi^1_2$-$\mathsf{CA}_0$.
  
  (iii)  We use the observation that if $\sigma$ is a $\Sigma_2$-sentence,
  then for $\beta \in S^1_{\alpha}$, that $J_{\beta} \models \sigma \Imp
  J_{\alpha} \models \sigma$ (since if $\sigma \equiv \ex u \psi (u)$ with
  $\psi \in \Pi_1$, if $u_0 \in J_{\beta} \models \psi [u_0]$ then by upwards
  persistence of $\Pi_1$ sentences $J_{\alpha} \models \psi [u_0]$) .
  
  (iv)  We recall the existence of a uniform (in $\nobracket \alpha)$
  $\Sigma_1$-definable skolem function $h^1$, so that when $h^1_{\alpha}$ is
  this function defined in $J_{\alpha}$ then for any $A \sset J_{\alpha}$,
  then $h^1_{\alpha}$``$\omega \times \left( A \cup \left\{ \emp \right\}
  \right)$ is the least $\Sigma_1$-skolem hull containing $A$ in $J_{\alpha}$.
  We write often $h_{\alpha}$ for $h^1_{\alpha}$. Thus $h_{\alpha}$``$\omega
  \times A \prec_{\Sigma_1} J_{\alpha}$. In particular if $\alpha < \alpha'$
  then $h_{\alpha} \sset h_{\alpha'}$.
  
  (v) \ We also have that if $\beta = \max S^1_{\alpha}$ then  (a)
  $h^1_{\alpha}$``$\omega \times \beta + 1 = J_{\alpha}$; and  (b)  if
  $J_{\alpha} \models$``$\beta$ {\em  is countable}'' then $h^1_{\alpha}$``$\omega
  \times \{ \beta \} = J_{\alpha}$; (c)  if
  $J_{\alpha} \models$``{\em every ordinal is countable}'' then if $X\prec_1 J_\alpha$ then $X=J_\beta$ for some $\b\leq\a$. \ \ ((b) follows from (a), since then
  there is an onto function $f \in J_{\alpha} \nocomma$, $f : \omega \imp
  \beta$, with $f$ $\Sigma_1$ definable in $J_{\alpha}$ from the parameter
  $\beta$; {\ie} $f \in h^1_{\alpha}$``$\omega \times \{ \beta \}$); (c) follows from (b) and the observation that $h^1_{\alpha}$``$\omega
  \times  X = h^1_{\alpha}$``$\omega
  \times On \cap X = J_{\b}$ where $\b=On\cap X$.
\end{note}
These facts will be used without further remark in the sequel. One concept that does require further remark is that we shall use, \eg in the definition of the Theory Machine,   {\em uniform $\Sigma_2$-Skolem functions,} for $\alpha$ such that $J_\alpha
\models$``{\em every set is countable}''. We say that $J_\a$ is ``locally countable". In particular such exist for $\a\leq \b$ where $L_\b$ is the least transitive model of $ZF^-$. The whole discussion of this paper takes place below this ordinal $\beta$ - and so could be formalised within analysis.
We say that a $J_\a$, (or $L_\a$) is {\em locally countable} if the structure $\pa{J_\a,\in}$ is a model of {\em ``every set is countable"}. In this paper all levels $J_\a, L_\a$ we deal with have $\a\leq\b$ and so will be locally countable.
In this region we shall show that we have uniform $\Sigma_n$-skolem functions. (This however notably fails in general for those $L_\alpha$ or $J_\a$ which do not think every set is countable.)

\subsubsection{Uniform $\S_{n+1}$-skolem functions}
We shall give a proof of the existence of these uniform $\S_{n+1}$-skolem functions at this level. Familiarity with these arguments is not necessary for the rest of the paper, and  the reader may simply wish to skip this subsection, and take the existence of uniform $\Sigma_2$-skolem functions on trust. 

The proof hinges on the fact that for locally countable $J_\a, L_\a$ every set $x \in J_\a$ say is $\Sigma_n$ definable in $J_\a$ from $0$ or from some singleton $\{\b\}$ where $\b \in S^n_\a$. In this paper we shall only require $\S_2$-skolem functions, so we give the argument in detail just for that case.
 This version of the argument has its roots in \cite{Fr08?}. Let $\pa{\psi_i}_{i<\omega}$ be a recursive enumeration of all $\Sigma_2$ formulae of the form $\ex v_0\vp_i(v_0,v_1,v_2)$ where in turn $\pa{\vp_i}_{i<\omega}$ is the  recursive enumeration of all $\Pi_1$ formulae of the form $\vp(v_0,v_1,v_2)$ with free variables as shown.
 Then we may think of any $\S_2$ relation $R(x,y)=R_i(x,y)$ as being given by some formula of the form $\ex u\vp_i(u,x,y)$. 
 
 We fix once and for all a standard Turing recursive pairing function $\pi:\w\times\w \equi \w$, with unpairing functions $(-)_0$ and $(-)_1$ so that $\pi^{-1}(u) = \pa{(u)_0,(u)_1}$. 
  
\begin{lemma} There is  a parameter free $\S_2$-definable $\S_2$-skolem function, $h^2(v_0,v_1)$ so that for any $J_\alpha$ which is a model of {\em ``Every set is countable"}, if $R(x,y)$ is any $\Sigma_2^{J_\a}$ definable relation then

$$  \all y[\ex x R(x,y) \imp R(h^2(i,y),y)]$$

where $R = R_i=\{\pa{x,y}\mid \ex u\vp_i(u,x,y)\}$.

\end{lemma}
\pf Let $h^1(v_0,v_1)$ be the uniformly definable $\Sigma_1$-skolem function, $\Sigma_1$-definable over any $J_\b$. (We write $h^1_\g$ to indicate the function so defined over $J_\g$.) Fix a $J_\a$ which is ``locally countable" in the above sense. \\

\nod{\em Claim 1} $\all u\in J_\a\ex \b\in S^1
_\a \ex n\in \w ( u = h^1_\a(n,\b))$.

\pf Case 1: $\ex \b < \a (\b = \max(S^1_\a))$. Then this follows from Note (v) above.\\
\indent Case 2: Otherwise. Let $\b$ be least in $S^1_\a$ with $\r_L(u)<\b$. Then let $\b_0= \max(S^1_\a\cap\b)=\max(S^1_\b)$. By Note (v) then, $h^1_\b$``$\w\times\{\b_0\}=h^1_\b$``$\w\times(\b_0+1)=J_\b$. So choose $n$ so that $(u=h^1(n,\b_0))^{J_\b}$. By upwards persistence of $\S_1$ formulae $(u=h^1(n,\b_0))^{J_\a}$.  \qed\, {\em Claim 1}\\

\nod Define:\\
$H(i,y)\simeq h^1_\a(m,\pa{\b,y})$ {\em where } $(\b,m)(=\b(i,y),m(i,y))$ {\em is lexicographically least so that }
$$(\b\in S^1_\a \vee \b=0) \wedge \ex u[u= h^1(m,\pa{\b,y})\wedge \vp_i((u)_0, (u)_1,y)]$$ {\em holds.} \\
 
 \nod{\em Claim 2} $H$ {\em  is $\S_2^{J_\a}$-definable.}\\
\pf The first conjunct here is $\Pi_1$ - this is because $S^1_\a$ is itself $\Pi_1^{J_\a}$.
The matrix of the second conjunct is itself a conjunction of, first, a $\S_1$ and then a $\Pi_1$ statement. Thence it is overall $\Sigma_2$.
 To say that $(\b,m)$ is lexicographically least for this to hold, is to say in addition to this that:\\

$\all\b'<\b\all m'<m[\ex u(u=h^1(m',\pa{\b',y})\imp \neg\vp_i((u)_0, (u)_1,y)] \vee$\\
$\vee (\b=\b' \wedge \all m'<m [\ex u(u=h^1(m',\pa{\b,y})\imp \neg\vp_i((u)_0, (u)_1,y)]$.\\

However the first disjunct here can be written:\\
\mbox{}\quad\quad$ J_\b\models$``$\all\b'\all m'<m[\ex u(u=h^1(m',\pa{\b',y})\imp \neg\vp_i((u)_0, (u)_1,y)]$''\\
because of the absoluteness of $\S_1$ formulae between $J_\b$ and $J_\a$. This makes the disjunct $\Delta_1$ (in $\b$). The second disjunct is some bounded natural number quantification in front of a $\S_2$ predicate, and so is again $\S_2$.\\ \mbox{ } \qed \,{\em Claim 2}\\

 We thus have:\\
\nod{\em Claim 3} $H(j,y) \simeq v$ iff $\vp_j((v)_0,(v)_1,y)\wedge v = h_1(m(j,y)),\pa{\b(j,y),y}) $.\\

Suppose $R_j(x,y) $ is a $\S_2$ relation with $R_j (x,y)\Equi \ex u\vp_j (u,x,y)$.
Now define \\

$
\begin{array}{rcl}
\hfill h^2(j,y)=x & \equi_{df} & \ex u H(j,y)=\pa{u,x}\\
               & \equi\hspace{1.0em} & \ex u \vp_j(u,x,y)\\
               &\equi\hspace{1.0em} & R_j(x,y)\\
               &\equi\hspace{1.0em} & R_j(h^2(j,y),y).

\end{array}
$

Thus $h^2$ is a uniform $\S_2$-skolem function as desired. \qed\\

By induction on $n$, for locally countably $J_\a$ we have uniform $\S_{n+1}$-skolem functions, just using the above argument to define $h^{n+1}(j,y)$ (substituting $h^n$ for $h^1$ to define $H$).
\begin{lemma}
For every $n\geq 0$ there is  a parameter free $\S_{n+1}$-definable $\S_{n+1}$-skolem function, $h^{n+1}(v_0,v_1)$ so that for any $J_\alpha$ which is a model of {\em ``Every set is countable"}, if $R(x,y)$ is any $\Sigma_{n+1}^{J_\a}$ definable relation then
$$  \all y[\ex x R(x,y) \imp R(h^{n+1}(i,y),y)]$$
\nod where $R=R_i=\{\pa{x,y}\mid \ex u\vp_i(u,x,y)\}$ with $\pa{\vp_i(v_0,v_1,v_2)}_{i<\w}$ a recursive enumeration of $\Pi_n$ the three free variable formulae with free variables indicated.
\end{lemma}
\pf Formally an induction on $n$, substituting ``$n+1$'' for ``2" and elsewhere {\em m.m.} Otherwise as above. \qed\\

We thus have:
\begin{lemma}
Let $J_\a$ be locally countable. Let $h^n$ be as above. Let $X\sset J_\a$. Then $\tilde X\dfs h^n$``$\w\times \mbox{}^{<\w}X$ is the least $\S_n$ skolem hull of $X$ in $J_\a$. Moreover $\mathop{Trans}(\tilde X)$ and $\tilde X=J_{\bar\a}$ for some $\bar\a\leq\a$.
\end{lemma}
Standard arguments show:
\begin{lemma}\label{transitivehull}
For any pair $\z\neq\S$ with $J_\z\prec_{\S_2} J_\S$ satisfy: 
$J_\z$ is $\S_2$-admissible, and $J_\S$ is a limit of $\S_2$-admissibles. Then for both ordinals we have $L_\z=J_\z$ and $L_\S=J_\S$.
Moreover
$h^2_\z$``$\w\times \{\z\} = L_\S.$ 
\end{lemma}
We shall only be using the last for $n=1,2$.
\begin{lemma}
The lexicographically least pair $\z\neq\s$ with (a) $J_\z\prec_{\S_2} J_\s$ is exactly the lexicographically least pair so that (b) $T^2_\z=T^2_\s$. Further this implies that (i)  $J_\s\models$ `` $V=HC$'' and (ii) $h^2_\z$``$\w =h^2_\s$``$\w = J_\z$; that is both of these  hulls are transitive and are the same $J_\z$. 
\end{lemma}
\pf Suppose  $\z,\s$ are lexicographically least satisfying (a) with $\z\neq\s$.  Note  $J_\z\models$ `` $V=HC$''. 
Since otherwise, if $\eta <\z$ is such that $\eta = (\omega_1)^{J_\z}$ then the skolem hull $X= h^2_\eta$``$\w$ satisfies $X\prec_{\S_2}J_\eta$.   But then by Lemma \ref{transitivehull} $X=J_\a$ for some $\a<\eta$, and this pair $(\a,\eta)$ contradicts our choice of $\z,\s$.  
Then by $\S_2$ elementarity we have $T^2_\z=T^2_\s$. 

Suppose the least pair $\a\neq\b$ with $T^2_{\a}=T^2_{\b}$ was lexicographically earlier than $(\z,\s)$ satisfying (a).  However all the statements of the form: ``$\ex y(y=h^2(i))$'' (for $i\in \w$); ``$\ex i,j<\w( h^2(i)\in\!\!\back\!\!=h^2(j))$''  are all $\Sigma_2$ sentences and so are in both $T^2_\a $ and $T^2_\b$ or in neither. But then $X_0\dfs h^2_\a$``$\w \cong X_1\dfs h^2_\b$``$\w$. 
Then the transitive collapses of $X_0,X_1$ are then the same $J_{\bar\a}$. Then $\bar \a =\a$ since otherwise $T^2_{\bar\a}=T^2_{\b}$, contradicting the leastness of $(\a,\b)$. But this means $J_\a\prec_{\S_2}J_\b$, contradicting the choice of $(\z,\s)$. 

By the same argument as just given for $\a,\b$ we have in (ii)  identical hulls  equalling $J_\zeta$.\qed

\subsection{The Theory Machine $\tmop{TM}$}

%%%%%%%%%%
\begin{definition}
  (i)  Let $T^n_{\alpha} \dfs \left\{ \ul \sigma \ur \in \omega \mid \right.
  \sigma$ is $\Sigma_2 \nobracket \wedge J_{\alpha} \models \sigma \}$ ; \
  
  (ii) $ \widehat{T_{\alpha}} \dfs \left\{ \ul \sigma \ur \in \omega \mid
  \right.$ \ $\sigma$ is $\Sigma_2$ $\wedge$ $\left.  \ex \beta < \alpha \all
  \tau \in (\beta, \alpha)( \sigma \in T^2_{\tau} )\right\}$.
\end{definition}

Then $\widehat{T_{\alpha}}$ is the set of $\Sigma_2$-sentences that are
`eventually true' below $\alpha$. Two preliminary lemmata are needed before
describing the ittm theory machine.

\begin{lemma}
  There is an  (ordinary)  Turing recursive function $f : \omega \times \omega
  \imp \omega$, so given by an index $e$, so that for any $\tmop{Lim}
  (\lambda)$ \ satisfying $J_{\lambda} \models$``Every set $x$ is countable'',
  if we set $T = \widehat{T_{\lambda}}$ \ then $T^2_{\lambda}$ is uniformly
  r.e. in $T$, {\via}$f$, that is \ \ \ \ \ \ $\ul \sigma \ur \in
  T^2_{\lambda} \equi \ex i f \left( i, \ul \sigma \ur \right) \in T$. \ \ 
\end{lemma}

{\pf} Let $\sigma \equiv \ex u \psi (u)$ be a sentence with $\psi \in \Pi_1$. We

{\tmem{Claim: $\sigma \in T^2_{\lambda} \equi \ex i \left[ \ex \tau_0 \all
\tau \in (\tau_0, \lambda) \right.$}}

\ \ \ \ \ \ \ \ \ \ \ \ \ \ \ \ \ \ {\tmem{$J_{\tau} \models$}}``$\ex \beta
\in S^1_{\tau} \left( (\beta \neq 0 \wedge \sigma^{J_{\beta}})  \, \vee \left(
h_{\tau} (i, \beta) \da \wedge \psi [h_{\tau} (i, \beta)]^{J_{\tau}} \right)
\right)$''$\nobracket]$

\

Note first that the expression in quotation marks on the right hand side,
$\eta_{\sigma} (i)$ say, here is, if true, a member of $T^2_{\tau}$, being
$\Sigma^{J_{\tau}}_2$ in $i$. We thus shall have $\sigma \in T^2_{\lambda} \lr
\ex i \ul \eta_{\sigma} (i) \ur \in T$ and the Lemma is proven.

\

{\pf}of Claim.

Case 1. $S^1_{\lambda}$ is unbounded in $\lambda$.

Suppose the left hand side holds of $\sigma$. Suppose $\psi
(u_0)^{J_{\lambda}}$ holds for $u_0$.Then for some sufficiently large $\beta
\in S^1_{\lambda}$, $u_0 \in J_{\beta}$, and then $\psi (u_0)^{J_{\beta}}$.
But $\beta \in S^1_{\lambda} \imp \beta \in S^1_{\tau}$ \ for any $\tau >
\beta$; consequently the first disjunct of the right hand side holds. For the
converse direction, fix the given $i$. By the Case hypothesis we can assume
that $\tau$ itself is in $S^1_{\lambda}$. But then if the first disjunct
holds, if $\sigma^{J_{\beta}}$ and $\beta \in S^1_{\tau}$ then
\tmtextbf{$\beta \in S^1_{\lambda}$} and thence $\sigma^{J_{\lambda}}$. If the
second disjunct holds for the supposed $i$ $\psi [h_{\tau} (i,
\beta)]^{J_{\lambda}}$ holds for the same reasons.

Case 2 $\beta_0 \dfs \max S^1_{\lambda} < \lambda$ exists.

By the bullet points above every $x \in J_{\lambda}$ is of the form
$h_{\lambda} (i, \beta_0)$. Again suppose the left hand side holds of $\sigma$
and $\psi (u_0)^{J_{\lambda}}$ holds for $u_0$. In particular now $u_0 =
h_{\lambda} (i, \beta_0)$ for some $i$. Let $\tau_0 \geq \beta_0 $ be
sufficiently large so that $\left( h_{\tau_0} (i, \beta) \da
\right)^{J_{\tau_0}}$ and thence, by the fact of $\psi$being $\Pi_1$, $(\psi
[h_{\tau_0} (i, \beta)])^{J_{\tau_0}}$. By the upwards persistence of
$\Sigma_1$ formulae in the first case and downwards persistence of $\psi$ in
the second case, these will hold in all larger $J_{\tau}$ for $\tau \leq
\lambda$ replacing $\tau_0$. But now the second disjunct of the right hand
side holds.

Conversely suppose the right hand side holds. Let $i$ be as supposed. By the
maximality of $\beta_0$ for unboundedly many $\tau' \in (\beta_0, \lambda)$
some new $\Sigma_1$-sentence about $\beta_0$ becomes true first in $J_{\tau' +
1}$. Pick such a $\tau = \tau' + 1$ of this form. Such a $\tau$ ensures that
$S^1_{\tau} = S^1_{\lambda}$ and thence $\max S^1_{\tau} = \beta_0$ too. So
suppose the first disjunct holds for such a successor $\tau$. Then if
$\sigma^{J_{\beta}}$ holds for a $0 \neq \beta \in S^1_{\tau} = S_{\lambda}^1$
we shall have $\sigma^{J_{\lambda}}$ and we are done. Thus we now suppose the
first disjunct fails for $\tau$ of this form; pick any such $\tau$, then the
second disjunct holds as witnessed by a $\beta \in S^1_{\tau}$.

Then if $\beta < \beta_0$ then $\left( \ex y (y = h_{\tau} (i, \beta))
\right)^{J_{\tau}}$ implies $\left( \ex y (y = h_{\beta_0} (i, \beta))
\right)^{J_{\beta_0}}$ by the uniformity of the definition of the
$\Sigma_1$-skolem function $h$, and the fact of $\beta_0 \in S^1_{\tau}$. But
then $( h_{\beta_0} (i, \beta) \da \wedge$ $ \psi [h_{\beta_0} (i,
\beta)]^{J_{\beta_0}} )$. But this entails that the first conjunct holds
for $\tau$, which we are assuming does not happen. Hence we must have $\beta =
\beta_0$. \ However then $\psi [h_{\tau} (i, \beta_0)]^{J_{\tau}}$ for any
$\tau$ of this form, and so for such $\tau$ arbitrarily large below $\lambda$.
By the upwards persistence of $h_{\gamma} (i, \beta_0)$ for $\gamma \in [\tau,
\lambda]$ we have \ $\psi [h_{\tau} (i, \beta_0)]^{J_{\lambda}}$ and hence
$\sigma^{J_{\lambda}}$. {\qed}

\begin{lemma}
  Let $\zeta$ be the least $\Sigma_2$-extendible ordinal and let $\Sigma$ be
  its extension: $L_{\zeta} \prec_{\Sigma_2} L_{\Sigma}$.
  
  (i)  There is a uniform procedure ittm-recursive in $T^2_{\alpha}$ in
  $\omega$ steps, for any $\alpha < \Sigma \nocomma \nocomma$, for computing a
  real $x_{\alpha}$ which is a code for the structure $\pa{J_{\alpha}, \in} .$
  
  (ii)  There is thus a partial onto function $f \in \Sigma^{J_{\alpha}}_2$,
  $f : \omega \twoheadrightarrow J_{\alpha}$.
\end{lemma}

(``Uniform'' here means the procedure is independent of $\alpha$.)  See \
{\cite{FrWe08}}.

{\pf}(Sketch) (i) \ It suffices to note that there is a uniform
$\Sigma_2$-skolem function $h = h^2_{\alpha}$ with domain a subset of $\omega
\times \omega$ which is onto $J_{\alpha}$ for those $\alpha < \Sigma$ (it is
not uniform for all $\alpha \in \tmop{On}$) . \ Granting this we can define
$\pa{i, n} \backsim \pa{j, m}$ iff $h (i, n) = h (j, m)$ and $\pa{i, n} E
\pa{j, m}$ iff $h (i, n) \in h (j, m)$. Both $\sim$ and $E$ are (ordinary)
recursive in $T^2_{\alpha}$. Then
$\pa{[{\pa{i,n}}]_{{\backsim}},E} 
\simeq \pa{J_{\alpha},
\in}$.  (ii)  For $i, n$ let $\pi : \omega \times \omega \equi \omega$ be a
recursive pairing function, and then set $f (\pi (i, n)) = h (i, n)$.{\qed}

\begin{lemma}
  There is an ittm programme $P_e = T M$ which does not converge, but
  continuously produces alternately codes $x_{\alpha} $ for levels
  $J_{\alpha}$ and their $\Sigma_2$-theories $T^2_{\alpha}$ for $\alpha <
  \Sigma$. At stage $\Sigma$ as $T^2_{\Sigma} = T^2_{\zeta}$ $T M$ loops back
  and reproduces the code $x_{\zeta}$ and continues this process thereafter
  repeating through $On$ codes and theories for $\alpha \in [\zeta, \Sigma) .$ 
\end{lemma}

{\pf} We describe the effective procedure to be formalised. The input to
$\tmop{TM}$ is presumed to be zero. We let $\pa{\varphi_n}_{n \in \omega}$ be
an effective enumeration of the sentences of $\mathcal{L_{\in}}$. We use the
$J$-hierarchy to avail ourselves of uniform $\Sigma_2$-Skolem functions. This
is not terribly important, but using the $L$-hierarchy is a bit more awkward.
Recall that if $\omega \cdot \alpha = \alpha$ then $J_{\alpha} = L_{\alpha}$.
On the output tape a theory $T$ is written with $\varphi_n \in T$ iff the
$n$'th cell contains a 1. In the first $\omega^2 + \omega \cdot 2$ stages
$\tmop{TM}$ writes the code of $J_1 = L_{\omega} = \tmop{HF}$ and \ its
$\Delta_0$-diagram to two reserved tapes, and its $\Sigma_2$-theory to the
output tape  (ot) .  (It takes less than this, but it keeps the induction
bookkeeping straight.)  We assume inductively that at time $\omega^2 \cdot
\alpha + \omega \cdot 2$ the ot contains the theory $T^2_{\alpha}$ of
$\pa{J_{\alpha}, \in}$ and the reserved tapes again the $\Delta_0$-diagram of
$J_{\alpha}$, $d_{\alpha}$, and a code for $J_{\alpha}$. With the theory
$T_{\alpha}^2$ of $J_{\alpha}$ TM can construct a code $x_{\alpha + 1}$ for
$J_{\alpha + 1}$ in $\omega^2$ additional steps together with its
$\Delta_0$-diagram $d_{\alpha + 1}$ (see next Lemma). We are are now at stage \
$\omega^2 \cdot \alpha + \omega \cdot 2 + \omega^2$. \ In an additional
$\omega \cdot 2$ steps $T^2_{\alpha + 1}$ is calculated from $d_{\alpha + 1}$
and written to OT as follows. This will take us to stage $\omega^2 \cdot
(\alpha + 1) + \omega \cdot 2$.  (This all takes some routine work to make
clear, but essentially once we have $T^0_{\beta}$ then $T^{n + 1}_{\beta}$ is
r.e. in $T^n_{\beta}$, and so in particular $T^{n + 1}_{\beta} \leq_T
(T^n_{\beta})'$. Each jump can be written out by an ittm in $\omega$-steps (in
fact the double jump can be so written, but we can ignore that), thus
requiring $\omega \cdot 2$ steps to write out the two jumps and thus obtain
the complete theory $T^2_{\beta} .$)

Of course we do this writing simply by changing the cells one by one according
to what has appeared or disappeared passing from $T^2_{\alpha} $ to $T_{\alpha
+ 1}^2$. If $\varphi_n$ is in both theories, then the 1 in the $n$'th cell is
not changed to a 0 and then back again to a 1. \ By this method of writing, at
a limit stage $\omega^2 \cdot \lambda$ for $\tmop{Lim} (\lambda)$,
$\widehat{T_{\lambda}}$ is on the ot, and thus the true $T^2_{\lambda}$ is
r.e. in the ot, by the first lemma. Hence in $\omega$ further steps it can
then write the correct $T^2_{\lambda}$ to the ot, thus by stage $\omega^2
\cdot \lambda + \omega$, and then by the last lemma the code $x_{\lambda}$ for
$J_{\lambda}$ on the scratch tape by stage $\omega^2 \cdot \lambda + \omega +
\omega$. A code for $J_{\lambda + 1}$ and the diagram $d_{\lambda + 1}$ is
written by stage $\omega^2 \cdot (\lambda + 1)$, and $T^2_{\lambda + 1}$ by
$\omega^2 \cdot (\lambda + 1) + \omega \cdot 2$.{\qed}

%\end{lemma}
\begin{lemma}
   A code for $J_{\alpha + 1}$ can be written in
  $\omega^2$ steps from a code for $J_{\alpha}$ simultaneously with its
  $\Delta_0$-diagram. 
\end{lemma}
\pf We give a sketch.
 We use the fact that $J_{\alpha + 1} = \tmop{rud}
  (J_{\alpha})$, and that there are 16 rudimentary basis functions under whose
  closure we can generate $\tmop{rud} (J_{\alpha})$. Having the
  $\Sigma_2$-theory, means we have the graph of $h^2_{\alpha}$ and so in
  effect, a partial onto map $f : \omega \imp J_{\alpha}$ with ``$f (n) \in f
  (m)$'' {\etc} recursive in $T^2_{\alpha} .$ Assume that we have a copy of $x_{\alpha}$ written on the Evens numbered cells of  (a recursive subtape of) the scratch tape:
  $\pa{2 n, 2 m} \in E_{x_{\alpha}} \equi f (n) \in f (m)$. We. may use the Odds as a
  space to build up the rudimentary closure of $J_{\alpha} \cup \{ J_{\alpha}
  \}$, by $\omega$ many passes through the basis functions applied to whatever
  has been created so far. This creates the domain of $J_{\alpha + 1}$. On
  another reserved tape simultaneously write the $\Delta_0$-diagram of the
  sets being created. To do this we use the fact that for any $\Delta_0$ $\varphi
  (v_0, \ldots v_n)$ there is a rudimentary function $F_{\varphi}$ (thus a
  combination of basis functions) so that $\varphi [x_0, \ldots x_n]$ iff
  $F_{\varphi} (x_0, \ldots x_n) \neq 0$ . 
\qed

\begin{corollary}
  For any $e$, any real that appears at some stage on the output tape of $P_e
  (0)$ is recursive in some $T^{2 }_{\alpha}$ for an $\alpha < \Sigma$, and
  thus is in $L_{\Sigma}$. Conversely for any real $y$ of $L_{\Sigma}$ there
  is an index $e$ with $y$ appearing on the OT of $P_e (0)$ at some stage
  (which perforce must be $< \Sigma$). 
\end{corollary}

{\pf}Any course of computation of a $P_e (0)$ is absolute to $L_{\Sigma}  (=
J_{\Sigma})$, with any snapshot at time $s_{\alpha}$ being
$\Sigma_2$-definable over $L_{\alpha} $,
($\pa{s_{{\beta}}{\mid}{\beta}<{\alpha}}$ is defined by a
$\Sigma_1$-recursion over $L_{\alpha}$) . Consequently the cell contents at
stage $\alpha$ are recorded by a certain recursive subset of the theory
$T^2_{\alpha}$. Conversely if $y \in L_{\Sigma}$ then $y$ is  (ordinary) 
Turing recursive in a code $x_{\alpha}$ for some $\alpha < \Sigma$: for some
$f \in \omega,$ $y = \{ f \}^{x_{\alpha}}$. Given $f$ it is easy to amend TM
so that the results of $\{ f \}^{x_{\alpha}}$ are instead continuously written
to the OT for increasing $\alpha$.{\qed}

\

Similar reasoning yields:

\begin{corollary}
  For any $e$, \ if $P_e (0) \da y$, then $y$ is recursive in some $T^{2
  }_{\alpha}$ for an $\alpha < \zeta$, and thus is in $L_{\zeta}$. Conversely
  for any real $y$ of $L_{\zeta}$ there is an index $e$ with $y$ \ $P_e (0)
  \da^{\alpha} y $ at some stage $\alpha < \zeta$. 
\end{corollary}

It is an exercise to show:
\begin{lemma}
  If $\lambda$ is admissible, then $\widehat{T_{\lambda}} =
  T^2_{\lambda}$.
\end{lemma}

\begin{definition}
  \label{lambda}Let $\lambda < \zeta$ be least so that $L_{\lambda}
  \prec_{\Sigma_1} L_{\zeta}$.
\end{definition}

\begin{lemma}
  For any $P_e (0)$, if this computation formally halts at time $\tau$, then
  $\tau < \lambda$. Conversely there are unbounded in $\lambda$ ordinals
  $\tau$ for which there is such an $e$ with $P_e (0)$ halting at time $\tau$.
\end{lemma}

{\pf} By considering the $\Sigma_2$-recursion in $L_{\zeta}$ that yields the
snapshots $s_{\alpha}$, that $s_{\tau}$ is a halting snapshot of $P_e (0)$ is
then a $\Sigma_1$ sentence in $T^1_{\tau + 1}$ which is true in $L_{\zeta}$
but first true at $L_{\tau + 1}$. By definition of $\lambda$ then $\tau <
\lambda$. Conversely there are unbounded in $\lambda$ levels $L_{\tau}$ of the
$L$-hierarchy where a new $\Sigma_1$ sentence $\sigma_{\tau}$ becomes true. 
(There could not be a bound $\lambda' < \lambda$ for such $\tau$ as then we
should have $L_{\lambda'} \prec_{\Sigma_1} L_{\lambda}$.)  However then we
could run a program that itself runs TM and halts when it finds that
$\sigma_{\tau} \in T^2_{\tau}$. This it can do by stage $\omega^2 \times (\tau
+ \omega) \nosymbol$ $< \lambda$ for example, as $\lambda$ is p.r. closed  (in
fact admissible). {\qed}

\begin{lemma}\label{writables} If $\tau <\lambda$ then $\tau$ is ittm writable: there is  $P_e (0) $ which halts with a code for a wellordering of type $\tau$ on its OT.
\end{lemma}
\pf  Amend the  program of the last lemma  \qed\\

Essentially as corollaries to the above we have:

\begin{theorem}
  (The ``$\lambda$-$\zeta$-$\Sigma$'' Theorem) . \ Let $\zeta$ be the least
  $\Sigma_2$-extendible ordinal, with $L_{\Sigma}$ the unique $\Sigma_2$
  end-extension of $L_{\zeta}$. Let $\lambda$ be as in Def. \ref{lambda}. Then
  we have that:
  
  (i) \ $(\lambda \nobracket$, $\zeta$, $\nobracket \Sigma)$ is the
  lexicographic least triple with $L_{\lambda} \prec_{\Sigma_1} L_{\zeta}
  \prec_{\Sigma_2} L_{\Sigma}$ ;
  
  (ii) $\lambda = \sup \{ \tau \mid \tau \nobracket$ is the halting time of
  some $\nobracket P_e (0) \}$
  
  \ \ \ \ \ \ \ \ \ \ $= \sup \{ \tau \mid \tau = \| y \|  \nobracket
  \nobracket \nobracket$ is the length of some ordinal code $y$ output by a
  halting program $\nobracket P_e (0) \}$.
  
  (iii) $\zeta = \sup \{ \tau \mid \tau \nobracket$ is the convergence time of
  some $\nobracket P_e (0) \}$
  
  \ \ \ \ \ \ \ \ \ $= \sup \{ \tau \mid \tau = \| y \|  \nobracket
  \nobracket \nobracket$ where $\ex e$ $\left. P_e (0) \da y\wedge y\in WO \right\}$.
  
  (iv) $\Sigma = \sup \{ \tau \mid \tau = \| y \|  \nobracket \nobracket
  \nobracket$ where $ \ex e$ with $y$ appearing on the OT of $P_e (0)$ at some
  stage $\nobracket \tau < \Sigma \} .$
\end{theorem}

{\bu} In the literature an ordinal is ``clockable'' if it is the halting time
of some $P_e (0)$. We thus have that (ii), using Ex.1, is asserting that all ``clockables
are writables''. The ordinals $\tau$ (or reals $y$)  in  (iv)  are called
``accidental'', and those in  (iii)  ``eventually writable''. See,
{\eg}{\cite{W09}}.\\

\pf{\em (i) } is immediate from the definition of the pair $(\zeta,\Sigma)$ and then $\lambda$.\\
{\em(ii) } If for a program we have ``$P_e(0)$ {\em halts}'' this is a $\Sigma_1$-statement (``There is $y$ that codes a wellordered sequence of snapshots, with the final snaphot in the halting state".) This yields the first equality. The second is Lemma \ref{writables}.\\
{\em(iii)} Now ``$P_e(m)\da$''  is a $\Sigma_2$-statement in $L_\Sigma$: $\ex \tau\all\g>\tau (${\em ``The OT of the machine $P_e(0)$ is unchanged at time $\g)$"}. This goes down to $L_\zeta$ and so the least such $\tau$ is less than $\zeta$. The second equality is similar.\\
{\em(iv)} For any computation $P_e(m)$ if $z\in WO$ appears at some stage $\tau$ in its computation, then $z$ is definable over $L_\tau$, and thus $||z||<\tau^+$ the next admissible above $\tau$ ($z$ being a wellorder in the admissible set $L_{\tau^+}$). We have seen that the TM produces at some stage or other, all ordinals $\tau$ less than the limit of admissibles $\Sigma$, and thus the $\sup$ in the statement of $(iv)$ is no more than $\S$. \mbox{ }
\qed

\subsection{Infinite Time Jump operator}

\begin{definition}
  (The infinite time jump ${\mathsf{iJ}}$)\\
  (i)  We write $\{e\} (\tmmathbf{m}, \tmmathbf{x}) \da$ if the $e$'th
  ittm-computable function with input $\tmmathbf{m}, \tmmathbf{x}$ has a fixed
  output $c \in \cant$, in which case we write $\{e\} (\tmmathbf{m},
  \tmmathbf{x}) = c$.\\
  (ii)  We then define $\mathsf{iJ}$ by:
  \[ {\mathsf{}} \mathsf{iJ} (e, \tmmathbf{m}, \tmmathbf{x}) = \left\{
     \begin{array}{ll}
       1 & \text{if } \{e\} (\tmmathbf{m}, \tmmathbf{x}) \hspace{0.17em} \da
       \hspace{0.17em} ;\\
       0 & \text{otherwise } (\text{for which we write } \{e\}(\tmmathbf{m},
       \tmmathbf{x}) \ua \hspace{0.17em}) .
     \end{array} \right. \]
  \ \ \ \ \ \ \ \ \ \ \ \ \ \ \ \ \ \ \ \ \ \ \ \ ${\mathsf{}}
  \mathsf{iJ} (y) = y$ if $y$ is not of the form $\pa{e, \tmmathbf{m},
  \tmmathbf{x}}$.
\end{definition}

The functional $\mathsf{iJ}$ then is the counterpart of the standard Turing machine
operator $\mathsf{oJ}$.

\begin{definition}
  For $x$ a real, the complete  (ordinary){\tmem{ ittm-semirecursive-in-$x$
  set}}, denoted by $\tilde{x}$, is the set of integers $\{e \mid \{e\}(e, x)
  \da \}$.
\end{definition}

(`Ordinary' here is simply to contrast with the higher type
ittm-semirecursive sets to come later.)
The following is a consequence of the  (relativized to $x$)
  $\lambda$-$\zeta$-$\Sigma$-Theorem above.
\begin{lemma}
  $\tilde{x}$ is recursively
  isomorphic to the complete $\Sigma_2$-Theory of $L_{\zeta^x} [x]$.
\end{lemma}

\section{Generalised type-2 ittm-recursion}

\subsection{Generalised ittm-recursion in a type 2 functional $\mathsf{I}$}

In the Kleenean recursion in type-2 functionals, in {\cite{Kl62b}},
{\cite{Kl62a}} (building up an equivalent approach to {\cite{Kl59}} and
{\cite{Kl63}})  a successful computation  (meaning one with output)  could be
effected by imagining tm's placed at nodes on a wellfounded tree, with
computations proceeding at nodes that make computation calls to a lower node,
seeking the value of some $x (k)$ say. The computation time at each node,
regarding each call to a lower node as being just one step in the computation
of the calling node, is then finite.  (For otherwise the computation at the
node is never completed and the whole overall computation will fail.)  An
overall computation may also fail by instituting a series of calls to
subcomputations that form an infinite descending path in the tree. In such
cases the machines on the path all hang after finitely many steps, all waiting
for data to be passed up from the immediate subcomputation it has called.

In the ittm case we may again conceive of an overall or master ittm
computation taking place at the top level; such a computation may take
infinitely many steps in time, and will be considered as successful if its
output tape is fixed from some point in time onwards. The master computation
may make queries of a type-2 functional $\mathsf{I}$ in which the computation
is considered ittm-recursive. It may call subcomputations of exactly the same
type: ittm's with the capability to make oracle queries of $\mathsf{I}$.

We give a more detailed description of this as a representation in terms of
underlying ittm's. $\{e\}^{{\mathsf{I}}} (\tmmathbf{m}, \tmmathbf{x})$
will represent the $e$'th program in the usual format, say Turing transition
tables, but designed with appeal to oracle calls possible. We are thus
considering computations of a partial function $\{e\}^{{\mathsf{I}}} :
\,^k \omega \times \hspace{0.17em}^l (^{\omega} 2) \rightarrow \omega$. Such a
computation has potentially computation time, or stages, unbounded in the
ordinals.

The computation of $P^{\mathsf{I}}_e (\tmmathbf{m}, \tmmathbf{x})$ proceeds in
the usual ittm-fashion, working as a tm at successor ordinals and taking
$\liminf$'s of cell values {\etc} at limit ordinals. \ At a time $\alpha$ an
oracle query may be initiated. We may conventionally fix that the real number
subject to query is that infinite string on the even numbered cells of the
scratch type. If this string is $(f, m, y_0, y_1 \ldots,)$ then setting $y =
y_0, y_1, \, \ldots$, the \tmtextit{query} or \tmtextit{oracle call} which we
shall denote $Q^{\mathsf{I}} (f, m, y)$ is the question: ?\tmtextit{What is
$\mathsf{I} (z)$ where } $P^{\mathsf{I}}_f (m, y) \da z$ ? and at stage
$\alpha + 1$ receives the value $\mathsf{I} (z)$. If it is not the case that
$P^{\mathsf{I}}_f (m, y) \da \hspace{0.17em} z$ for any $z$, {\ie}, it fails
to have a fixed output, then there is no $z$ to which $\mathsf{I}$ can be
applied, and the overall computation fails.  (We could try to stay closer to
the Kleenean setting, where a tree branches infinitely often downwards, to
potentially compute some $z \in\,^{\omega} \omega$, {\tmem{via}} $z (0), z (1),
\ldots$ in turn, and then can ask for $\mathsf{I} (z)$. There, if any one of
the single computations $z (k)$ failed, then the query to $\mathsf{I}$ did not
take place, and the overall computation failed. But one thing we have with
ittm computation is plenty of time, so we can, and do, amalgamate the
individual computations $z (k)$ as simply one computation of all of $z$.)

We can determine its effect as follows {\via}an inductive operator $I$. Just
as the Kleene equational calculus can be seen to build up in an inductive
fashion a set of indices and equational strings $\Omega [\mathsf{I}]$ for
successful computations recursive in $\mathsf{I}$ (see Hinman {\cite{Hi78}},
pp. 259-261), so we can define the fixed point of a monotone operator $I =
I^{\mathsf{I}}$ on $(\omega \times \omega^{< \omega} \times
(\omega^{\omega})^{< \omega}) \times \omega^{\omega}$ which will give us the
successful ittm-computations recursive in $\mathsf{I}$. (We blur distinctions
between Cantor and Baire space. Note that we have defined the outputs here as
reals in Baire space, rather than just integers in $\omega$ which the notion
of a Type-2 functional would seem to require. However with ittm's, just as in
the comment above, to compute a $z \in \omega^{\omega}$ is just to compute the
sequence $z (0), z (1) \nocomma, \ldots$ which we can do here, and may
consider the characteristic function of the graph of $z$ as an element of
Cantor space. So by doing this we simply acknowledge that fact of life for
ittm's.

\begin{definition}
  \label{3.8} We set $I (X) =$:
  \[ \begin{array}{l}
       \left\{ \langle \pa{e, \tmmathbf{m}, \tmmathbf{x}}, z \rangle \right|
       P_e^X (\tmmathbf{m}, \tmmathbf{x}) \da z \text{ is an ittm-computation
       making only oracle calls }\\
       \quad \quad \quad Q^X (e', \tmmathbf{m'}, \tmmathbf{x'})
       \hspace{0.17em} \hspace{0.17em} \text{and receiving back } \mathsf{I}
       (z')  \hspace{0.17em} \text{where } X (\pa{e', \tmmathbf{m'},
       \tmmathbf{x'}} \nobracket \nobracket) = z' \hspace{0.17em} \}.
     \end{array} \]
  {\nod}As this is monotone, we may let
  
  $I^{-1} = \emp=I^{<0}$; $I^{< \alpha} = \bigcup_{\beta < \alpha} I^{\beta}$ \&
  $I^{\alpha} = I (I^{< \alpha})$ in the usual way, and reach a least fixed
  point $I^{\infty}$.
\end{definition}

\begin{definition}
  The {\tmem{rank}} of a defined computation, $\rho^{\mathsf{I}} \left(
  \langle \pa{e, \tmmathbf{m}, \tmmathbf{x}}, z \rangle \right)$ is the least
  $\alpha$, if it exists, \ such that $\langle \pa{e, \tmmathbf{m},
  \tmmathbf{x}}, z \rangle \in I^{\alpha}$. We often abbreviate this as
  $\rho^{\mathsf{I}} (e, \tmmathbf{m}, \tmmathbf{x})$ with the $z$ understood
  but unspecified.
\end{definition}

Then:

\begin{definition}
{\sc[The \tmtextbf{$\{e\}$}'th function partial
generalised-ittm-recursive in \tmtextbf{$\mathsf{I}$}]}\\
  Using $I^{\infty}$:
  $$\{ e \}^{\mathsf{I}} (\tmmathbf{m}, \tmmathbf{x}) \mbox{ is {\tmem{defined}}, or
  {\tmem{convergent}}, with output $z$ iff  }I^{\infty} (\pa{e, \tmmathbf{m},
  \tmmathbf{x}}) = z.$$
  
  In which case we set $ \{ e \}^{\mathsf{I}} (\tmmathbf{m}, \tmmathbf{x}) =
  z$ or write $\{ e \}^{\mathsf{I}} (\tmmathbf{m}, \tmmathbf{x}) \da z$.
  Otherwise it is {\tmem{undefined}}, or  {\tmem{divergent}}, and we write $ \{ e \}^{\mathsf{I}} (\tmmathbf{m}, \tmmathbf{x}) \ua. $ $\{e\}^{\mathsf{I}}$ is
  {\tmem{generalised-ittm-recursive in}} \tmtextbf{$\mathsf{I}$} \ if it is
  partial generalised-ittm-recursive in \tmtextbf{$\mathsf{I}$} and total.
\end{definition}

\begin{definition}
  For functionals ${\mathsf{I}, \mathsf{J}}$ we say $\mathsf{I} \leq \mathsf{J}
  $ (``$\mathsf{I}$ is {\tmem{(ittm-generalised) partial recursive in}}
  $\mathsf{J}$'')  if there is $e \in \mathbb{N}$ so that \ $\mathsf{I} = \{ e
  \}^{\mathsf{J}}$. We write $\mathsf{I} \equiv \mathsf{J}$ \ if both
  $\mathsf{I} \leq \mathsf{J}$ and $\mathsf{J} \leq \mathsf{I}$ hold.
  
  The functional $\mathsf{I}$ is {\em recursive in} $\mathsf{J}$  if
  it is partial recursive in $\mathsf{J}$ and total. A relation $\mathsf{R}$ is
  {\tmem{recursive in}} $\mathsf{J}$ if the characteristic function $\mathsf{K_R}$
  is recursive in $\mathsf{J}$.
  
  (ii)  A relation $\mathsf{R} $ is {\tmem{semi-recursive in}} $\mathsf{J}$ if it
  is the domain of a function partial recursive in $\mathsf{J}$.
\end{definition}

\nod The following are straightforward.
\begin{lemma}
  \label{reimage} If $B\sset\omega$ is semi-recursive in a functional
  $\mathsf{I}$ then $B= Im(f)$ for a  partial function $f:\omega \imp \omega$ that is partial recursive in $\mi$.  
  \end{lemma}

\begin{lemma}
  \label{Ex2.2a} The class of relations semi-recursive in a functional
  $\mathsf{I}$ is closed under  universal  number
  quantification $\all^{\omega}$, and in particular under
  $\cap$.
  \end{lemma}
\begin{definition}
For $\mj$ a functional we define the (lightface) {\em halting set for} $\mj$ as 
$$ H^\mj\dfs \{e\mid {e}^\mj(e)
\da \}.$$
%and a boldface version:
%$$ {\mbox{\boldmath{$H^\mj$}}}\dfs \{(e,\vec m, \vec X)\mid {e}^\mj
%(e,\vec m, \vec X)\da\}.$$
\end{definition}

\begin{definition}
  The  (top-level)  length of a computation $P_e^{\mathsf{I}} (\tmmathbf{m},
  \tmmathbf{x})$ in a type-2 oracle $\mathsf{I}$ is the least $\sigma_0 =
  \sigma_0^{\pa{\mathsf{I}, e, \tmmathbf{m}, \tmmathbf{x}}}$ (when defined) 
  so that the snapshot at time $\sigma_0$ of $P_e^{\mathsf{I}} (\tmmathbf{m},
  \tmmathbf{x})$ is the repeat of some earlier snapshot $\zeta_0 =
  \zeta_0^{\pa{\mathsf{I}, e, \tmmathbf{m}, \tmmathbf{x}}} < \sigma_0$, and so
  that the snapshot at $\sigma_0$ recurs unboundedly in $\mathrm{On}$.
\end{definition}

Again, by an easy L\"owenheim-Skolem argument, the ordinal $\sigma_0$ is
countable. Thus the snapshots of the cell distributions in $(\zeta_0,
\sigma_0]$ form the final loop which infinitely repeats thereafter. \ Actually
this top-level length of loop is sometimes is of less interest than the {\tmem{overall length}} of
the computation - to be defined below. Both of these will be undefined if the
computation tree describing $P_e^{\mathsf{I}} (\tmmathbf{m}, \tmmathbf{x})$ is
illfounded. We give here a more detailed description of these trees.

\

Continuing the discussion above, the $\pI{f}$ 'th function on input $m, y$
say, has the opportunity to make oracle calls, and we shall thus have a
\tmtextit{tree} representation of calls made. \ We wish to represent the
overall order of how such calls are made, and indeed the ordinal times of the
various parts of the computation as it proceeds.$\mathsf{}$

\subsubsection{Computation trees \ $\mathfrak{T} =
\mathfrak{T^{\mathsf{I}}} (e, \tmmathbf{m}, \tmmathbf{x})$}

Overall we have a `linear' mode of evaluation of the \tmtextit{computation
tree} - also called a {\em tree of subcomputations}. In particular we
should like to keep track of an \tmtextit{overall length of computation}. This
overall length will be the length not of the top node only,  (which we may
think of as the `master computation', and receives its replies to oracle
queries immediately in one step only)  but as of the whole computation when
the lengths of the computations at lower nodes of the tree, which we regard as
actually performing the sub-computations of the form $P^{\mathsf{I}}_f (y) \da
z$ in order to obtain $\mathsf{I} (z)$, are then also taken into
consideration. It will rapidly be seen that the structure of the tree
$\mathfrak{T} = \mathfrak{T}^{\mathsf{I}} (\tmmathbf{e, m}, \tmmathbf{x})$ of
a convergent computation $P^{\mathsf{I}}_e (\tmmathbf{m}, \tmmathbf{x}) \da z$
reflects how subcomputations arrive into the fixed point of the induction Def.
\ref{3.8}, and thus the rank of this wellfounded tree will be the ordinal
$\left. \rho^{\mathsf{I}} \left( \tmmathbf{\la \la e, m}, \tmmathbf{x}
\right.,\ra z \ra \right) . \nosymbol$ Thus although the computation is most
easily represented by a tree, we may think of the computation as a linear
process as we visit each node of the tree in turn.

We therefore make the following conventions. During the calculation of $\pI{e}
(\tmmathbf{m}, \tmmathbf{x})$ the initial calculation takes place at the
topmost node $\nu_0$ which we declare to be \tmtextit{at Level }$0$, in our
computation tree $\mathfrak{T} = \mathfrak{T}^{\mathsf{I}} (\tmmathbf{e, m},
\tmmathbf{x})$.  (We set $e_0 = e \nocomma, n_0 = \tmmathbf{m}, y_0 =
\tmmathbf{x}$ and pretend that this computation and all the oracle calls below
are only for single number and real variable, merely for ease of
presentation.) \ Let us suppose the first instruction for an oracle query
concerning $\pI{e_1} (n_1, y_1)$ is actioned at stage $\delta_0$ in the
computation of $\pI{e_0} (n_0, y_0)$. The tree $\mathfrak{T}$ will then have a
node $\nu_1$ below $\nu_0$, labelled with $\pa{e_1, n_1, y_1}$ and we declare
the computation $\pI{e_1} (n_{1,} y_1)$ to be performed at this Level $1$.
Thus `control' of the overall process is defined to be at the level of the
node $\nu_1$ at stage $\delta_0 + 1$. The `time' for this sub-computation,
starting thus at $\delta_0 + 1$of course starts locally at its `$t = 0$' -
although each stage is also thought of as one more step in the overall length
of the computation above: namely of $\pI{e} (\tmmathbf{m}, \tmmathbf{x})$. \
Suppose $\pI{e_1} (n_1, y_1)$ makes no further oracle calls and the least
stage at which it exhibits looping behaviour, according to Def. \ref{3.1} is
$\sigma_1$. If there is a settled output of $\pI{e_1} (n_1, y_1)$, $z$ say,
then the correct value $\mathsf{I} (z)$ is then passed back up to Level 0 at
the next stage, that is $\delta_0 + 1$ {\tmem{but only in terms of the stages
of computation at the top level}}, and the master computation proceeds to its
next step at this Level 0.

However we deem that $\delta_0 + 1 + \sigma_1 + 1$ steps have occurred so far
towards the final \tmtextit{overall length}, or $H = H (e, \tmmathbf{m},
\tmmathbf{x})$ of the calculation, that is, of what will be $\pI{e}
(\tmmathbf{m}, \tmmathbf{x})$ if it is successful. To be clear: at stages
$\delta \in (\nobracket$ $\nobracket \delta_0 + 1 + \sigma_1]$ the computation
is at $\nu_1$, whilst in the interval $[0, \delta_0] $and at $\delta_0 + 1 +
\sigma_1 + 1$ it is at $\nu_0$.

However if $\mathsf{e_1} (n_1, y_1)$ itself has made an oracle query, let us
suppose the first such was $? Q^{\mathsf{I}} (e_2, n_2, y_2)$?,
then a new node $\nu_2$ is placed below $\nu_1$ labelled with $\pa{e_2, n_2,
y_2}$ (the label is also part of $\mathfrak{T}$) . If this piece of
computation at $\nu_2$ is successful, that it has settled output $z'$ say, \
and if we suppose it made no oracle
calls,  and took $\sigma_2$ steps to
exhibit looping behaviour, then the value $\mathsf{I} (z')$ is passed back up
to $\nu_1$; lastly the overall length of $\pI{e_2} (n_2, y_2)$ is $\sigma_2$
and then those $\sigma_2$ steps will have to be added to the overall length of
calculation for $\pI{e} (\tmmathbf{m}, \tmmathbf{x})$, being added as they
are, to the top-level length of $\pI{e_1} (n_1, y_1)$.

If the computation $\pI{e} (\tmmathbf{m}, \tmmathbf{x})$ is defined then we
shall have as its computation tree $\mathfrak{T} = \mathfrak{T} (e,
\tmmathbf{m}, \tmmathbf{x})$ a finite path tree  (with potentially infinite
branching)  and some countable rank. $\mathfrak{T}$ will be labelled with
nodes $\{\nu_{\iota} \}_{\iota < \eta (\mathfrak{T})}$ that are visited by the
computation in increasing order  (with backtracking up the tree of the kind
indicated) . Thus $\nu_{\iota}$ is first visited only after all $\nu_{\tau}$
have been visited for $\tau < \iota$.  (This is the sense in which the
computation can be considered as linear after all.)  The $\beta$'th oracle
call to Level $k$ will generate a node we picture as placed to the right of
those so far at Level $k$ (meaning to the right of those with lesser indices
$\alpha < \beta$ at that level). When a subcomputation at a node successfully
finishes, then control of the overall computation is envisaged as passing one
level up to the node immediately above. As the computation progresses it
traverses the tree in the order of the indices on the nodes just described. We
could say that \tmtextit{`control'} of the process is \tmtextit{at a node}
$\nu_{\iota}$ (or is \tmtextit{at a level}) \tmtextit{at time $t$} in the
overall length, if the current sub-computation is running at the node  (or at
a node at that level)  at that time $t$.  (See the next definition.)

The tree will thus have a linear leftmost branch, before any branching occurs.
Further, for a well-founded tree we may define the \tmtextit{overall length \linebreak
function} $H = H (\mathsf{I}, e, \tmmathbf{m}, \tmmathbf{x})$ as above,
formally by recursion on the rank of nodes as the length of the computation.

\begin{definition}
  \label{H}(The overall length function) $H = H (\mathsf{I}, e, \tmmathbf{m},
  \tmmathbf{x})$ is defined by recursion on the rank $\rho^{\mathsf{I}} (e,
  \tmmathbf{m}, \tmmathbf{x})$. Let $\sigma_0 = \sigma_0^{\pa{\mathsf{I}, e,
  \tmmathbf{m}, \tmmathbf{x}}}$ and suppose that $\pI{e} (\tmmathbf{m},
  \tmmathbf{x})$ makes sub-computation calls $\pI{e_{\iota}} (n_{\iota,}
  y_{\iota})$ at times $\tau_{\iota} < \sigma_0$ for $\iota < \theta \leq
  \sigma_0 $. Then
  
  $H (\mathsf{I}, e, \tmmathbf{m}, \tmmathbf{x}) \dfs$ \ $\sum_{\iota <
  \theta} ((\tau_{\iota} - \sup \{ \tau_{\xi} \mid \xi < \iota \}) \noplus + 1
  + H (\mathsf{I}, e_{\iota}, n_{\iota}, y_{\iota})) + (\sigma_0  - \sup \{
  \tau_{\iota} \mid \iota < \theta \}) .$
\end{definition}

Then $H$ gives simply the total ordinal length of the whole computation
together with its subcomputations as if laid out in a linear fashion. Note that
$H (\mathsf{I}, e, \tmmathbf{m}, \tmmathbf{x})$ is defined as long as 
$ \mathfrak{T}^{\mathsf{I}} (\tmmathbf{e, m},
\tmmathbf{x})$ is wellfounded.

\begin{definition}
  \label{D2.5}(i) \ The {\tmem{level}} of a node $\nu_{\iota}$ is the length
  of the path in the tree from $\nu_0$ to $\nu_{\iota}$.\\
  (ii)  By {\tmem{Level}} $n$ we accordingly mean the set of nodes in the tree
  with level $n$.
  
  {\nod}(iii)  The {\tmem{node of the computation}} $\pI{e} (\tmmathbf{m},
  \tmmathbf{x})$ {\tmem{at time}} $\alpha < H (\mathsf{I}, e, \tmmathbf{m},
  \tmmathbf{x})$, denoted {\linebreak}$\nu (\alpha) = \nu (\mathsf{I}, e,
  (\tmmathbf{m}, \tmmathbf{x}), \alpha)$, is the node $\nu_{\iota}$ at which
  the overall computation is being performed at time $\alpha$, and has label
  $\pa{e_{\nu (\alpha)}, \tmmathbf{m}_{\nu (\alpha)}, \tmmathbf{x}_{\nu
  (\alpha)}}$. The {\tmem{level}} $\Lambda (\alpha) = \Lambda (\mathsf{I}
  \nocomma, e, (\tmmathbf{m}, \tmmathbf{x}), \alpha)$ is the level of $\nu
  (\alpha)$.
  
  {\nod}(iv) The {\tmem{current snapshot at time}} $\alpha < H (\mathsf{I},
  e, \tmmathbf{m}, \tmmathbf{x})$ is denoted $\pa{I (\alpha), R (\alpha),
  \pa{C^{\nu (\alpha)}_i (\alpha) \mid i < \omega}}$, and equals the snapshot
  $\pa{I (\bar{\alpha}), R (\bar{\alpha}), \pa{C_i (\bar{\alpha}) \mid i <
  \omega}}$ of the sub-computation $\{ e_{\nu (\alpha)} \}^\mathsf{I} (\tmmathbf{m}_{\nu
  (\alpha)}, \tmmathbf{x}_{\nu (\alpha)})$, where $\bar{\alpha} = \alpha -
  \pi_{\nu (\alpha)}$, and which was initiated at the overall time $\pi_{\nu
  (\alpha)}$.
  
  If $\mathsf{I} = \mathsf{E}$ (see Def. \ref{DefE} below) it can be omitted.
\end{definition}

Thus for a {\tmem{defined}} (or `{\tmem{successful}}') computation, at any
time the level is a finite number  (`depth' would have been an equally good
choice of word). An {\tmem{undefined,}} or {\tmem{failed}}, or
\tmem{unsuccessful}, \tmtextit{ computation} is one in which a
sub-computation call resulting in a calculation at some node fails to produce
an output $z$ (and so no value $\mathsf{I} (z)$ can be returned to the level
above)  either (a) because some subcomputation produced no convergent output
or  (b) $\mathfrak{T} (\mathsf{I}, e, \tmmathbf{m}, \tmmathbf{x})$ is
illfounded  (with a rightmost path of order type then $\omega$) ; or else (c)
the topmost computation itself fails to have convergent output, {\ie} to have
a stable output tape. If either of these kinds of failure occur we denote this by 
 $\pI{e} (\tmmathbf{m},
  \tmmathbf{x})\ua$.

In (iv) the current snapshot is thus the snapshot of the
machine that is running at time $\alpha$ (thus computing $\{ e_{\nu (\alpha)}
\} (\tmmathbf{m}_{\nu (\alpha)}, \tmmathbf{x}_{\nu (\alpha)})$). This machine
started running when it's local time was $t = 0$ of course, but in the overall
picture of things, it starts at time $\pi_{\nu (\alpha)}$, and $\bar{\alpha}$
simply gives how many steps it has run at overall time $\alpha$.

The following lemma incorporates level and cell value facts from the
description of the trees and how control passes from level to level just
given. The point to notice, {\eg}in {\tmem{(i)}}, is if a computation at a
node $\nu (\alpha)$ locally reaches its first repeating point
$\sigma_{\alpha}$ say, then control is immediately passed back to one level
above; thus $\nu$ and $\Lambda$ are decreased. So the liminf in question in
{\tmem{(i)}} for $\Lambda$ is over levels of computation at the current level
or below in the tree, and thus cannot contribute unboundedly in $\lambda$
smaller integers to the $\liminf \nosymbol$. Similar considerations justify
that $\nu (\lambda) = \liminf_{\alpha \rightarrow \lambda} \nu (\mathsf{I}, e,
(\tmmathbf{m}, \tmmathbf{x}), \alpha)$: any subcomputation calls are to nodes
with higher nodal index: to $\nu (\beta)$ greater than what will be $\nu
(\lambda)$ on a tail of $\beta$ below $\lambda$.

\begin{lemma}
  \label{L2.6}Let $\tmop{Lim} (\lambda)$. The computation $\pI{e}
  (\tmmathbf{m}, \tmmathbf{x})$, if not divergent by stage $\lambda$,
  satisfies:
  
  (i) $\nu (\lambda) = \nu (\mathsf{I}, e, (\tmmathbf{m}, \tmmathbf{x}),
  \lambda) = \liminf_{\alpha \rightarrow \lambda} \nu (\mathsf{I}, e,
  (\tmmathbf{m}, \tmmathbf{x}), \alpha)$;
  
  and so: $\Lambda (\lambda) = \Lambda (\mathsf{I}, e, (\tmmathbf{m},
  \tmmathbf{x}), \lambda) = \liminf_{\alpha \rightarrow \lambda} \Lambda
  (\mathsf{I}, e, (\tmmathbf{m}, \tmmathbf{x}), \alpha)$
  
  (ii)  If $\nu = \nu (\mathsf{I}, e, (\tmmathbf{m}, \tmmathbf{x}), \lambda)$,
  and if $\pI{e_{\nu}} (\tmmathbf{m}_{\nu}, \tmmathbf{x}_{\nu})$ is the
  subcomputation at Level $k = \Lambda (\lambda )$, currently being run at
  stage $\lambda$, then if $\pa{C^{\nu}_i \mid i < \omega}$ are the cell
  values of this subcomputation, then $C^{\nu}_i (\lambda) = \liminf_{\alpha
  \rightarrow \lambda \tmcolor{black}{\tmcolor{red}{\tmcolor{black}{, \nu
  (\alpha) = \nu}} {}}} C_i^{\nu (\alpha)} (\alpha)$.
\end{lemma}

\nod It is easy to construct indices $e_1, e_2$ so that:\\
  (i)  $\{ e_1 \}^{\mathsf{I}} (0) \da$ in $< \,
  \omega_{}$ steps  (at the top level)  for any $\mathsf{I}$, but $H (e_1, 0,
  0) \geq \zeta$.\\
  (ii)  $\pI{e_2} (0) {} \ua$ (for any
  $\mathsf{I}$).\\

\begin{lemma}
  \label{L2.18}There is a p.r. function $f$ such that for any $\mathsf{I}$,
  $e, k \nocomma, \tmmathbf{m}, \tmmathbf{x}$
  
  $$\{ f (e, k)
  \}^{\mathsf{I}} (\tmmathbf{m}, \tmmathbf{x}) \simeq \{ e \}^{\mathsf{I}}
  (\tmmathbf{m}, \tmmathbf{x}, \lambda n. \{ k \}^{\mathsf{I}} (n,
  \tmmathbf{m}, \tmmathbf{x}))$$
  
\nod  and hence the functions partial recursive in $\mathsf{I}$ are closed under
  functional substitution.
\end{lemma}

{\pf}The index $f (e, k)$ is for the procedure that does the following: (A) it
first simulates for $n = 0, 1, 2, \ldots$ in turn on a scratch tape the calculation $\{ k \}^{\mathsf{I}} (n,
\tmmathbf{m}, \tmmathbf{x}) \da k_n$, for some $k_n \in \omega,$ and if so
records the value on a scratch tape. \ If $\{ k \}^{\mathsf{I}} (p,
\tmmathbf{m}, \tmmathbf{x}) \ua$ for some $p\in\omega$, then the R.H.S. fails to compute
anything (as $\lambda n. \{ k \}^{\mathsf{I}} (n, \tmmathbf{m}, \tmmathbf{x})$
is not total). Lastly if for some $p \, \{ k \}^{\mathsf{I}} (p, \tmmathbf{m},
\tmmathbf{x}) \da y \notin \omega$ then we perform some fixed trivial program code
with an illfounded tree to ensure the non-totality of $\lambda n. \{ k
\}^{\mathsf{I}} (n, \tmmathbf{m}, \tmmathbf{x})$.

Otherwise all is well and we have eventually written some real $z = (k_0, k_1, \ldots)
\in\,{^{\omega} \omega}$ on the scratch tape. (B) Simultaneously the current contents of the scratch tape, $z'$ say, are used as input to a computation  $\{ e \}^{\mathsf{I}} (\tmmathbf{m}, \tmmathbf{x}, z')$. But this latter computation is reworked each time $z'$ changes; however eventually $z'$ is the intended $z$ above, and   the
calculation $\{ e \}^{\mathsf{I}} (\tmmathbf{m}, \tmmathbf{x}, z)$ has its output as the final value on the overall OT.  \mbox{ } {\qed}

\begin{lemma}
  (i)  There exists a (ittm-)recursive function $\mathsf{H}$, and a
  function $\mathsf{F}$ partial recursive in $\mathsf{H}$ such that $\mathsf{F} $ is not
  partial recursive.  (We should not be surprised at this.)
  
  (ii) If for any $\mathsf{H}$ recursive in $\mathsf{I}$,
  and $\mathsf{F}$ partial recursive in $\mathsf{H}$ then $\mathsf{F} $ is partial
  recursive in $\mathsf{I}$.
\end{lemma}

{\nod}Usual methods prove an $S^n_m$-Theorem and more particularly:

\begin{theorem}[The generalised ITTM -Recursion Theorem]
 \mbox{ }\\ If $F (e, \tmmathbf{m, x})$ is ittm-recursive in ${\mathsf{I}}$,
  there is $e_0 \in \omega$ so that
  \[ \{e_0 \}^{{\mathsf{I}}} (\tmmathbf{m, x}) = F (e_0, \tmmathbf{m,
     x}) . \]
\end{theorem}

\begin{lemma}
  \label{snapshot}The computation $\nobracket \{ e \}^{\mathsf{I}}
  \tmmathbf{(m \nobracket}, \tmmathbf{x})$ is absolute to $L [\mathsf{I},
  \tmmathbf{x}]$. \ In the above notation, there is a function $S (\alpha) = S
  (\mathsf{I}, \tmmathbf{\tmmathbf{m}, x}, \alpha)$ for $\alpha \in
  \tmop{On}$, (the `snapshot function') so that
  
  (i) $S (\beta) = \pa{\nu (\beta), \Lambda (\beta) \nocomma, \la I (\beta), R
  (\beta), \pa{C^{\nu (\beta)}_i \mid i < \omega} \ra}$ ; 
  
  (ii) {\pa{S
  ({\beta}){\mid}{\beta}<{\alpha}}}$\in \Delta_1^{J_{\alpha} [\mathsf{I},
  \tmmathbf{x}]}$ ; 
  
  (iii) $S (\alpha) {\in \Sigma^{J_{\alpha} [\mathsf{I}, \tmmathbf{x}]}_2} $.
\end{lemma}

\subsection{The functional $^2 \mathsf{E}$}

This is the simple functional of existential quantification. Recall that we
are representing elements of Baire space, so $\text{ } x \in \; \bai,$ in
$^{\omega} 2$ on the tape as the infinite sequence of $1$'s interspersed with
a string of $0$'s of length $x (n) + 1$. An integer $k \in \mathbb{N}$ is
represented by $\bar{k} \dfs n_k \in \bai$, where $n_k$ is a string of $k + 1$
1's followed by only $0$'s.

\begin{definition}
  \label{DefE} (The functional $ \mathsf{^2 \mathsf{E}}$) \
  (i)  We define $\mathsf{E} = \mathsf{^2 E}:  \bai \imp \bai$
  by:
  \[ {\mathsf{E}} (x) = \left\{ \begin{array}{ll}
       0 & \text{if } x \in \: \bai \wedge \ex n x (n) = 0\\
       1 & \text{} \text{if } x \in \: \bai \wedge \all n x (n) \neq 0
     \end{array} \right. \]
  (ii)  For $R \sset \omega$ we set $\mathsf{E} (R) \dfs \mathsf{E} (K_R)$.
\end{definition}

Recursions in $\mathsf{^2 E}$ are already quite powerful: Kleene showed that
for Kleene recursion $\mathsf{{oJ}} $ and $\mathsf{^2 E}$ are mutually Kleene
recursive in each other. The functional $\mathsf{E}$ was important for much of
the development of this recursion, and a type-2 functional $\mathsf{I}$ for
Kleene was {\tmem{normal}} if $\mathsf{^2 E}$ was Kleene-recursive in
$\mathsf{I}$. Many of the theorems of theory were only valid for normal
functionals. For Kleene recursion the functionals $\mathsf{{oJ}}$ and
$\mathsf{^2 E}$ being equivalent, these were the simplest useful functionals.
However here in the ittm setting normality is trivial with $\mathsf{E}(x)$ being simply computable in $\omega$ steps, 
without the use of any oracle or functional input at all.
But we must build our type 2 recursions recursive in something, so we can take
$\mathsf{{iJ}}$ and $\mathsf{^2 E}$ as being the simplest useful
functionals. We shall see that for ittm-recursion $\mathsf{^2 E}$ and
$\mathsf{{iJ}}$ are likewise mutually recursive. 
% (Note that $\mathsf{^2 E} \leq \mathsf{{iJ}}$ is trivial. Check!)

\

For performing a computation together with all its subcomputations as a
tree, and seeing how the length of computation relates to extendibility in the
$L$ hierarchy, even if the function $\mathsf{I}$ is quite simple, and
constructibly definable, this may have to be done in suitably large admissible
sets. However note that any computation $\{ e \}^{\mathsf{E}} (\tmmathbf{m},
\tmmathbf{x})$ is absolute to $L [\tmmathbf{x}]$.

We have relativized in Definition \ref{Def3} to reals $x$ in the obvious way,
the concept of $x$-$(\Sigma_2)$-extendible pairs $(\xi, \sigma)$. Note that
for such a pair, since $\mathsf{E}$ is $\Delta_1$-definable over $L_{\xi}$, so
$(\xi, \sigma)$ is also an $x$-$\mathsf{E}$-$(\Sigma_2)$-extendible pair in an
obvious sense. We use this without further mention.

\begin{lemma}
  \label{L2.10}If $(\xi, \sigma)$ is an $\tmmathbf{x}$-extendible pair, then for $\left. \ptwoe{e} \tmmathbf{(m
  \nobracket}, \tmmathbf{x} \right)$:
  
  (i) $\nu (e, (\tmmathbf{m}, \tmmathbf{x}), \xi) = \nu (e, (\tmmathbf{m},
  \tmmathbf{x}), \sigma)$ and so $\Lambda (e, (\tmmathbf{m}, \tmmathbf{x}),
  \xi) = \Lambda (e, (\tmmathbf{m}, \tmmathbf{x}), \sigma)$;
  
  (ii)  If $\nu (e, (\tmmathbf{m}, \tmmathbf{x}), \sigma) = \nu$ then in the
  notation above $C^{\nu}_i (\xi) = C^{\nu}_i (\sigma)$ for $i < \omega$;
  
  (iii)  If in  (ii) $\nu = \nu_0$, then $\left. \ptwoe{e} \tmmathbf{(m
  \nobracket}, \tmmathbf{x} \right)$ has entered final looping behaviour by
  stage $\xi$.
\end{lemma}

{\pf}These all follow from the $\Sigma_2 \tmop{Liminf}$ nature listed in
Lemma \ref{L2.6} of the properties of Def. \ref{D2.5}.{\qed}\\

Conversely:

\begin{lemma}
  \label{L2.11} 
  If $\left.
  \ptwoe{e} \tmmathbf{(m \nobracket}, \tmmathbf{x} \right)$ is convergent,
  then there is $(\xi, \sigma)$ an $\tmmathbf{x}$-extendible pair with
  $$\Lambda (e, (\tmmathbf{m}, \tmmathbf{x}), \xi) = \Lambda (e, (\tmmathbf{m},
  \tmmathbf{x}), \sigma) = 0.$$
\end{lemma}

\begin{definition}
  A type-2 functional $\mathsf{I}$ is called {\tmem{suitable}} if the
  $\tmop{ran} \left( \mathsf{I} \rest  \omega \right)$ is not a singleton,
  where we represent $k \in \omega$ as $\bar k$: the infinite sequence of $k $ $1$'s
  followed thereafter by $0$'s: $1^k\smallfrown 0^\omega$.
\end{definition}

In general, $^{< \omega} \omega$ can be interpreted as those infinite strings
from Cantor space that are zero from some point onwards and whose initial part codes a finite sequence
in some manner that the reader may care to provide. Thus by abuse of
notation $\mathbb{N}$ ``$\sset$'' $^{< \omega} \omega$. Then ${\mathsf{iJ}}$ is
suitable. It is easy to see that for any $\mathsf{K} $ there is a suitable $\mathsf{K}'
\equiv \mathsf{K}$.

\begin{lemma}\label{no queries}
The set $Z\dfs \{(e,m,x)\mid \{ {e}
  \}^{\mathsf{J}}(m,x) \mbox{ makes no query calls} \}$  is a recursive in ${\mj}$ set.
\end{lemma}
\pf The following procedure describes a total recursive in $\mj$ function $K_Z
: \omega \times \omega \times\omega^\omega \imp 2$ that is the characteristic function of $Z$. $P$ simulates a run of $ \{ {e}
  \}^{\mathsf{J}}(m,x)$ with a zero on its output tape. As soon as $ \{ {e}
  \}^{\mathsf{J}}(m,x)$ makes a query request,  $P$ changes that zero to a one and does nothing further.\\ \mbox{ } \qed\\

The set $Z$ then catalogues those indices and inputs which are equivalent to those of standard ittm-recursive machine computations.

\def\concat{\smallfrown}
\begin{lemma}  \label{L2.15} Assume  $\mathsf{J} $ is suitable.
  (a)  There is an index $t_0$
  so that for any $s \in \,^{\omega} 2$: \\
  
  $
  \begin{array}{rcll}
  P^{\mathsf{J}}_{t_0} (e, m, y, s ) & =&
  1  \mbox{ if $s$ is the first finally  repeating snapshot in the computation } \{ e \}^{\mathsf{J}} (m,
  y)\\
  & = & 0 \mbox{ otherwise }.
\end{array}
$

  (b)  There is $\bar{e}$ so that for any $e,m,y$ if 
  $\mathfrak{T} (\mathsf{J}, e, {m}, {y})$ wellfounded,
  $P_{\bar{e}}^{\mathsf{J}}(e, m, y) = 1\concat z$ if $\{ e \}^{\mathsf{J}} (m, y) \da z$, and $=
  \, 0$ if 
  $\{ e \}^{\mathsf{J}} (m, y) \ua$.   
  
  \end{lemma}
{\pf} We assume without loss of generality that $\mathsf{\mathsf{J}} \left(
\overline{0} \right) = 0$ and $\mathsf{\mathsf{J}} \left( \overline{1} \right) =
1$.  (It will be apparent what to do if we need to appeal to other values
under suitability of $\mathsf{\mathsf{J}} .$)  

\nod (a) The action of $T_0 \dfs P^{\mathsf{J}}_{t_0} (e, m, y, s)$
does the following:

(i) $T_0$ first ensures the  cells of the OT ($T_0$) \ (the output
tape of $T_0$)  are $0$;

(ii) $T_0$ itself runs the code of $\{ e \}^{\mathsf{J}} (m, y)$ on a scratch
tape,  but instead starting from the presumed snapshot $s$ onwards, using for
this the input  snapshot data $s$, all the while
inspecting the later snapshots $s_{\beta}$ that it
generates.

(iii) If the snapshot $s$ recurs later, with the proviso that no cell
of $s$ which has a $1$, switches $1 \imp 0 \imp 1$, \ $T_0$ sets the
first cell of $\tmop{OT} (T_0)$ to 1. 
(And in this case we can require the machine to actually halt at this point.)
 If the snapshot $s$ does not recur then the first cell remains set on zero (but without necessarily halting).

Thus we see that $P^{\mathsf{J}}_{t_0} (e, m, y, s)$
is convergent with the correct values to finish (a).

\nod (b) We define the following procedure
$P$ which will be realised as a programme $P_{\bar{e}}^{\mathsf{J}}$. The
process $P$ does the following:  (A)  it looks for a snapshot of the top level
of the computation of $P_e^{\mathsf{J}} (m, y) $ that is the first repeating
snapshot in a final loop, using $P_{t_0}^\mathsf{J}$ from part (a) above; 
(B)  A
snapshot that passes these checks can also at the same time be inspected
to see if the OT that it encodes has a convergent value or not. 

In more
detail: $P$ first runs a copy of the program $P_e^{\mathsf{J}} (m, y)$ on a
scratch tape, running the programme instructions coded in $e$. For each top level time $\alpha$ of the simulated run of this
$P_e^{\mathsf{J}} (m, y)$, the simulated snapshots $s_{\alpha}$ at simulated stages 
 $\alpha$ are written by $P$ to a reserved piece of tape $R$
(each snapshot overwriting the contents of $R$ as this is repeatedly done).
Then  $P$ makes the query ?$Q^{\mathsf{J}} (t_0, (e, m, y, s_{\alpha}))?$

If the outcome of the above was that $s_\alpha$ was not  the first final looping
snapshot of the computation, (our assumptions on
$\mathsf{J}$, \ $\mathsf{J} (1 \smallfrown 0) = 0$ \etc, are set to ensure the correct value is returned) then we return to the original simulation which then proceeds to calculate $s_{\a +1} $ and continue as above. However if $s_\a$ is the first snapshot of a final looping computation, we can  run the check again that it is such a first snapshot, but check now whether the OT coded recursively into that
snapshot ever changes. If so, the the first cell of the OT of $P$ is set to zero, and if not $1\concat z$ is written to OT($P$), where $z$ was the contents of the OT part of $s_\a$. (In both these cases we can also require $P$ to halt.)

(Note that at limit
stages $\mu$ of this process the liminf process in the master program 
naturally records the correct liminf
snaphot $s_{\mu}$ in $R$ of the simulated program.) 
Eventually $P$ will reach a snapshot $s_{\xi}$ that will indeed be the start
of a final loop and $\mathsf{J}$ will be returning the correct $0 / 1$ value. Only if 
$\mathfrak{T} (\mathsf{J}, e, {m}, {y})$ is illfounded will this process fail to halt.
\mbox{ }\qed\\

{\nod}Although we build into our framework, following Kleene, that a query
$?Q{\mathsf{^{{\mathsf{K}}}}} (e, m, x) ?$ first does some computation, namely
$\{ e \}^{\mathsf{K}} (m, x)$, and if this is convergent, submits the result to
$\mathsf{K}$, actually we can shortcircuit the process, and the above Lemma shows, that if $\mathsf{K}$ is suitable, we can obtain
by appropriate  queries during a
computation,  convergence/divergence facts, nd indeed output facts. (However there can be no total recursive function $F(e,m,x)$ that returns $0/1$  depending on whether $\mathfrak{T}^\mathsf{K}(e,m,x)$ is well founded or illfounded.) 

\begin{corollary}\label{Cor3.22} Let $\mj$ be suitable.
There is a p.r. function $k$ so that for any $e,\tmmathbf{m}, \tmmathbf{x}$ if $P^\mj_e(\tmmathbf{m}, \tmmathbf{x})\da y$ then   $P^\mj_{k(e)}(\tmmathbf{m}, \tmmathbf{x})\da y$, with the difference that the latter is a  {\em halting} computation.
\end{corollary}

\pf This is a particular case of the preceding Lemma \ref{L2.15}. \qed\\

The generalised ittm recursive computations are thus quite flexible: we can extend any query $Q^\mathsf{J}(e,m,x)$?  to return not just $\mathsf{J}(y)$, where ${e}^\mathsf{J}(m,x)\da y$, to the calling computation, but also the whole real $y$ can be effectively returned. One way is to use the last lemma, or argue as follows.

\begin{lemma}\label{Ex2.6}Let $\mathsf{J}$ be suitable. 
There is a p.r. function $\bar h$ so that
   $P^{\mathsf{J}}_{{\bar h}(e)}
  (\tmmathbf{m}, \tmmathbf{x})$ runs the computation  $P^{\mathsf{J}}_e
  (\tmmathbf{m}, \tmmathbf{x})$ but 
   where every query instruction
  ?$\nobracket Q^{\mathsf{J}} (t, n, y)) ?$ in a computation $P^{\mathsf{J}}_e
  (\tmmathbf{m}, \tmmathbf{x})$, can be replaced with some finite set of
  instructions which initiates a recursive sequence of queries ?$\nobracket
  Q^{\mathsf{J}} (t_i, n, y)) ?$ which effects the writing to a recursive
  slice of the scratch tape $R$, of the sequence of digits $z (i)$ where, if
  it exists, $\{ t \}^{\mathsf{J}} (n, y) \da z$ . \end{lemma}

Thus instead of the query
  returning just the single integer $\mathsf{J} (z)$ we can think of the
  amended program as returning $z$ itself to (a recursive slice of) the
  scratch tape. This `subroutine' is independent of $e$, but with $t_i$ primitive
  recursively dependent on $t$ and $i$ only.

\rem Corollary \ref{Cor3.22} demonstrates that we could have used as part of our basic choice of architecture, properly halting computations as our notion of convergence rather than those with eventually settled OT: the resulting class of generalised computable functions would have been the same.\\

The above  gives a clue as to how to present yet another way of reconfiguring the architecture of generalised ittm recursion. The template we adopted was to look at something very close to Kleene's ideas, in particular keeping the oracle's  $\mj$ as maps into $\w$ rather than $\bai$; hence a subcomputation call or query ends up returning a single digit, $\mj(z)$ for some $z$ which was the convergent result of some query. However as we saw above we could rewrite the program to enable us to get all the digits of $z(i)$ together.  We could have simply specified a template whereby infinitary objects, \ie the reals $z$ as considered, were returned. The following argument shows what to do.  
\begin{lemma}
  \label{Ex2.9}Let $\mathsf{J}$ be suitable. 
   There is a p.r. function $h$ so that in any computation $P^{\mathsf{J}}_e
  (\tmmathbf{m}, \tmmathbf{x})$ any query call ?$\nobracket Q^{\mathsf{J}} (t,
  n, y)) ?$ that occurs throughout the computation tree, is replaced in
  $P^{\mathsf{J}}_{h (e)} (\tmmathbf{m}, \tmmathbf{x})$ by one which first
  returns $z \simeq P^{\mathsf{J}}_t (n, y)$ itself to some recursive slice of
  the scratch tape of the calling computation, before returning the integer $\mj
  (z)$. Hence this enhanced `model of computation' which returns  whole reals in this way,
  is only a form of generalised ittm-recursion.
  \end{lemma}
 \pf We build on Lemma \ref{Ex2.6} using the p.r. function $\bar h$,  which had the effect in a
  computation $P^{\mathsf{J}}_e (\tmmathbf{m}, \tmmathbf{x})$ of any query call
  ?$\nobracket Q^{\mathsf{J}} (t, n, y)) ?$ at the top level being replaced in
  $P^{\mathsf{J}}_{\bar h (e)} (\tmmathbf{m}, \tmmathbf{x})$ by one which
  instead returns the  value $z$. 
   Now, somewhat
  trivially do one more modification: thinking of $z$ as input to the identity
  programme $P^{\mathsf{J}}_{\iota} (z) = \tmop{id} (z)$, make the
  call $? Q^{\mathsf{J}} (\iota, z) ?$ This returns $\mathsf{J} (z)$ after all of $z$, and
  which can be written to a scratch area.
  
  This yields a p.r. function $h_0$ which substitutes this code for the top
  level queries occurring in $P^{\mathsf{J}}_e (\tmmathbf{m}, \tmmathbf{x})$,
  yielding a computation $P^{\mathsf{J}}_{h_0 (e)} (\tmmathbf{m},
  \tmmathbf{x})$. Lastly we use the recursion theorem to show that there is a
  p.r. function $h$ so that $P^{\mathsf{J}}_{h (e)} (\tmmathbf{m},
  \tmmathbf{x})$ has the effect of applying $h_0$ throughout the computation
  tree for $P^{\mathsf{J}}_e (\tmmathbf{m}, \tmmathbf{x})$, to the indices of
  all sub-computation calls. By design we shall have that $P^{\mathsf{J}}_e
  (\tmmathbf{m}, \tmmathbf{x}) \simeq$ $P^{\mathsf{J}}_{h (e)} (\tmmathbf{m},
  \tmmathbf{x})$.
\qed

\begin{lemma}
  Let $\mathsf{J} $ be suitable, then there is an index $p_1$,
    such that $?
  Q^{\mathsf{J}} \left( p_1, \left( e, m, y, \pa{k, l} \right) \right) ?$
  returns 1 if $\{ e \}^{\mathsf{J}} (m, y) \da$ and $\{ e \}^{\mathsf{J}} (m,
  y) (k) = l$ \ \ and 0 otherwise. In particular a programme $P^{\mathsf{J}}$
  can compute the graph of convergent $\{ e \}^{\mathsf{J}} (m, y) .$
\end{lemma}
\pf Straightforward. \qed

\begin{theorem}
  \label{trans}If  $ \mathsf{I} \leq \mathsf{J} \leq \mathsf{K} $ and $\mathsf{K} $ is
  suitable, then $\mathsf{I} \leq \mathsf{K}$. More generally there is a p.r.
  function $f$, so that $\{ e \}^{\mathsf{J}} = \{ f (e) \}^{\mathsf{K}}$.
\end{theorem}

{\pf} We again assume without loss of generality that $\mathsf{K} (0) = 0$ and
$\mathsf{K} (1) = 1$.  (It will again be apparent what to do if we need to appeal
to other values under suitability of $\mathsf{K} .$)  We are given that there are
$g_1 \nocomma, g_2 \in \omega$ so that $\mathsf{I} = \{ g_1 \}^{\mathsf{J}}$ and
$\mathsf{J} = \{ g_2 \}^{\mathsf{K}}$. We show there is a method for finding an index
$g_3$ for a recursion $\mathsf{I} = \{ g_3 \}^{\mathsf{K}}$. \ But more generally,
we actually show how to rewrite, in a p.r. way, a program computing $\{
\bar{e} \}^{\mathsf{J}}$ into one $\{ f (\bar{e}) \}^{\mathsf{K}}$ computing
the same function. This will be done in a uniform manner that is independent of $\mathsf{J}, \mathsf{K}$
(as long as any other $\mathsf{K}'$ under consideration agrees with
$\mathsf{K}$ on $0$ and $1$).

Just note first that any index $e$ codes a finite ittm programme that syntactically may be run as a
programme $P^{\mathsf{K}}_e$ or indeed with any other oracle functional
$\mathsf{K}$ \ - the grammar of the programme does not impose any conditions
on the oracle - it merely asks for values.

 A query at local time $\alpha$, $Q^{\mathsf{J}}
(e_1, m_1, y_1) \nocomma \nocomma$ at level 0, has two phases: it asks if $\{
e_1 \}^{\mathsf{J}} (m_1, y_1) \da z$ \ for some $z$ and secondly, if so, it asks
for $\mathsf{J} (z)$ which is then returned at local time $\alpha + 1$ again at
level 0. To effect the translation of this as a recursion in $\mathsf{K}$ we
slightly modify the mechanism of Lemma \ref{Ex2.9}. There the value of
$\mathsf{J} (z)$ was introduced by the device of the identity function, as the
value returned, following on the subcomputation call $P_{\iota} (z)$. The
original program, $P_e$, for which this was introduced was modified to $P_{h
(e)}$ (for a p.r. $h$). (And all subcomputation calls $Q^{\mathsf{J}} (e_1, m_1,
y_1) \nocomma \nocomma$, by recursion, were modified to $P_{h (e_1)}$ {\etc}
We here just replace $h$ by $h (e_{1,} g_2) $ the index function arising in
the same way, but replacing the identity query \ $? Q^{\mathsf{J}} (\iota, z)
?$ by the code of the program $\{ g_2 \}^{\mathsf{K}} (z)$ using $z$. Doing this recursively throughout the computation, yields a new program, now recursive in $\mathsf{K}$, an
index for which can be found by the Recursion Theorem once more, uniformly in $g_2$;
we let this be then $h (e, g_2)$. Then $f (e) \dfs h (e, g_2)$ is
the p.r. function of the theorem. \hfill{\qed}\\

The construction here only depended in trivial ways upon the suitability of
$\mathsf{K}$.

\begin{lemma}
  \label{L2.16}$\mathsf{{iJ}} \leq \mathsf{E}$. Hence as $\mathsf{E}$ is
  trivially recursive in $\mathsf{{iJ}}$, we have $\mathsf{E} \equiv
  \mathsf{{iJ}} .$
\end{lemma}

{\pf}Adapt the methods above. Exercise.{\qed}\\

From this point on we shall assume our functionals are suitable, unless
otherwise stated. The following is only a starting sample of what we shall later prove (\cf for example, Lemma \ref{4.30}).

\begin{lemma}\label{Ex2.8} There are indices $e_0, e_s$ so that $\{ e_0
  \}^{\mathsf{K}} (e, m, x) \da z^x$, \ $\{ e_s \}^{\mathsf{K}} (e, m, x) \da
  s^x$, where $z^x$ is the $L [x]$-least code for $L_{\xi^x} [x]$ and $s^x$ is
  the $L [x]$-least code for $L_{\sigma^x} [x]$, the least level of the $L
  [x]$ hierarchy with a proper $\Sigma_2$-elementary substructure $L_{\xi^x}
  [x]$. 
\end{lemma}

\subsection{Stage Comparison}

For Kleene Recursion the next move would be to prove a Stage Comparison
Theorem and use this to develop some of the theory of the semi-recursive on
$\mathsf{I}$ sets. We have potentially two options here along two different
axes, to assign ordinals to computations rather as for Kleene recursion,
namely as to which stage $I^{\alpha}$ a computation $\pI{e} (m, x)$ enters
into the monotone inductive definition of all successful computations in
$\mathsf{I}$, or else we may look at the overall convergence time $H
(\mathsf{I}, e, \tmmathbf{m}, \tmmathbf{x})$ it takes for $\pI{e} (m, x)$ to
converge. They yield somewhat differing prewellorderings on computations. We
stay with the former possibility.
The extra power of being able to compute ordinals for lengths of successful
computations, some of the results for Kleene Recursion become then simple in
this context. But we find that the ranking also has implications for computation lengths.

We first note that there are suitable {\tmem{universal semi-recursive in
$\mathsf{I}$ sets}}.

\begin{definition}
  \label{Def2.16}$\mathsf{U}^{\mathsf{I}} (e, \nobracket m, x \nobracket) \equi
  \pI{e} (m, x) \da$;
  
  \ \ \ \ \ \ \ \ \ \ \ \ \ \ \ \ \ \ \ \ $U_y^{\mathsf{I}} (e, m) \equi
  U^{\mathsf{I}} (e, m, y) ;$
  
  \ \ \ \ \ \ \ \ \ \ \ \ \ \ \ \ \ \ \ \ $U^{\mathsf{I}} (e, m) \equi
  U^{\mathsf{I}} (e, m, \bar{0}) .$
  
  We say that $\mathsf{U^{\mathsf{I}}}$ and $U_y^{\mathsf{I}}$ are{\tmem{
  parametrized by}} $\omega$: as the index $e \in \omega$ varies we obtain all
  semi-recursive in $\mathsf{I}$ (or $\mathsf{I}$ and $y$)  sets. 
\end{definition}

\begin{definition}
  We set $\rho^\mi (e, m, x) = \left| \la e, \pa{m, x} \ra \right|^\mi$ to
  be the least $\alpha$ such that $ {\pa{\pa{e,m,x}},z} \in I^{\alpha}$ for
  some $z.$
\end{definition}

Thus we rank convergent computations by that $\alpha$ where they appear in the
inductive definition of convergent computations. 
This gives a {\em  norm} on $\mathsf{U}^{\mathsf{I}}$ and so on any semi-recursive in
$\mathsf{I}$ set $A$.

If $ \pa{e, m, x}
 \notin \dom(I)$, it is convenient to set $\rho^\mi (e, m, x) = \omega_1$.
 %==============FROM GANDYFEEDBACK
 We further adopt the following abbreviating notation: we let greek letters such as
$ \gamma, \delta$ in the sequel stand in for computations $\gamma = \{e\}^\mi
(m,x)$ which we shall also abbreviate as simply $\gamma = (e, m, x)$. If $\gamma =(e, m, x)$
abbreviates the computation $\{ e \}^\mi (m,x)$, which makes a query
call $Q^\mi(e_0, m_0,x_0) ?$ then we say that  $\{e_0\} ^\mi (m_0,x_0) $ is a
{\tmem{subcomputation}} of \ $\{ e \}^\mi (m,x)$, and we abbreviate such a
subcomputation as $\gamma_0 =_{\tmop{df}} (e_0, m_0, x_0)$.

We let  $\gamma_p = (e_p, m_p, x_p) (p < k (\gamma))$ %(where $k (\gamma) \leq \omega$)
enumerate $? Q^\mi (e_p, m_p,x_p) ?$, the query calls occurring in turn in the run of
a computation $\gamma$. Then note that $\rho^\mi (e, m, x) = sup^+\{\rho((e_p, m_p, x_p))\mid \gamma_p = (e_p, m_p, x_p) (p < k (\gamma)) \,\mbox{\em  a subcomputation of } \g\}$.

We then have:

\begin{theorem}
  \label{OrdComp}(Stage Comparison) There is a functional {$\mh = \mh^\mi$},
  partial recursive in $\mathsf{I}$, such that for all $\gamma = \pa{e^0, m^0,
  x^0}$ and $\delta = \pa{e^1, m^1, x^1} :$
  
  $
  \begin{array}{lrcl}
  (i) &\rho^\mi (\gamma) < \omega_1 \wedge \rho^\mi (\gamma) \leq \rho^\mi (\delta)
  &\imp &\mathsf{H} (\gamma, \delta) \simeq 0 ;\\
  
  (ii)& \rho^\mi (\delta) < \omega_1 \wedge \rho^\mi (\delta) < \rho^\mi (\gamma)
  &\imp& \mathsf{H} (\gamma, \delta) \simeq 1.
\end{array} $
\end{theorem}

{\pf} Define a partial recursive in $\mi$ function $F$ as follows, \\

$
\begin{array}{lccl}
F (e, \gamma, \delta) &= &0& \mbox{ if } \rho^\mi (\gamma) = 0\\

&=&1&\mbox{ if }\rho^\mi (\gamma) \neq 0 \wedge \rho^\mi (\delta) = 0\\

&= &0& \mbox{ if }\all \mbox{ subcomps.} \gamma_p \tmop{of}
\gamma \exists \mbox{ a subcomp. } \delta_q\mbox{ of }\delta\\
&&& %\mbox{ such that }
\{ e \}^\mi (\gamma_p, \delta_q) \downarrow 0{\hspace{12em}}(A)\\

&= & 1&\mbox{ if }\exists \mbox{ a subcomp. } \gamma_p \mbox{ of }\gamma
\forall \mbox{ subcomps. } \delta_q \mbox{ of } \delta \\
&&& %\mbox{ such that }
\{ e \}^\mi (\gamma_p, \delta_q) \downarrow 1.{\hspace{12em}}(B)
\end{array}
$

\

$F (e, \gamma, \delta)$ is undefined otherwise.

\

\nod {\tmem{Claim}} $F$ is partial recursive in $\mathsf{\mi}$.

{\pf}of Claim. We give a procedure $P$, recursive in $\mathsf{\mi}$, for
computing $F$.

1) $P$ has inputs $\gamma, \delta$ as above and first checks if they are of
the right form, and then becomes undefined if not.

2) The set $\{(e,m,x)\mid \rho^\mi(e,m,x)=0\}$ is a recursive in $\mi$ set. (See Lemma \ref{no queries}.) $P$ then checks if $\rho^\mi(\g) =0$ and if so $F$ will be made to output $0$, and then otherwise, if   $\rho^\mi(\d) =0$, $F$ will be made to output $1$.  These are then recursive in $\mi$ conditions.

We proceed to:

3)
 As follows: $P$  (i) first simulates a run of $\gamma$ until ?$Q (\gamma_0) ?$ occurs. 
 (As the first base case conditions do not hold, then since $\rho^\mi(\g)>0$, $\g=P_{e^0} (m^0, x^0)$ must make a first query call for a subcomputation
$\gamma_0$ and so the latter is defined) Then we pause before evaluating the subcomputation.
Before
instigating this query, $P$ then starts a simulation of $\delta$, searching
through the called subcomputations $\delta_0, \delta_1, \ldots$ in turn, writing down the values $i_\alpha \dfs  \{ e \}^{\mathsf{I}} (\gamma_0, \delta_{\alpha})$ to an output tape, until, if ever, some $\delta_{\alpha}$, for some $ \alpha < k(\delta)$, is
reached with $\{ e \}^{\mathsf{I}} (\gamma_0, \delta_{\alpha})$ undefined, or
$\{ e \}^{\mathsf{I}} (\gamma_0, \delta_{\alpha}) \downarrow z \notin 2$, or
$\{ e \}^{\mathsf{I}} (\gamma_0, \delta_{\alpha}) \downarrow 0$. If either of
the two former cases hold then $F (e, \gamma, \delta)$ is deemed undefined. If no
such $\delta_{\alpha}$ is found then for all $\alpha < k(\delta)$, \ $\{ e
\}^{\mathsf{I}} (\gamma_0, \delta_{\alpha}) \downarrow 1$. In other words the eventual value on this segment of the scratch tape is $1$. Then  $1$ is written on the output tape of $F (e, \gamma,
\delta) $ which then halts. 
Otherwise ({\ie}with $\{ e \}^{\mathsf{I}} (\gamma_0, \delta_{\alpha})
\downarrow 0$ for some $\alpha <k(\delta)$) $P$ writes $0$ to the OT for $F$, and then  repeats the process with $\gamma_1$ (if the latter is defined), searching for the
least $\delta_{\alpha}$ with $\{ e \}^{\mathsf{I}} (\gamma_1,
\delta_{\alpha})$ undefined or $\{ e \}^{\mathsf{I}} (\gamma_1,
\delta_{\alpha}) \downarrow i_0 \notin 2$, or $\{ e \}^{\mathsf{I}} (\gamma_1,
\delta_{\alpha}) \downarrow 0$. Again in the first two cases $F (e, \gamma,
\delta)$ is undefined. If no such $\delta_{\alpha}$ is found then for all
$\alpha < k(\delta)$, \ $\{ e \}^{\mathsf{I}} (\gamma_1, \delta_{\alpha})
\downarrow 1$. As before $1$ is written to the OT for  $F (e, \gamma, \delta) $ which then halts..
If $\{ e \}^{\mathsf{I}} (\gamma_1, \delta_{\alpha})
\downarrow 0$ for some $\alpha <k(\delta)$, then  $0$ is written to the OT of $P$, which then  continues with the simulation of $\gamma$ until a ?$Q (\gamma_2) ?$ possibly occurs \etc
\\

If for all $\beta < k(\gamma) $ there is $\alpha (\beta) < k(\delta)$ with $\{
e \}^{\mathsf{I}} (\gamma_{\beta}, \delta_{\alpha (\beta)}) \downarrow 0$ then
we have ensured $P$ has settled output $0$ as the value for $ F (e, \gamma, \delta)$.{\qed} (Claim)\\

By the $\mathsf{I}$-recursion theorem there is $\bar{e}$ with $\{ \bar{e}
\}^{\mathsf{I}} (u, v) = F (\bar{e}, u, v)$, and we set
$${\mh} (u, v) = \{ \bar{e} \}^{\mathsf{I}} (u, v).$$

We now claim that ${\mh}$ satisfies (i) and (ii) of the Theorem. Set $\rho = \rho^\mi$.
We have that\\

$ \begin{array}{rcll}
\mh (\gamma, \delta)& = &0 &\mbox{ if }\rho (\gamma) = 0\\

& = &1& \mbox{ if }\rho (\gamma) \neq 0 \wedge \rho
(\delta) = 0\\

& =&0& \mbox{ if }  \all \gamma_{\alpha} \ex \delta_{\beta} \mh
(\gamma_{\alpha}, \delta_{\beta}) \simeq 0\\

& =&1& \mbox{ if }\ex \gamma_{\alpha} \all \delta_{\beta} \mh
(\gamma_{\alpha}, \delta_{\beta}) \simeq 1.

\end{array}$
\\

\nod Which we can restate as:\\
\nod{\tmem{Claim 2}} {\em If  not both $\rho (\gamma) = \rho (\delta) = \omega_1$, then }${\mh}
(\gamma, \delta) \downarrow $  {\em and we have:}\\

(1) $\all \alpha < k(\gamma) \ex \beta <
k(\delta) $ ${\mh} (\gamma_{\alpha}, \delta_{\beta}) \simeq 0 \imp$ ${\mh} (\gamma, \delta) \simeq 0$ 

(2) $\ex \alpha < k(\gamma) \all \beta <
k(\delta) $ ${\mh} (\gamma_{\alpha}, \delta_{\beta}) \simeq 1 \imp$  ${\mh} (\gamma, \delta) \simeq 1$.\\

\nod \pf of {\em Claim 2} By induction on $\sigma
=_{\tmop{df}} \min \{ \rho (\gamma), \rho (\delta) \} < \omega_1$. If
either $\rho (\gamma), \rho (\delta) = 0$ the main claim is clear.
Suppose without loss of generality $0<\rho(\g)=\sigma<\omega_1$. Then $\all \g_\a  ( \a<k(\gamma) \imp \r(\g_\a)<\s)$ so by induction $\mh(\g_a	,\d_\b)\da$ for any $\a<k(\gamma),\b<k(\delta)$. The argument for $\delta$ is similar. (1) and (2) then follow from the definition of $\mh$ as  $\{ \bar{e} \}^{\mathsf{I}}$. \mbox{ }\hfill \qed (Claim 2)\\

\nod{\tmem{Claim 3}}
 ${\mh}$ {\em satisfies }{\em (i)  and  (ii) }.\\
 \pf Again by induction on $\sigma
=_{\tmop{df}} \min \{ \rho (\gamma), \rho (\delta) \} < \omega_1$. (If both
$\rho (\gamma) = \rho (\delta) = \omega_1$ there is nothing to do.) Thus at
least one of these ranks is countable.

If $\sigma = 0$: then either $\rho (\gamma) = 0$ and ${\mh} (\gamma, \delta) = 0$ or
$0 = \rho (\delta) < \rho (\gamma)$ and so ${\mh} (\gamma, \delta) = 1$.

\

If $0 < \sigma$: as inductive hypothesis we assume {\em (i)} and {\em (ii)} hold for all
$\bar{\gamma}, \bar{\delta}$ with $\min \{ \rho(\bar{\gamma}), \rho(\bar{\delta} ) \} <
\sigma$.
Suppose first that $\rho (\gamma) < \omega_1 \wedge \rho (\gamma) \leq \rho
(\delta)$. Then $\all \alpha < k(\gamma) \: \all \beta < k(\delta)$\,$\min \{
\rho (\gamma_{\alpha}), \rho (\delta_{\beta}) \} < \sigma$, and the inductive hypothesis for {\em (i)
} and {\em (ii)} applies and we have:\\

(3) $\rho^I (\gamma_{\alpha}) \leq \rho^I (\delta_{\beta}) \imp \mathsf{H}
(\gamma_{\alpha}, \delta_{\beta}) \simeq 0 ;$

(4) $\rho^I (\delta_{\beta}) < \rho^I (\gamma_{\alpha}) \imp \mathsf{H}
(\gamma_{\alpha}, \delta_{\beta}) \simeq 1.$\\

As these antecedents are mutually exclusive, we have that $\all \alpha <
k(\gamma) \: \all \beta < k(\delta) $ ${\mh} (\gamma_{\alpha}, \delta_{\beta})\downarrow.$ As $\rho (\gamma) \leq \rho (\delta)$ by assumption, (1) and (3)
imply ${\mh} (\gamma, \delta) = 0$.

Suppose secondly, that $\rho (\delta) < \rho (\gamma) \leq \omega_1$. Then
$\all \beta < k(\delta) \rho (\delta_{\beta}) < \rho (\delta) = \sigma$, and
by the inductive hypothesis again on {$ (i)$} and {$(ii)$} we have that for any $\alpha < k(\gamma),
\beta < k(\delta)$ (3) and (4) again hold. As $\rho (\delta) < \rho (\gamma)$
there is $\alpha < k(\gamma)$ with $\rho (\delta) \leq \rho
(\gamma_{\alpha})$. This implies $\all \beta < k(\delta) \rho (\delta_{\beta})
< \rho (\gamma_{\alpha}) .$ By (2) and (4) again ${\mh} (\gamma, \delta) = 1.$
{\qed} Claim 3 and Theorem\\

\begin{theorem}[Gandy Selection]\label{GandySelection} There exists a functional
  {$\mathsf{{Sel^I}}$} partial recursive in ${\mi}$ such that for all $e,
  m, x$ the following are equivalent:
  
  (i) $\ex p \in \omega . \{ e \} ^{\mi }(p, m, x)$ ;
  
  (ii) $ \{ e \}^\mi  (\mathsf{\tmop{Sel^I}}  (e, m, x), m, x)$.
\end{theorem}

{\pf}There is a p.r. function $G$ so that for any index $e$, $G (e)  = e^+$ is
an index so that
$$\{ e^+ \}^\mi  (p, m, x) \simeq \{ e \} ^\mi (p + 1, m, x)\mbox{ for }p \in \omega.$$ 
We
take $\mh(u, v)$ from the Stage Comparison Theorem \ref{OrdComp}, and define
$F (f, p, m, x)$ to be the functional which accords with the following
procedure.

Step 1 It computes $\mh( \la e, 0, m, x \ra, \pa{f, e^+, m, x})$.

Step 2 If this value is $0$, then that is the final value for $F (f, p, m,
x)$. If this value is $1$, it then computes $\mh( \la e, 0, m, x \ra, \la e,
\{ f \}^{\mi}(e^+, m, x) + 1, m, x \ra )$.

Step 3 If this latter value is 0, then that is the final value. If this value
is $1$, it then computes $\{ f \}^{\mi}  (e^+, m, x) + 1$ as its final value.

\

By the Recursion Theorem there is $\bar{f}$ so that $\{ \overline{f} \} ^{\mi} (p,
m, x) \simeq F (\bar{f}, p, m, x)$. We take $\mathsf{{Sel}^{\mi}}  = \{
\overline{f} \}^{\mi} $. We now show that this works.

Let $p (e, m, x)$ be the least $p$ such that $\{ e \}^{\mi}(p, m, x)\da$ if there is
such a $p.$\quad We show the equivalence of ($i$) with $(i i)$ by induction on
$\nobracket p (e \nobracket, m, x)$. $(i i) \Imp (i)$ is immediate.

\nod {\em Case 1} \ $p (e, m, x) = 0$.

Then $\{ e \}  (0, m, x)\da$, and so $\pa{e, 0, m} \in U _x$. Hence for some $i <
2$
$$\mh( \la e, 0, m, x \ra, \pa{\bar{f}, e^+, m, x} ) =
i.$$

\nod If $i = 0$: then $\tmop{Sel}^{\mi}  (e, m, x) = 0$ as required.

\nod If $i = 1$: this implies $| (\nobracket \bar{f}, e^+, m, x) | \nobracket < |
(e, 0, m, x) |$ which in turn means that there is a $q$ with $\{ \bar{f} \}^{\mi}
(e^+, m, x) = q$. As $\pa{e, 0, m} \in U _x^\mi$ we have that for some $j < 2$,

$$\mh\left( \la e, 0, m, x \ra, \pa{e , q + 1, m, x} \right) = j.$$

Again if $j = 0$ then $\tmop{Sel}^{\mi}(e, m, x) = 0$ as required.

If $j = 1$ then we have $| e , q + 1, m, x | < | e, 0, m, x | \nobracket $and
thence $\{ e \} ^{\mi} (q + 1, m, x)$. However then $\tmop{Sel}^{\mi}  (e, m, x) \simeq \{
\bar{f} \}^{\mi}  (e^+, m, x) + 1 = q + 1$, as we wanted.

\

\nod {\em Case 2} Now suppose that $p (e, m, x) > 0$. By our definitions $p (e, m, x) = p (e^+,
m, x) + 1$. So we can apply the induction hypothesis to $p (e^+, m, x)$ and
then there will be $q$ such that
$$\{ \bar{f} \} ^\mi (e^+, m, x) \simeq \mathsf{Sel^I}  (e^+, m, x) = q 
\mbox{ with }\{ e \}  (q + 1, m, x) \simeq \{ e^+ \}  (q, m, x)$$
\noindent and with both the latter defined. In particular $\la \bar{f}, e^+, m, x \ra \in U
_x$. By our assumption on $p$ we have $\pa{e, 0, m, x} \notin U _x$, and thus
\ $| (\bar{f}, e^+, m, x) | < | (e, 0, m, x) | |$ and thus
$$\mh\left( \la e, 0, m, x \ra, \pa{e , q + 1, m, x} \right) = 1.$$
Hence $\mathsf{Sel^I}  (e, m, x) = \{ \bar{f} \}^\mi  (e^+, m, x) + 1 = q + 1$
again as required for $\mathsf{Sel}^\mi $. \\ \mbox{ }{\qed}

\

With the Selection Theorem established we now can use it to get a number of
results about the structure of partial recursive functions and semi-recursive
sets. The first is immediate from this theorem.

\begin{lemma}
  \label{existclosure}For any relation {${\mr}$} that is semi-recursive
  in $\mathsf{I},$ for any $m, x$ we have:
  $$\ex q {\mr}(q, m, x) \equi {\mr}(\mathsf{Sel} ^\mi (m, x), m, x).$$
\end{lemma}
It might be  tempting to claim that the union of two semi-recursive
in  relations, $\mr$ and $\mathsf{S}$ say, is semi-recursive in $\mi$, is
established by running a procedure $P$ that simulates both functions $\{ e \}^\mi
$ and $\{ f \}^\mi $ simultaneously whose domains are $\mr$ and $\mathsf{S} $
respectively until, if possible, an $\pa{m, x}$ falls into one domain or the
other (or neither). But this will not work since we may have $\pa{m, x} \in
\tmop{dom} (\{ e \}^\mi ) \back \tmop{dom} (\{ f \}^\mi )$, but the computation $\{ f
\}^\mi  (m, x)$ being simulated is not just divergent, but has $\mathfrak{T } ^\mi (m,
x)$ illfounded. And moreover it is divergent before $\{ e
\}^\mi  (m, x)$ is convergent,  This would render $P$ undefined. The use of the Selection
Theorem neatly gets around this:

\begin{lemma}
  \label{closure}
  
  The class of relations semi-recursive in $\mi$ is closed under: \\(a) finite
  unions;\\ (b) bounded existential, and existential number quantification `$\ex
  n$' ;\\ (c) definitions by cases.
\end{lemma}

{\pf}(a) Suppose $\mr = \tmop{dom} (\{ e \}^\mi )$ and $\mathsf{S} =
\tmop{dom} (\{ f \} ^\mi)$ are two semi-recursive relations on $\omega \times
\omega^{\omega}$.
Let $F$ be the (ordinary) Turing function defined by $F (p) = e$ if $p = 0$
and $F (p) = f$ if $p > 0$. Set $\mh(p, m, x) = \{ F (p) \} ^\mi(m, x)$. Then $\pa{m,
x} \in \mr \cup \mathsf{S}$ iff $\ex p \mh(p, m, x)$. If $\{ h \} =\mh$, then $\ex p
\mh(p, m, x) \equi \mh(\tmop{Sel}^\mi  (h, m, x), m, x)$. Let $\{ s \} ^\mi (m, x) \simeq
\tmop{Sel}^\mi  (h, m, x)$. Then $\mr \cup \mathsf{S} = \tmop{dom} (\{ s \}^\mi)$. 

(b) This follows from Lemma \ref{existclosure}. (c) Exercise. \qed

\begin{lemma}
  A relation $\mr$ is recursive in ${\mi}$ iff both $\mr$ and $\neg \mr$ are
  semi-recursive in $\mi$.
\end{lemma}

{\pf}$\left( \Imp \right)$: If $\mr$ is recursive in ${\mi}$, so is $\neg \mr$, and
both are semi-recursive in $\mi$.

$(\Leftarrow)$: As both $\mr$ and $\neg \mr$ are the domains of some partial
recursive in ${\mi}$ functions, we can choose indices $a, b$, modifying those partial
recursive functions if need be, so that
$$
\begin{array}{rcl}
\mr (m, x) & \equi & \{ a \} ^{\mi} (m, x)\da 1\\
\neg \mr (m, x) & \equi & \{ b \} ^{\mi} (m, x)\da 0.
\end{array}
$$
Let $f$ be recursive with $f (0) = a$, $f (1) = b$ and $f (k)\ua$ for $k > 1$.
Let $\mg(i, m, x) = \{ f (i) \} ^\mi (m, x)$. Suppose that $\mg = \{ g \}^{\mi} $. Then for
all $(m, x)$ ${\mk}(m, x)=_{df}$ \\$\mg (\mathsf{Sel^I}  (g, m, x), m, x) \nobracket)$ is defined,
and then {\mk} is the total recursive characteristic function of $\mr$. Hence the
latter is recursive in $\mi$. {\qed}

\begin{lemma}\label{L3.28}
  For any partial functional,
  
  (i) ${\mf}$ is partial recursive in  $\mi$ iff $\tmop{Gr({\mf})}$ is semi-recursive in $\mi$.
  
  (ii) ${\mf}$ is recursive in $\mi$  iff $\mf$ is total and $\tmop{Gr} ({\mf})$ is recursive in $\mi$.
\end{lemma}

{\pf} Straightforward.{\qed}

\begin{lemma}
For any $B\sset \omega$, $B$ is semi-recursive in $\mi$ if and only if $B=\mathop{Im}(f)$ for some  partial function $f$ partial recursive in $\mi$.
\end{lemma}
\pf $(\Imp)$ Exercise (\cf Lemma \ref{reimage}.)
$(\Leftarrow)$: by Lemmas \ref{closure} and \ref{L3.28}.
\qed\\

Let $\Gamma^{\mathsf{I}}$ be the pointclass of semi-recursive in $\mathsf{I}$
relations. We have a norm $\rho^{\mathsf{I}}$ on the universal semi-recursive
in $\mathsf{I}$ set $\mathsf{U^I} $ and hence on any $A \in
\Gamma^{\mathsf{I}}$. We have that $\rho^{\mathsf{I}}$ is a
$\Gamma^{\mathsf{I}}$-norm in the sense of Moschovakis, {\cite{Mosch4}}p.153, or
\cite{KeMo72}, Sect 2.1, in that there is a semi-recursive in $\mathsf{I}$ relation
$\leq^{\rho}_{\Gamma}$ and a co-semirecursive in $\mathsf{I}$ relation
$\leq^{\rho}_{\neg \Gamma}$ with the property $(\ast)$ that:
$$ \rho (y) < \omega_1 \imp \all x \left\{ [\rho (x) < \omega_1 \wedge
\rho (x) \leq \rho (y)] \equi x \leq^{\rho}_{\Gamma} y \equi x
\leq^{\rho}_{\neg \Gamma} y \right\}.$$

\begin{definition}
  A pointclass $\Gamma$ is said to be {\tmem{normed}} or to have the
  {\tmem{Pre-wellordering Property}} if every pointset $A \in \Gamma$ admits a
  $\Gamma$-norm.
\end{definition}

\begin{lemma} \label{PWO}
  For any functional $\mi$, the class $\Gamma$ of relations semi-recursive in  $\mi$ has
  the prewellordering property.
\end{lemma}

{\pf} It suffices to show this for the universal semi-recursive in $\mi$  set
$U^\mathsf{I}$.  We verify $(\ast)$ above for $\rho = \rho^\mi $ and $\Gamma =
\Gamma ^\mi$. Using the Stage Comparison Theorem \ref{OrdComp} we can take 
  $$
\begin{array}{rcl}
x\leq^{\rho}_{\Gamma} y &\Equi & \mh(x, y) = 0 \mbox{ and }\\
x \leq^{\rho}_{\neg \Gamma} y & \Equi &\mh(x, y) \neq 1.
\end{array}$$
\nod as {\mh} is partial recursive in $\mi$ , we are done. {\qed}

\begin{lemma}\label{Spectorpointclass}
  For any functional $\mi$ the semi-recursive in $\mi$  relations form a Spector class.
\end{lemma}

{\pf} To form a Spector class the relations must be closed under $\wedge,
\vee$, bounded and unbounded existential and universal number quantification,
$\ex^{\omega}_{\leq}, \ex^{\omega}$, $\all^{\omega}_{\leq}, \all^{\omega} .$ It
must have $\omega$-parametized universal classes in each arity, be closed under recursive in $\mi$ substitutions, and have the
Prewellordering Property. These have now all been established above. {\qed}

\begin{corollary}
  The class of sets semi-recursive in $\mathsf{I}$ has the Reduction property,
  and its dual class the Separation Property.
\end{corollary}

From this many structural results follow, in a straightforward manner.
We mention one of these here without proof as we shall not be using it directly in the sequel.

\begin{definition}
  Let $\Delta^{\mathsf{I}} (x)$ those sets that are both
  ittm-semi-recursive-in-$x$ and $\mathsf{I}$ and
  co-ittm-semi-recursive-in-$x$ and $\mathsf{I}$. When $\mathsf{I} =
  \mathsf{E}$ we omit it.
\end{definition}

\begin{lemma}
  ({\cf}{\tmem{Kleene {\cite{Kleene59b}}}} If $Q \sset \bai \times \bai$ is 
  (ittm) -semi-recursive in $\mathsf{I}$, then defining
  
  \ \ \ \ \ \ \ \ \ \ \ \ \ \ \ \ \ \ \ \ \ \ \ \ \ \ \ \ \ \ \ $P (x) \equi
  \ex y \in \Delta \, Q (x, y)$,
  
  {\nod}we have that $P$ is  (ittm) -semi-recursive in $\mathsf{I}$. By
  relativisation, if $Q \sset \bai \times \bai \times \bai$ is  (ittm)
  -semi-recursive in $\mathsf{I}$, then defining \
  
  \ \ \ \ \ \ \ \ \ \ \ \ \ \ \ \ \ \ \ \ \ \ \ \ \ \ \ $P (x, z) \equi \ex
  y \in \Delta (z) \, Q (x, z, y)$,
  
  {\nod}then $P$ is  (ittm) -semi-recursive in $\mathsf{I}$.
\end{lemma}

We now investigate the norm given by the rank function $\rho^I$, and compute its length in terms of the wellorderings ittm recursive in $\mi$.

\begin{definition}
  We let  $| - |_0^{\mathsf{I}} \nobracket$  be the regular norm on
  $\omega^2$ induced by restricting the prewellordering $\leq^\rho$ derived from $\rho$, on
  $U^{\mathsf{I}}$ to sequences of type $\left( e, m, \overline{0} \right)$ ($
  \overline{0}$ is the zero function); thus 
  $$| e, m |^{\mathsf{I}}_0 =
  \sup \,\!\! ^+  \{ | g, n |^{\mathsf{I}}_0 : \rho(g, n, \bar{0}) < \rho(e, m,
  \bar{0}) \} \nobracket  .$$
  
  \nod Let $\kappa^{\mathsf{I}} = \sup^+ \{ | e, m |_0^{\mathsf{I}} :
  U^{\mathsf{I}} (e, m) \} \nobracket$. We slightly abbreviate $|f,n|_0^\mi$ as $|f,n|^\mi$, or as $|f,n|$ when the context is clear.
\end{definition}

\begin{definition}\label{alphazero}
  For $\mathsf{I}$ a functional, let $\alpha^{\mathsf{I}}_0$ be the least
  ordinal not the order type of a wellordering of $\omega$ which is recursive
  in $\mathsf{I}$.  If $\mi = \me$ we may simply write $\alpha_0$.
\end{definition}

For $\mathsf{I}$ then, $\alpha^{\mathsf{I}}_0$ is the analogue of
$\omega_1^{\mathsf{I}, \tmop{ck}}$.

\begin{theorem}\label{kappaequalsalpha}
  $\kappa^{\mathsf{I}} = \alpha_0^{\mathsf{I}}$.
\end{theorem}
\pf
$\kappa^{\mathsf{I}} \leq \alpha_0^{\mathsf{I}}$: following \cite{Hi78}, note that for $(e,m)\in U^\mi$ the sets
$$  \{ (f,n)\in U^\mi\mid |f,n|^\mi < |e,m|^\mi\} \mbox{ and }  \{ (f,n)\in U^\mi\mid |f,n|^\mi = |e,m|^\mi\}
$$
 are both  generalised itttm recursive in $\mi$, by using the prewellordering property holding for semi-recursive in $\mi$ sets. The first set is a prewellordering but we refine it to a wellordering by choosing unique notations for each ordinal less than $|e,m|^\mi$ and placing them in the set $A_{(e,m)}$:
$$ A_{(e,m)}\dfs \{ |f,n|\in U^\mi \mid |f,n|^\mi < |e,m|^\mi\wedge |g,p|^\mi = |f,n|^\mi \imp (f,n)\leq_{lex} (g,p)\}.$$
Then define the ittm recursive in $\mi$ wellorder of type $|e,m|^\mi$: $$ R_{(e,m)}((f,n),(g,p))\edfs (g,p), (f,n) \in A_{(e,m)}\wedge (f,n)\leq_{\Gamma^\mi}^{\rho^\mi}(g,p).$$

 To show that $\kappa^{\mathsf{I}} \geq \alpha_0^{\mathsf{I}}$, let $\d <\alpha_0^{\mathsf{I}}$ and let $y\in WO \wedge ||y||=\d$, with $y$ generalised recursive in $\mi$.  Let $|n|_y$ denote the order type of $y\rest n$. Suppose also that $n_0$ was such that $|n_0|_y=0$. Define two functions partial recursive in $\mi$.\\

$
\begin{array}{rll}
& n&\mbox{  if }n<_y m; \\
F(n,m) = & & \\
&  n_0 &\mbox{ otherwise; }\\

& & \\
\end{array}
$

And:

$
\begin{array}{rll}
& 0 & \mbox{ if } j=n_0 \\
G(e,j) = & & \\
& \mi(\l n.\{e\}^\mi(F(n,j))),& \mbox{ otherwise.}

\end{array}
$

\

Using the $\mi$-Recursion Theorem there exists an index $\bar e$ so that $G(\bar e,j)\simeq \{\bar e\}^\mi(j)$. We now show by induction on $|j|_y$  that $(\ast) :\pa{\bar e,j}\in U^\mi \wedge |\bar e, j|^\mi \geq |j|_y$ holds, thereby establishong that there is an $\mi$-prewellordering of length at least $\d$. 

If $|j|_y=0$, then $j=n_0$ and so $\{\bar e\}^\mi(j)=G(\bar e,n_0) =0$.  Hence $(\ast)$ trivially holds in this case. So suppose $|j|_y>0$. Note that for any $n$, $|F(n,j)|_y <|j|_y$, and so by the induction hypothesis $\pa{\bar e, F(n,j)}\in U^\mi$. This implies that $\l n.\{e\}^\mi(F(n,j))$ is total, and so $G(\bar e,j)$ is defined, and finally $\pa{\bar e,j}\in U^\mi$, which is the first half of $(\ast)$.

From the definition of $G$ let $a$ be a programme index for $G$ derived from the above definition; thus with $|a, e,j|^\mi \geq sup^+\{| e, F(n,j)|^\mi\mid n < \w\}$. Thus taking $\bar e$ as $e$ we shall have $|a, \bar e,F(n,j)|^\mi < |a,\bar e, j|^\mi$. By the definition of the fixed point $\bar e$ in the Recursion Theorem we shall further have $|a,\bar e, j|^\mi < |\bar e, j|^\mi$. By the induction hypothesis then, for any $n<_yj$,
$|n|_y \leq |\bar e, n|^\mi =|\bar e, F(n.j)|^\mi < |\bar e, j|^\mi$. But this implies that $|j|_y \leq |\bar e, j|^\mi$ which is the second half of $(\ast)$, as required.
\qed\\

In particular the following question is pertinent for the simplest
functional:

\

{\tmem{Question: What is}} $\alpha_0 \dfs \alpha_0^{\mathsf{E}}$?

\

By our constructions to date, we have seen that it must be much larger than
$\zeta$ (see Ex.\ref{Ex2.8}) . It is part of our task to identify this
ordinal. The next lemma is straightforward.

\begin{lemma}
  Let $x$ be ittm-recursive in $\mathsf{I}$. Then $\zeta^x <
  \alpha^{\mathsf{I}}_0$. Conversely for any $\tau < \alpha^{\mathsf{I}}_0$
  there is $y$ ittm-recursive in $\mathsf{I}$ with $\tau < \zeta^y$. \ Thus
  $\alpha^{\mathsf{I}}_0$ is a limit of extendibles.
\end{lemma}

We state the next two lemmata for completeness: they express properties that would be expected of universal sets and their prewellorderings;  we do not use them in the sequel and so omit their proofs here.

\begin{lemma}
  (Boundedness Theorem)  If $A$ is co-semi-recursive in $\mathsf{I}$ and $A
  \sset U^{\mathsf{I}}$ then
  
  \ \ \ \ \ \ \ \ \ \ $\sup^+ \{ | u |^{\mathsf{I}} : u \in A \} <
  \kappa^{\mathsf{I}} \nobracket \nobracket$.
\end{lemma}

\begin{lemma}
  (Hierarchy Theorem)  Let $U^{\mathsf{I}}_{\tau} \dfs \{ u : u \in
  U^{\mathsf{I}} \wedge | u |^{\mathsf{I}} < \tau \} \nobracket \nobracket$.
  For any $R \sset^k \omega$,
  
  $R$ is recursive in $\mathsf{I}$ iff $\ex \tau < \kappa^{\mathsf{I}}  (R
  \nobracket$ is many-one reducible to $\nobracket U^{\mathsf{I}}_{\tau})$
  
  \ \ \ \ \ \ \ \ \ \ \ \ \ \ \ \ \ \ \ \ \ \ \ \ \ \ \ \ iff \ $\ex \tau <
  \kappa^{\mathsf{I}}  \left( R = \left\{ \tmmathbf{m} : \pa{e, \tmmathbf{m}}
  \in U^{\mathsf{I}}_{\tau} \right\} \right)$
\end{lemma}

\section{Computation Lengths}

As with the basic ittm's the strength of the model is tied up with the length
of computations possible or needed on, or by, the model; in turn that is
mutually tied up with the class of reals so recursive (the slogan ``clockables
are writables'' is apposite). The same considerations are true of the
generalised ittm recursions here. Even for recursions in $\mathsf{E}$ we see
(Lemma \ref{Ex2.8}) there are computations recursive in $\mathsf{E}$ that
on input $x$ compute (a code for) $L_{\Sigma^x} [x]$, or for $L_{\zeta^x}
[x]$. Combining this with a programme like the Theory Machine we shall see
that much longer sections of the $L$-hierarchy can be computed on integer
input, indeed done in the simplest fashion such a programme would loop at the
first extendible ordinal that is a limit of such. \ So we first investigate
this hierarchy and then see how to compute through long initial segments.

\subsection{Extendability Hierarchy}
For $E$ a class of ordinals, let $E^{\ast}$ denote the class of its limit
points.

\begin{definition}
  We classify $\Sigma_2$-extendible ordinals as follows. Define by recursion
  on $\alpha$ the class $E^{\alpha}$ of $\alpha ($-$\Sigma_2)$-extendible
  ordinals:
  \[ \begin{array}{lcl}
       E^0 & = & \{ \zx{0} \mid \hspace{0.17em} \zx{0} \text{ is
       $\Sigma_2$-extendible} \} ;\\
       E^{\alpha + 1} & = & \{ \zx{\alpha + 1} \mid \hspace{0.17em} \zx{\alpha
       + 1} \in (E^{\alpha})^{\ast} \cap E^0 \} ;\\
       E^{\lambda} & = & \bigcap_{\alpha < \lambda} E^{\alpha} \tmop{for} \tmop{Lim} (\lambda) .
     \end{array} \]
  {\tmem{}}{\tmem{{\tmem{(For $\alpha = 0$ we may simply say
  ``$\Sigma_2$-extendible'' or even just ``extendible''.)}}}}
\end{definition}

Note that $\alpha < \beta \imp E_{\alpha} \supset E_{\beta}$. Here we
decorate the variable $\zeta$ with the prefix indicating its level of
extendibility. We shall let $\sx{k}$ indicate that for some $\zx{k}$,
({\zx{k}}, {\sx{k}}) is a $k$-extendible pair (we shall be mainly interested
in finite $\alpha = k$). \ Obviously for any $\gamma$ the least element of
$E^k$ greater than $\gamma$ is always an element of $E^k \setminus E^{k + 1}$,
{\ie} is $k$-extendible but not $k + 1$-extendible. We extend the definition,
relativising to reals $x$ the notion of $\alpha$-$x$-extendibles, $E^{\alpha}
[x]$, where an $x$-extendible pair is a $\xi, \sigma$ with $L_{\xi} [x]
\prec_{\Sigma_2} L_{\sigma} [x]$ {\etc}

Notice that for any {\zx{1}} if $L_{\sx{1}}$ is the natural $\Sigma_2$-end
extension of $L_{\zx{1}}$, then $\sx{1}$ is also in $(E^0)^{\ast}$ (but not
necessarily in $E^1$). This holds by simple $\Pi_2$-reflection of the
statement that there are arbitrarily large elements of $E^0$ below $\zx{1}$ up
to $\sx{1}$. So in fact there must be $\Sigma_2$-extendible pairs of the form
$\left( \zx{0}, \sx{0} \right)$ which are {\tmem{nested}} in the interval
$\left( \zx{1}, \sx{1} \right.$), that is $\zx{1} < \zx{0} < \sx{0} < \sx{1}$.
This is suggestive of the kind of linearized computation that allows one depth
of subroutine call, Thus the top node of its computation tree has rank only
$1$. However the computation continues sufficiently far that it only enters a
final loop at $\zx{1}$, all the while making subroutine calls to nodes
$\nu_{\alpha}$, (all at $\nobracket \Lambda = 1)$ for $\alpha$ unboundedly in
$\zx{1}$, and thence by reflection it must be doing so for $\alpha$
unboundedly in $\sx{1}$. But at $\sx{1}$ it drops back to $\zx{1}$. We thus have:

\begin{lemma}
  Let $P_e^{\mathsf{E}} (j)$ be such that no query $Q^{\mathsf{E}} (e_1, m_1,
  y_1)$ is made so that $\{ e_1 \}^{\mathsf{E}} (m_1, y_1)$ itself invokes a
  query (in other words the rank of $\mathfrak{T} (e, j) =\rho^\me(e,m) \leq
   1$). Then if
  $\zx{1}$ is the least element of $E^1$, \ $P_e^{\mathsf{E}} (j)$ enters a
  final loop by stage $\zx{1}$.
\end{lemma}

This picture propagates: if queries are made to a greater depth in
calculations of $P_e^{\mathsf{E}} (n)$ a greater rank of
$\Sigma_2$-extendibles may be needed to represent the ordinal length of time
for the overall computation, with in turn, the rank corresponding to the rank
of the computation tree $\mathfrak{T} (e, n)$. (Or equivalently the rank
$\rho^\me$.)

\begin{lemma}
  Suppose $P_e^{\mathsf{E}} (j) \da z$. Let $\rho = \rho^{\mathsf{E}} \left(
  \pa{\pa{e, j}, z} \right)$. If $\zx{\rho}$ is the least element of
  $E^{\rho}$, then \ $P_e^{\mathsf{E}} (j)$ enters a final loop by at latest
  stage $\zx{\rho}$. More generally, if $P_e^{\mathsf{E}} (j, x) \da z$, $\rho
  = \rho^{\mathsf{E}} \left( \pa{\pa{e, j, x}, z} \right)$ and $\zx{\rho}^x$
  is the least element of $E^{\rho} [x]$, then $P_e^{\mathsf{E}} (j, x)$
  enters a final loop by stage $\zx{\rho}^x$.
\end{lemma}

{\pf}By induction on $\rho$. Suppose the Lemma holds for any $\rho' $of the
form $\rho^{\mathsf{E}} \left( \pa{\pa{e', j'}, z'} \right) < \rho =
\rho^{\mathsf{E}} \left( \pa{\pa{e, j}, z} \right)$ as specified. Let $\Lambda(\tau)$ abbreviate $\Lambda({\me}, e,j,\tau)$.

{\tmem{Claim:}} $\Lambda (\nobracket \left. \zx{\rho} \right) = \Lambda \left(
\sx{\rho} \right) = 0$.

{\pf} Suppose not, and $\Lambda (\nobracket \left. \zx{\rho} \right) = \Lambda
\left( \sx{\rho} \right) = k > 0$. As $\Lambda \left( \zx{\rho} \right) > 0$
there is some `current query' $Q (e_i, m_i, x_i) $ in process at stage
{\zx{{\rho}}}. We note that this subcomputation could not have been called at
stage $\zx{\rho}$ itself, but must have been called at an earlier stage
$\alpha_0 < \zx{\rho}$. (For otherwise, $\zx{\rho}$ would be
$\Sigma_2$-definable in $L_{\sx{\rho}}$, as the last point at which a
subcomputation at this level was invoked - because if the query $Q (e_i, m_i,
x_i)$ was completed, control would have passed up to level $k - 1$. But then
by $\Sigma_1$ reflection that would happen unboundedly in $\zx{\rho}$, and
thence \ $\Lambda \left( \zx{\rho} \right) < k$ - a contradiction.) \ By
similar reasoning we have a property $(\ast)$: that this subcomputation is run at levels $k' \geq
k$ at all stages $\beta \in \left[ \alpha_0, \sx{\rho} \right)$ - again it
could not be completed at level $k$ at a stage $\tau < \sx{\rho}$, as then
control would pass to level $k - 1$, and result in a contradiction once more.

We next note that if $\{ e_i \} (m_i, x_i) \da w_i$, then $\rho' \assign
\rho^{\mathsf{E}} \left( \pa{\pa{e_i, m_i, x_i}, w_i} \right) < \rho$ (being a
subcomputation of $P_e^{\mathsf{E}} (j)$). However $x_i \in L_{\zx{\rho}}$ (as
the query $Q (e_i, m_i, x_i) $was invoked at stage $\alpha_0 < \zx{\rho}$).
Moreover $\zx{\rho}$ is a limit of $x_i$-$\zx{\rho'}$-extendibles. Let $(\xi,
\sigma)$ be such a pair, with $\alpha_0 < \xi < \sigma < \zx{\rho}$. By the
(more general case of the) inductive hypothesis $\{ e_i \} (m_i, x_i)$ has
$(\xi, \sigma)$ as a looping pair, and consequently the query $Q (e_i, m_i,
x_i)$ is completed by stage $\sigma$. But then $\Lambda (\sigma + 1) = k - 1$,
contradicting $(\ast)$ that the level is  constantly $\geq k$ in $\left[ \alpha_0, \sx{\rho} \right)$.
Contradiction!{\qed \,Claim}

\

Hence the {\tmem{Claim}} holds: but this is the Lemma: the snapshots of the
master computation $P_e^{\mathsf{E}} (j)$ at times $\zx{\rho}$ and  ${\sx{{\rho}}}$ are identical.{\qed} Lemma\\

\begin{corollary}
With the notation of the last lemma, if $\zx{\rho}^x$
  is the least element of $E^{\rho} [x]$, and suppose it has $\sx{\rho}^x$ as its least  $\Sigma_2$ extension,
 then $P_e^{\mathsf{E}} (k, x)$
  has $H(\mathsf{E}, e,k,x   ) \leq \sx{\rho}^x$.

\end{corollary}
\pf The Lemma said that the procedure entered a final loop by stage no later than $\zx{\rho}^x$, and so the overall length of the computation is no greater than $\sx{\rho}^x$. \qed\\

The methods above always allow us to calculate the $k$-extendibles (and even
$\alpha$-extendibles) above any ordinal, as we shall now turn to.

We first collect together some of the above Facts and results, in order to
abbreviate our descriptions of algorithms. This will help to have a library of
basic algorithms which we shall simply quote as being `recursive in
${\mathsf{K}}$' without further justification. (We just use the
adjective `basic' to classify them; we are not intending that they form a
basis for any class.)

\begin{definition}
  (Basic Computations-$\mathrm{BC}$) Let $\mathsf{K}$ be suitable.
  
  (i) Any standard ittm-computation $P_e (n, x)$ is Basic.
  
  (ii) If a code $y$ for an $\alpha$ ordinal is given, then the computations
  that compute: for any $x$ (a code for) $L_{\alpha} [x]$ and the satisfaction
  relation for $L_{\alpha} [x]$ are both Basic (in $x, y$). (These
  computations show those objects are $\mathsf{} \mathsf{K}$-recursive, if
  $\alpha$ is).
  
  The following functions are all $\mathsf{K}$-recursive, and Basic ({\cf}Lemma \ref{Ex2.8}):
  
  (iii)  $x \rightarrowtail \tilde{x}$; and thus $x
  \rightarrowtail T^2_{\zeta^x}$ 
  
  (iv) $x \rightarrowtail \zeta^x$, the least
  $x$-$\Sigma_2$-extendible;
  
  (v) $x \rightarrowtail \Sigma^x$, the larger of
  the next extendible pair in $x$;
  
  (vi)  $x \rightarrowtail \Sigma^{x+}$ (the next
  admissible beyond $\Sigma^x$).
\end{definition}

Stronger ordinals than simply $\Sigma^{x+}$ are
${\mathsf{K}}$-recursive.
The following is only an indication of what is possible (we shall see that there are much longer ordinals that are computable) but it illustrates further the method of book-keeping by computing successive  $L_\a$-theories.

\begin{lemma}
  \label{4.30} Let $\mk$ be suitable. There is a (Turing) recursive sequence of indices $\pa{e_i |i <
  \omega}$ so that for any $\alpha < \omega_1$ with a code $x \in \cant$,
  $P^{\mathsf{K}}_{e_i} (x)$, {\tmcolor{red}{ \tmcolor{black}{with $\tmop{rk}
  (\mathfrak{T}^{\mathsf{K}} (e_i, x)) = i,$}}} computes a code for
  the next $i$-$x$-extendible $\zx{i} > \alpha$.
\end{lemma}
{\pf} For $i = 0$ this has been done
using Basic Computations, and is Lemma \ref{Ex2.8}. Suppose $e_i$ has been defined, and we describe the
programme $P^{\mathsf{K}}_{e_{i + 1}}$ . Assume without loss of generality
that $\alpha = 0$, $x = \mathrm{const}_0$. Then $P^{\mathsf{K}}_{e_i} (0)$
computes a code for the least $i$-extendible, $\zeta_0 \assign \zx{i}$ say. By
a basic computation let a slice of the scratch tape $R$ be designated to hold
$T_{\zeta_0}^2$; $R \assign T^2_{\zeta_0}$. A code for $\zeta_0$,
$W_{\zeta_0}$ say, is recursive in $T^2_{\zeta_0}$. Now compute
$\tmcolor{red}{\tmcolor{black}{P^{\mathsf{K}}_{e_i}}} (W_{\zeta_0})$. This
yields the next $i$-extendible $\zeta_1 : = \hspace{0.17em}^i \zeta_1$. Now,
using Basic Computations, write successively to $R$ the theories
$T^2_{\zeta_0}, T_{\zeta_0 + 1}^2, \ldots, T^2_{\zeta_0 + \beta}, \ldots$ for
$\beta < \zeta_1$. We note that at limit stages $\lambda \leq \zeta_1$, $R$
will contain ``liminf'' theories $\hat{T}_{\lambda} = \mathrm{Liminf}_{\alpha
\rightarrow \lambda} T^2_{\alpha}$ (by the usual automatic ittm liminf
process) but that $T^2_{\lambda}$ is uniformly r.e. in \ $\hat{T}_{\lambda}$
as we saw above. And again a code $W_{\lambda}$ for $\lambda$ is then
arithmetic in $T^2_{\lambda}$ - uniformly in $\lambda$. The point of this
exercise of writing theories to $R$ is to ensure continuability of the
computation, and that we do not start to loop too early. (Another way to put
this is to say that it ensures sufficient `universality'.) The writing out of
all levels of the theories to $R$ is a precautionary step: in general we do
not have $\hat{T}_{\zx{i + 1}} = \liminf_{^i \zeta \rightarrow \zx{i + 1}} 
\hat{T}_{^i \zeta}$. However $\hat{T}_{\zx{i + 1}} ( = T^2_{\zx{i + 1}}
)$ is what we shall need to calculate $\zx{i + 1}$.

Set $R \assign$ $T^2_{\zeta_1} \in L_{\zeta_1 + 1}$; now compute
$P^{\mathsf{K}}_{e_i} (W_{\zeta_1})$ and repeat this process. \ As there is no
means for the process we are describing to halt, there is a least looping pair
of ordinals for it, $(\xi, \sigma)$ say. \ Let $(\zx{i + 1}, \sx{i + 1})$ be
the least $i + 1$-extendible pair. We claim that this is the pair $(\xi,
\sigma)$. By construction both $\sigma, \xi \in (E^i)^{\ast}$. Suppose $\xi <
\hspace{0.17em} \zx{i + 1}$. By the repetition of the contents of $R$ in the
loop points, we have $\hat{T}_{\xi} = \hat{T}_{\sigma}$ in the above
algorithm, hence $T^2_{\xi} = T^2_{\sigma}$, \ and thus:

\

{\tmem{Claim: }} $L_{\xi} \prec_{\Sigma_2} L_{\sigma}$.

\

{\pf}of Claim: Notice that $\xi$ is least with $T^2_{\xi} = T_{\sigma}^2$ and
$\xi \in (E^i)^{\ast}$ (otherwise there'd be an earlier beginning of our
loop). Then the $\Sigma_2$-skolem hull of $\emp$ in $L_{\sigma}$ is
$L_{\xi}$.{\qed} Claim

\

If the Claim holds then $\xi \in (E^i)^{\ast} \cap E^0 = E^{i + 1}$. This
contradicts the minimality of $\zx{i + 1}$. Hence $\xi$ equals the latter, and
$\sigma = \hspace{0.17em} \sx{i + 1}$ follows.

Hence we may compute $\hat{T}_{\zx{i + 1}} = T^2_{\zx{i + 1}}$, as the
eventually fixed output by means of the above looping procedure
\tmcolor{red}{\tmcolor{black}{(and determining this convergent output requires
the extra $+ 1$ depth to $i + 1$ to the overall calculation})}. It is thus
recursive in ${\mathsf{K}}$ (and $x$). We let $P^{\mathsf{K}}_{e_{i +
1}}$ be the programme of the procedure just described followed by the basic
computation that finds a code $W_{\zx{i + 1}}$ for $\zx{i + 1}$ by a method
uniformly arithmetic in the now $\mathsf{K}$-recursive $T^2_{\zx{i
+ 1}}$.

Finally note that the continuing description of the programme
$P^{\mathsf{K}}_{e_{i + 2}}$ from \ $P^{\mathsf{K}}_{e_{i + 1}}$ merely
repeats the above but altering only a few suffices. We may thus determine a
recursive function $i \mapsto e_{i + 1}$. \\ \mbox{ }\hfill \qed

\

Similar arguments to those of Lemma \ref{4.30} show:

\begin{lemma}
  \label{nextext}There is a recursive sequence of indices $\pa{e'_i \mid i <
  \omega}$ so that $P^{\mathsf{K}}_{e'_i} (\tmmathbf{m}, \tmmathbf{x})$ writes
  a code for $L_{^i \Sigma (\tmmathbf{x})} [\tmmathbf{x}]$, the least
  $\Sigma_2$-extension of $L_{^i \zeta} [\tmmathbf{x}]_{\nosymbol}$ where
  $\left( \zx{i}, \sx{i} \right)$ is the least $i$-$\tmmathbf{x}$-extendible
  pair in $E^i (\tmmathbf{x}) .$
\end{lemma}

The last lemmas shows only that we can recursively find, for example, the
least $\Sigma_2$-extendible in a real $x$, namely $\zeta^x$, or $\Sigma^x$, or
some \ {\zx{i}}. \ \  However more is possible: given \ $(e, m, x)$ we may,
recursively in $\mathsf{K}$, compute a code for $\sigma_0$ where $(\zeta_0,
\sigma_0)$ is the least looping pair of ordinals for the computation $\{ e
\}^{\mathsf{K}} (m, x)$ (assuming of course the latter has a wellfounded tree
$\mathfrak{T}^{\mathsf{K}} (e, m, x)$).

For this to work we,
as a minimum, need to be able to compute (recursively in terms of $\mathsf{K}, e,
m, x$) the ordinal $H (\mathsf{K}, e, m, x)$ which is the overall length of
this computation. To put it in terms often used for ordinary ittm's, the
`clockable' ordinal $H (\mathsf{K}, e, m, x)$ needs to be `writable', which is
what the next Lemma asserts. We thus regard the ancestor of the Lemma as being
the ``{\tmem{Clockables = Writables}}'' Lemma of the basic machines. We stick
with $\mathsf{E}$.

\begin{lemma}
  \label{H-lemma}There is a p.r. function $k$ so that if $\{ e \}^{\mathsf{E}}
  (m, x)$ has a wellfounded computation tree $\mathfrak{T}^{\mathsf{E}} (e, m,
  x)$ then $\{ k (e) \}^{\mathsf{E}} (m, x)$ computes a code for $H
  (\mathsf{E}, e, m, x)$.
\end{lemma}

{\rem}Note that we do not need $\{ e \}^{\mathsf{E}} (m, x) \da$ to assert
this. $H (\mathsf{E}, e, m, x)$ is defined even if $\{ e \}^{\mathsf{E}} (m,
x) \ua$ as long as $\mathfrak{T}^{\mathsf{E}} (e, m, x)$ is wellfounded.

\

{\pf}Given a successful computation we define another computation that sums
up all the lengths of the loops of all the subcomputations in the tree
$\mathfrak{T}^{\mathsf{E}} (e, m, x)$. Suppose then that
$\mathfrak{T}^{\mathsf{E}} (e, m, x)$ is wellfounded.

We outline a process for calculating codes for ordinals that are the length
of the computation in hand, and when we get to a query we again apply the
self-same process to that query to obtain the length of {\tmem{its}} loop. \
Our calculation in the end should output $W_{\sigma}$ a code for $H
(\mathsf{E}, e, m, x)$. Just as earlier we used $\Sigma_2$-theories and the
liminf property to keep track of levels of the $L$ hierarchy, (and so
ordinals) we do the same here, and so avoid the problem of registers designed to
contain written codes for wellorders,  containing garbage at limit
stages. \ We suppose that $(\zeta_0, \sigma_0) $ is the least looping pair of
ordinals for the computation $\{ e \}^{\mathsf{E}} (m, x)$. And thus $H
(\mathsf{E}, e, m, x) \geq \sigma_0$.

\

\noindent We describe then first an $\mathsf{E}$-recursive $F (f, e, m, x)$. 

\

a) It first simulates running $P_e^{\mathsf{E}} (m, x)$ on a scratch tape;
for each step $\xi$ we denote the current snapshot of the simulation
$P_e^{\mathsf{E}} (m, x)$ by $s_{\xi}$.

\

b) at each step $\xi$ of the simulation it also writes using another scratch
tape $T^2_{\delta (\xi)} [x]$ for some $\delta (\xi) \geq \xi$ (to be
specified below) to a register $\calr$. (We are going to drop the uniform $x$
in what follows and write here $T^2_{\delta (\xi)}$.)

\

c) If at a step of the simulation of $P_e^{\mathsf{E}} (m, x) $ a query
?$Q^{\mathsf{E}} (e_1, m_1, x_1) ?$ is encountered it replaces the query with
?$Q^{\mathsf{E}} \left( f, \la e_1, m_1 \ra, x_1 \right) ?$ and inserts
following it, the program lines as specified in $g_0$ of Lemma \ref{Ex2.6}. This has the effect that if the latter
computation, $\{ f \}^{\mathsf{E}} \left( \la e_1, m_1 \ra, x_1 \right)$, is
convergent to say $w_1$, this process will return as well as $\mathsf{E}
(w_1)$ as answer to the query, $w_1$ itself, to some reserved register on the
scratch tape $\mathcal{S}$ in the computation $F (f, e, m, x)$. We set:
\

\
$
\begin{array}{rcl}
 \eta(e_1, m_1, x_1) &=& 0, \mbox{ if $w_1$ is not the code of an
ordinal;} \\
 &=& \| w_1 \| \mbox{ otherwise.}

\end{array}
$

\

\noindent Then $F (f, e, m, x)$ continues:\\

d) {\em For a non-query instruction at stage $\xi$ in }$P_e^{\mathsf{E}} (m, x)$:

(i) Let $\bar{\delta} (\xi) \dfs \sup  \{ \delta (\xi') \mid \xi' < \xi \}$;

($R$ will contain the theory $\widehat{T }_{\bar{\delta} (\xi)}$ where
$\bar{\delta} (\xi) \dfs \sup \{ \delta (\xi') \mid \xi' < \xi \}$. A code for
$\bar{\delta} (\xi)$ is then uniformly arithmetic in $\widehat{T
}_{\bar{\delta} (\xi)}$.\quad If $\xi = \bar{\xi} + 1$, then again we
similarly have a code for \ $\bar{\delta} (\xi)$.)

Set $\delta (\xi) = \bar{\delta} (\xi) + 1$; then:

(ii) writes  the theory $T^2_{\delta (\xi)}$ to $\calr$.

\

{\em For a query instruction $?Q^\mathsf{E}(e_1, m_1, x_1   )?$ at stage $\xi$,}  then:

(i)$_q$ it sets $\delta (\xi) = \bar{\delta} (\xi) + \eta (e_1, m_1, x_1)$;

(ii)$_q$ writes successively the theories $T^2_{\alpha}$ to $\calr$, for
$\alpha \in (\bar{\delta} (\xi), \delta (\xi)] .$

\

e) It then checks whether $s_{\xi}$ was the least final looping snapshot of
$P^{\mathsf{E}}_e (m, x)$. This is effected by calling a generalised query as
we saw above in the proof  of Lemma \ref{L2.15} (ii) using the program there
$P_{t_0} (e, m, x, s_{\xi})$.

\

IF $s_{\xi}$ was NOT a final looping snapshot, then a new step in the simulation is taken,
and RETURNS to b), resulting in a new theory $T^2_{\delta (\xi) + 1}$ being
written to $\calr$, and, taking cognizance of c), proceeds to d) and the next
snapshot $s_{\xi + 1}$ is then rechecked at e) {\etc}

IF $s_{\xi}$ was a final looping snapshot:  then, the first (and only) time this is encountered is when  $\xi = \xi_0 < \sigma_0$ with $s_{\xi} = s_{\sigma_0}$. We set
$s_{\xi}$ to one side.
We have that, arithmetic from $T^2_{\delta (\xi)} [x]$ we have
a wellorder of type $\xi_0$. We need one of type $\sigma_0$.  $F$ sets in
motion queries to ascertain the set of cells of the top level tape of eventually constant value:

$C \dfs \left\{ i \in \omega \mid \ex \alpha_0 \all \beta < \sigma_0 \left(
\alpha_0 < \beta \imp C_i (\beta) = C_i (\alpha_0) \right\} \right.$.

Thus $\omega \back C$ is the set of cells that change cofinally in $\sigma_0$.
We continue the run of $F$, calculating further snapshots $s_{\xi'}$ and
$T^2_{\delta (\xi')} [x]$ for $\xi' \geq \xi$. Then there is a strictly increasing sequence of least
$\xi_n $ with $\xi = \xi_0$ so that, for $n > 0$:

(1) $s_{\xi_n} \rest n$ agrees with the cell values for those $i \in C \cap
n$;

(2) Every cell in $\left( \omega \back C \right) \cap n$ has changed value at
least once in $[\xi_{n - 1}, \xi_n)$.

We obtain codes for the $\xi_n$ from $T^2_{\delta (\xi_n)} [x]$ which are set
to one side. By design setting $\xi_{\omega} \dfs \sup_n \xi_n$, then
$s_{\xi_{\omega}} = s_{\sigma_0}$, by appealing to the liminf rule;  thus $\xi_{\omega} = \sigma_0$ (by
definition of $\sigma_0$). We then assemble a code $w_{\sigma_0} \in
\tmop{WO}$ from the cofinal $\omega$ sequence of (set aside codes for) the $\xi_n$.
Then we specify $F (f, e, m, x) \da w_{\sigma_0}$. \\

This completes the description of $F$. By the Recursion Theorem there is $\bar{f}$ so that $F (\bar{f}, e, m, x) =
\{ \bar{f} \}^{\mathsf{E}} (e, m, x)$.
One may
then check by induction on $\tmop{rk} (\mathfrak{T}^{\mathsf{E}} (e, m, x)
\nobracket$ that  $\{ \bar{f}
\}^{\mathsf{E}} (e, m, x)$ is indeed the desired $H (\mathsf{E}, e, m,
x)$. Then by the $S^m_n$-Theorem there is a p.r. $k$ so that
$P^{\mathsf{E}}_{k (e)} (m, x) = P^{\mathsf{E}}_{\bar{f}} (e, m, x)$.\\ \mbox{ } \qed\\

Similar, but simpler, arguments show we can retrieve just the top level looping length.
\begin{corollary}
  \label{lengths} Suppose $\etwoe{e} (\tmmathbf{m}, \tmmathbf{x})$ is any
  recursion effected by the program $P^{\mathsf{E}}_e$, 
Suppose
  $P^{\mathsf{E}}_e$ has as least pair of looping ordinals $\xi_0 <
  \tau_0$. Then, uniformly in $e$ we may find $e'$ so that $\etwoe{e'}
  (\tmmathbf{m}, \tmmathbf{x})$, computes a code for $\tau_0$.
\end{corollary}

\begin{lemma}
  \label{R-lemma}There is a p.r. function $k$ so that if $\{ e \}^{\mathsf{E}}
  (m, x)\da z$ 
  
 then $\{ k (e) \}^{\mathsf{E}} (m, x)\da r$ where $r\in WO \wedge ||r||=\rho^\mathsf{E}(\pa
 {\pa{e,m,x},z})$.
\end{lemma}

\pf This is not too dissimilar to the last argument. We just sketch the main outline and let the reader fill in the 
computational aspects if desired.
The intended procedure $P^\mathsf{E}_{k(e)}$ again runs a simulation of $\{ e \}^{\mathsf{E}}
  (m, x) $, but instead of summing the lengths of the subcomputation calls, inspects the ranks of the subcomputations which we can assume passed up to the calling query as reals coding ordinals in the same manner as the last proof (and going back to Lemma \ref{Ex2.6}). A scratch tape keeps a running value $\rho$ of the supremum of ordinals so far received from below. Each time an ordinal (real) $\eta$ is passed up it is compared to the running tab value $\rho$. If $\eta < \rho$ then nothing is done and the simulation of $\{e\}$ proceeds. Otherwise, if  $\eta = \rho$, then $\rho$ is increased by one;  if $\eta > \rho$ then $\eta$ becomes the new $\rho$. We use the same trick of keeping track at limit stages of the ordinal supremum of these $\rho$-values, by having in a scratch tape the theory $T^2_\rho[x]$; each time $\rho$ is increased by an ordinal amount $\delta$ say, then an additional $\delta$ stages are run by an $L[x]$-theory machine, starting from $T^2_\rho[x]$. The register $R$ holding these theories now has $\widehat T^2_{\rho+\delta}[x]$ and arithmetic in this is the $L[x]$-least code for a wellorder of length $\rho+\delta$. This device ensures that at a limit stage of the process the liminf theory in $R$ codes the supremum of the previous $\rho$ values.
\qed

\begin{corollary} There is a p.r. function $h$ so that for any $(e,m,x)$, if ${\frak{T}}^\mj(e,m,x ) $ is wellfounded,  then $\{h(e)\}^\mathsf{J}(m,x)\da r \in WO$ where 
$||r|| = rk({\frak{T}}^\mj(e,m,x ))$.
\end{corollary}
\pf This is just a specialisation of the previous lemma. \qed\\

\begin{theorem}[Spector-Gandy  I]\label{SpectorGandy} The following are equivalent for an $A\sset
\omega$:\\
(i) $A$ is semi-recursive in $\mi$;\\
(ii) There exists a $\Sigma_1$ $\vp(v_0)\in {\call}_{\dot I}$ so that
 $$m\in A \equi L_{\alpha^\mi_0}[\mi]\models\vp[m,\mi];$$
(iii) There exists a $P$ recursive in $\mi$ so that 
$$m\in A \equi \ex y \mbox{ recursive in } \mi \, (P(m,y)).$$
\end{theorem}
\pf  {\em (i)$\Imp $(ii)} Let $A= \dom(\{e\}^\mi(m,0))$.  Let $h$ be such that  $\{h\}^\mi(m,0)\da r$ where $||r||=H^\mi(e,m,0)$. Then $r$ is $\mi$-recursive, and so $||r||<\alpha_0^\mi$. However then the computation $\{e\}^\mi(m,0)$ can be run inside of 
$ L_{\alpha^\mi_0}[\mi]$. So we have $m\in A$ iff $ L_{\alpha^\mi_0}$ is a model of the $\Sigma_1$-statement: ``{\em There is a code for the run of a convergent computation witnessing }
$\{e\}^\mi(m,0)\da$''.\\
 {\em (ii)$\Imp $(iii)} Note for a fixed $\Sigma_1$ $\vp$, the relation  $P(m,y)\equi y\in WO\wedge L_{||y||}[\mi]\models\vp(m)$ is $\Pi^1_1(\mi)$ and so recursive in $\mi$. \\
 
 $
 \hspace{-2em}\begin{array}{rcl}
  L_{||y||}[\mi] \models\vp(m) & \Equi & \ex \gamma < \alpha_0^\mi \, L_{\gamma}[\mi]\models\vp(m)\\
  & \Equi & \ex y( y \mbox{ recursive in } \mi \wedge y\in WO\wedge  L_{||y||}[\mi]\models\vp(m)); \\
  & \Equi & \ex y   \mbox{ recursive in } \mi \, (P(m,y)).
 \end{array}
 $\\
{\em (iii)$\Imp $(i)} Suppose 

$
\begin{array}{rcl}
 m \in A & \Equi &  \ex y   \mbox{ recursive in } \mi \, (P^\mi(m,y)) \mbox{ then } \\
 & \Equi &\ex e P^\mi
(m,\lambda n.\{e\}^\mi(n,0)).
\end{array}
$

As the class of functionals recursive in $\mi$ is closed under a) substitution by partial recursive functions, and b) existential number quantification, this makes $A$ semi-recursive in $\mi$. \qed \\

There is also the parametrized version   of the above:
\begin{theorem}[Spector-Gandy Theorem II] The following are equivalent for an $A\sset
\omega\times \bai$:\\
(i) $A$ is semi-recursive in $\mi$;\\
(ii) There exists a $\Sigma_1$ $\vp(v_0,v_1)\in {\call}_{\dot I}$ so that
 $$(m,z)\in A \equi L_{\alpha_0^{\mi,z}}[\mi,z]\models\vp[m,z,\mi];$$
(iii) There exists a $P$ recursive in $\mi$ so that 
$$(m,z) \in A \equi \ex y \mbox{ recursive in } \mi,z \, (P(m,z,y)).$$
\end{theorem}

\nod Note the analogy with (lightface) Kleene recursion; for (ii)  we have that:
  $A \sset \omega\times^{\omega} \omega$ is Kleene-$\mathsf{I}$-semirecursive  iff there is
  some $\Sigma_1$ $\varphi$ so that
  
  $$(m,x) \in A \equi L_{\omega^{\mathsf{I}, x}_1} [\mathsf{I}, x] \models \varphi [x,m,
  \mathsf{I}].$$

For $\mi = {\me}$ we can state something simpler.
\begin{corollary} 
There exists a $\Delta^1 _1\, \mathbb{P}$  so that
 for any $A\sset
\omega\times \bai$:\\
 $A$ is semi-recursive in $\me$ 
iff there is $e$ with $A= \dom(\{e\}^\me)$ and
$$\all m\all z [(m,z) \in A \equi \ex y \mbox{ recursive in } \me,z \, (\mathbb{P}(e,m,z,y))].$$
\end{corollary}
\pf Take $\mathbb{P}$ here to say that $y$ codes a complete wellfounded computation of the convergent $\{e\}^\me(m,z)$. \qed

\subsection{Nested Extendible Pairs}

Bound up with the notion of levels of a computation tree is that of the depth
of nesting of ordinals which we proceed to analyse. This will be crucial in
our investigation of the ordinal $\alpha_0^{\mathsf{E}}$.

\begin{definition}
  For $m \geq 1$ an $m$-{\tmem{depth $\Sigma_2$-nesting}}, or just
  $m$-{\tmem{nesting}}, of an ordinal $\alpha$ is a sequence with: $$(\zeta_n,
  \sigma_n)_{n < m} = \zeta_0 \leq \z_1 \leq \cdots \leq  \zeta_{m - 1} <\a < \sigma_{m - 1} < \cdots <
  \sigma_0$$ 
  with: if $k < m$ then $L_{\zeta_k} \prec_{\Sigma_2} L_{\sigma_k}$.
\end{definition}

\subsubsection{The $\Sigma_2$-extendibility tree}

(This subsection is not needed for what follows.) There are a number of things
to discover about nestings. One can define a tree of such pairs, as follows.

\begin{definition}
  (The $\Sigma_2$-extendibility tree) We let $(\mathcal{T}, \prec)$ be the
  natural tree on such pairs under inclusion: as follows: \ if $(\zeta',
  \Sigma')$, $(\bar{\zeta}, \bar{\Sigma})$ are any two countable
  $\Sigma_2$-extendible pairs, then set \ $(\bar{\zeta}, \bar{\Sigma}) \prec
  (\zeta', \Sigma')$ iff $\zeta' \leq \bar{\zeta} < \bar{\Sigma} < \Sigma'$.
\end{definition}

\noindent {\em Note:} If we had allowed the inequality $\bar{\Sigma} \leq \Sigma'$ rather than
a strict inequality in the last definition we could have defined a larger
strict partial order $\prec'$, and a larger tree $(\mathcal{T}', \prec')$;
however this would not have been wellfounded: if $L_{\Sigma} \models
\Sigma_2$-$\mathrm{Sep}$ then it is easy to see that $\left( \mathcal{T}'
\rest \Sigma + 1, \prec' \right)$ is illfounded.

\begin{lemma}
  Let $\delta$ be least such that $L_{\delta} \models \Sigma_2$-Sep. ; let
  $\alpha$ be maximal so that $\left( \mathcal{T}' \rest \alpha, \prec'
  \right)$ is wellfounded, where $$\mathrm{Field} \left( \mathcal{T}' \rest
  \alpha \right) \dfs \{(\zeta, \Sigma) \mbox{ an extendible pair } \mid
  \Sigma < \alpha\}.$$ Then $\delta = \alpha$.
\end{lemma}

{\pf} $(\leq)$: Suppose $\delta > \alpha$. Then \ $\left( \mathcal{T}' \rest
\delta, \prec' \right)$ is illfounded. So there is an infinite sequence of
extendible pairs $(\zeta_n, \Sigma_n)$ with $(\zeta_{n + 1}, \Sigma_{n + 1})
\subset (\zeta_n, \Sigma_n)$. By wellfoundedness of the ordinals there is some
$m_0$ so that for all $n > m_0$ all $\Sigma_n$ are equal to a fixed $\Sigma$,
whilst $\zeta_n < \zeta_{n + 1}$. Let $\zeta^{\ast} = \sup_n \zeta_n$. Then we
have for $n > m_0$ $L_{\zeta_n} \prec_{\Sigma_2} L_{\zeta_{n + 1}}
\prec_{\Sigma_2} L_{\zeta^{\ast}}$. Then $\zeta^{\ast}$ is not
$\Sigma_2$-projectible, and hence $L_{\zeta^{\ast}} \models \Sigma_2$-Sep. But
$\zeta^{\ast} < \delta$. Contradiction.

$(\geq) : L_{\delta} \models \Sigma_2$-Sep. \ Then $S^2_{\delta}$ is unbounded
in $\delta$. Let $\delta_i < \delta_{i + 1}$ be a cofinal sequence, for $i <
\omega$. Then $\la (\delta_i, \delta) \mid i < \omega \ra$ is a
$\prec'$-descending sequence in $\mathcal{T}' \rest \delta + 1$. So $\alpha
\leq \delta$. \qed\\ 

\subsubsection{Infinite depth nestings}

We shall want to consider non-standard admissible models $(M, E)$ of
$\mathrm{KP}$ together with some other properties. We let $\mathrm{WFP} (M)$
be the wellfounded part of the model. By the so-called `Truncation Lemma' it
is well known (\tmtextit{v.}{\hspace{0.17em}}{\cite{Bar}}) that this well
founded part must also be an admissible set. Usually for us the model will
also be a countable one of ``$V = L$'' (or $L [\tmmathbf{x}]$). \ Let $M$ be
such and let $\alpha = \mathrm{On} \cap \mathrm{WFP} (M)$. By the above
$\alpha$ is thus an `admissible ordinal', {\ie}$L_{\alpha}$ will also be a
$\mathrm{KP}$ model. An `$\omega$-depth' nesting cannot exist by the
wellfoundedness of the ordinals. However an ill founded model $M$ when viewed
from the outside may have infinite descending chains of `$M$-ordinals' in its
illfounded part. These considerations motivate the following definition.

\begin{definition}
  \label{Def3.11}An {\tmem{infinite depth $\Sigma_2$-nesting}}{\tmem{ of
  $\alpha$ based on}} $M$ is a sequence $(\zeta_n, s_n)_{n < \omega}$ with :
  
  (i) $\zeta_n \leq \zeta_{n + 1} < \alpha \subset s_{n + 1} \subset s_n$ ;
  (ii) $s_n \in \mathrm{On}^M$; (iii) $(L_{\zeta_n} \prec_{\Sigma_2}
  L_{s_n})^M$.
\end{definition}

Thus the $s_n$ form an infinite descending $E$-chain (where, as above, $E$ is
the membership relation of the illfounded model) through the illfounded part
of the model $M$.

In order for there to exist a non-standard model with an infinite depth
nesting (of the ordinal of its wellfounded part) then the wellfounded part
will already be a relatively long countable initial segment of $L$ (it is easy
to see that if $\zeta = \sup_n \zeta_n$ then already $L_{\zeta} \models
\Sigma_1$-Separation, hence there can be no infinite depth nesting of
$\omega_1^{\tmop{ck}}$ for example). The next lemma shows the existence of
such nested ordinals.

\begin{lemma}
  \label{Ex3.2} Let $\delta$ be least so that $\bai \cap L_{\delta}$ is a
  model of $\Pi^1_3$-$\mathsf{\tmop{CA}_0}$, or more generally $\delta$
  countable with $L_{\delta} \models \Sigma_2$-Separation, and let $(M, E)$ be
  any admissible non-wellfounded end extension of $L_{\delta}$ with
  $L_{\delta}$ as its wellfounded part. Then there is an infinite depth nesting of $\delta$
  based on $M$.
  \end{lemma}
\pf First we just note that such $(M, E)$ exist by the Barwise
  Compactness Theorem. 
   If $\zeta_0 \in S^2_{\delta} \dfs \{\tau < \delta \mid
  L_{\tau} \prec_{\Sigma_2} L_{\delta} \}$, then $\zeta_0$ has arbitrarily
  large $\Sigma_2$-end extensions $L_{\tau}$ for $\tau < \delta$ - namely any
  $L_{\tau}$ with $\tau \in S^2_{\delta}$. So by overspill considerations,
  $L_{\zeta_0}$ must have an ill-founded $\Sigma_2$-end extension
  $(L_{t_0})^M$. The we can repeat this for $\zeta_1 > \zeta_0$ {\etc}
  \qed
  
\begin{definition}
  \label{Def3.13}Let $\beta_0$ be the least ordinal $\beta$ so that
  $L_{\beta}$ has an admissible end-extension $(M, E)$ based on which there
  exists an infinite depth $\Sigma_2$-nesting of $\beta$.
\end{definition}
 It is not too difficult to show that $\beta_0 < \delta$, indeed that $\beta_0 < \gamma_0$ where $\gamma_0$ is such that $L_{\g_0}\prec_{\S_2}L_{\gamma_1}\models KP$. (The claim is false if the requirement that $L_{\g_1}$ be a model of $KP$ were dropped.) Indeed there are many facts that one can demonstrate about such nesting ordinals, but we only prove here what we need.

\subsection{The length of computations}

We want to investigate the course of possible ittm computations. We shall
consider just the functional $\mathsf{E}$, and most of the time just consider
computations of the form $\etwoe{e} (m)$ for simplicity's sake as the methods
parameterize to reals uniformly. We've seen (Lemma \ref{Ex2.6}) that a computation $\etwoe{e} (m)$
can, by making use of suitable queries, `import' into its scratch tape the OT
of any $\{ f \}^{\mathsf{E}} (k, y)$ (when convergent) for any $f, m, y$ for
which it can formulate the question. It can also calculate a code (Lemma \ref{lengths}) for the
order type of the upper end, $\sigma$, of the looping (perhaps $\Sigma_2$-extendible)
pair $(\xi, \sigma)$ which witness the final looping status of $\{ f
\}^{\mathsf{E}} (k, y)$. The latter indicates that the lengths of convergent
$\mathsf{E}$-computations are likely to be all $\mathsf{E}$-recursive
ordinals. This indeed will turn out to be the case. However that does not
yield a characterisation of such ordinals or even a bound on those lengths. \
This we now want to investigate.

Recall that the snapshot at time $\alpha < H (\mathsf{E}, e, \tmmathbf{m},
\tmmathbf{x})$ is  $$s (\alpha) = \pa{I (\alpha), R (\alpha), \pa{C^{\nu
(\alpha)}_i (\alpha) \mid i < \omega}}.$$ Here for $P_e^{\mathsf{E}} (m)$, the
functional $\mathsf{E}$ is simplicity itself, hence the snapshot function $s
\rest \eta$ for limit $\tau > \eta$ is $\Delta_1^{J_{\tau}}$, whilst $s
(\tau)$ itself is $\Sigma_2^{J_{\tau}}$ by use of the liminf rule.
Consequently if $(\xi, \sigma)$ is an extendible pair, then we shall have that
$s (\xi) = s (\sigma)$, and in particular the computation $P_e^{\mathsf{E}}
(m)$ is being carried on at the same node $\nu (\xi) = \nu (\sigma)$ at these
two times. If this node is $\nu (0)$, the topmost node, then this is a pair of
final looping snapshots. If this node is some other $\nu (\alpha)$ at some
level $\Lambda = k > 0$, then this is a pair of final looping snapshots in
some subcomputation $P_{e_{\iota}}^{\mathsf{E}} (m_{\iota}, y_{\iota})$ -
immediately after which at time $\sigma + 1$ some value, and control, is
passed up to the node immediately above $\nu (\alpha)$ in the tree. Similar
considerations are at play in the following.

\begin{lemma}
  Suppose $(\xi, \sigma)$ is an extendible pair, and $P_e^{\mathsf{E}} (m)$ a
  computation. Then for all $\alpha \in (\xi, \sigma)$, $\Lambda (e, m,
  \alpha) \geq \Lambda (e, m, \xi) = \Lambda (e, m, \sigma)$.
\end{lemma}

{\pf}The latter equation follows by $\Sigma_2$-reflection and the $\liminf$
rule, as above, this will mean the snapshots at $\xi$ and $\sigma$ are
identical. However if $\ex \alpha \in (\xi, \sigma)$, $k = \Lambda (e, m,
\alpha) < \Lambda (e, m, \xi)$, then again by $\Sigma_1$-reflection, there are
unboundedly many $\beta < \xi$ with $k = \Lambda (e, m, \beta)$. Again by
Liminf applied to the levels at stage $\xi$, $k \geq \Lambda (e, m, \xi)$ - a
contradiction. {\qed}

\

Similarly:

\begin{lemma}
  \label{3.13} (i) Suppose we have a 2-nesting $\zeta_0 < \zeta_1 < \Sigma_1 <
  \Sigma_0$. Suppose no $\alpha < \zeta_0$ of the overall computation of
  $P_e^{\mathsf{E}} (m)$ is the start of a final loop and $\Lambda (e, m,
  \zeta_0) = k$. Then no $\alpha < \zeta_1$ is the starting point of a final
  loop, and moreover $\Lambda (e, m, \zeta_1) \geq k + 1$.\\
  (ii) More generally if we have a $p$-nesting $\zeta_0 < \cdots < \zeta_{p -
  1} < \Sigma_{p - 1} < \cdots < \Sigma_0$ and we suppose again that no
  $\alpha < \zeta_0$ is the start of a final loop in the computation of
  $P_e^{\mathsf{E}} (m)$, and that $\Lambda (e, m, \zeta_0) = k$. Then
  $\Lambda (e, m, \zeta_{p - 1}) \geq k + p - 1$.
\end{lemma}

{\pf} We consider first just a $p = 2$-nesting. By $\Sigma_2$-reflection and
the $\liminf$ rule, as above, this will mean the snapshots at $\zeta_0$ and
$\Sigma_0$ are identical; hence $P_e^{\mathsf{E}} (m)$ is still running at
depth $k$ at $\Sigma_0$ at the same node $\nu (\zeta_0) = \nu (\Sigma_0)$. \
Suppose $k = 0$. Thus $P_e^{\mathsf{E}} (m)$ has as its first repeating loop
$[\zeta_0, \Sigma_0]$. Suppose for a contradiction that control is at level 0
also at time $\zeta_1$ (and again also at $\Sigma_1$). So again
$P_e^{\mathsf{E}} (m)$ has looping snapshots at $(\zeta_1, \Sigma_1)$. However
this is a $\Sigma_1$-fact about $P_e^{\mathsf{E}} (m)$ that $L_{\Sigma_0}$
sees: ``{\em There exists a 2-extendible pair $(\bar{\zeta},
\bar{\Sigma}$) with $P_e^{\mathsf{E}} (m)$ having identical
snaphots at level 0 at } $(\bar{\zeta}, \bar{\Sigma})$.'' \ But then there is
such a pair $\bar{\zeta} < \bar{\Sigma} < \zeta_0$ and $P_e^{\mathsf{E}}
(m)$'s computation is again looping at $\bar{\zeta}$ contrary to the
supposition on $\zeta_0$.

The argument for $k \geq 1$ is very similar: if $\liminf_{\alpha \rightarrow
\zeta_0} \Lambda (e, m, \alpha) = \Lambda (e, m, \zeta_0) = k$, then
$\liminf_{\alpha \rightarrow \Sigma_0} \Lambda (e, m, \alpha) = k$ also.
Again, if it entered the interval $(\zeta_1, \Sigma_1)$ at this same level $k$
it would loop there with identical snapshots at $\zeta_1, \Sigma_1$, and by
the same reflection argument applied repeatedly would do so not just once but
unboundedly below $\zeta_0$ at the same level $k$. But after each successful
loop at level $k$, control passes up to level $k - 1$. However then
$\liminf_{\alpha \rightarrow \zeta_0} \Lambda (e, m, \alpha) \leq k - 1$.
Contradiction! \hfill \qed

\begin{corollary}
  \label{Nesting-levels}Suppose $\alpha$ is $n$-nested (for some $n \geq 1$).
  Then $P_e^{\mathsf{E}} (m)$ has not entered a final repeating loop 
  before $\alpha$ only if $\Lambda (e, m, \alpha) \geq n$.
\end{corollary}

{\pf}If $n = 1, \zeta_0 < \alpha < \Sigma_0$ with $(\zeta_0, \Sigma_0)$ an
extendible pair, and were $\Lambda (e, m, \alpha) = 0$, then by $\Sigma_2$
reflection we should have $\Lambda (e, m, \beta) = 0$ for unboundedly many
$\beta < \zeta_0$. But then $\Lambda (e, m, \zeta_0) = 0$ by the Liminf rule,
and $(\zeta_0, \Sigma_0)$ is a final looping pair in the computation, and
$P_e^{\mathsf{E}} (m)$ has entered a  final repeating loop at or before
$\alpha$. For $n > 1$ use induction. \ (Exercise) {\qed}

\begin{lemma}
  \label{Boundedness}(Boundedness Lemma for computations recursive in
  $\mathsf{E}$) Let $\beta_0$ be the least infinitely nested ordinal in some
  ill-founded model $M$ with $\mathrm{WFP} (M) = L_{\beta_0}$. Let
  $\bar{\alpha}_0$ be least with $L_{{\bar\alpha}_0} \prec_{\Sigma_1}
  L_{\beta_0}$. 
 If $\mathfrak{T}^{\mathsf{E}} (e, m)$ is wellfounded then it has tree rank less than $ \bar\alpha_0$.  
\end{lemma}

{\nod}{\pf} Let $\zeta_0 < \cdots < \zeta_n < \cdots \hspace{0.17em} \beta_0 
\hspace{0.17em} \cdots \subset s_n \subset \cdots \subset s_0$ witness the
infinite nesting at $\beta_0$ in $M$. By the definition of $\bar{\alpha}_0$
there is no least $\alpha \in [\bar{\alpha}_0, \beta_0)$ so that $L_{\alpha}$
sees that $P_e^{\mathsf{E}} (m)$ has a repeating looping snapshot as this
would be a $\Sigma_1$-fact true in $L_{\beta_0}$; but then by
$\Sigma_1$-reflection, it is true in $L_{\bar{\alpha}_0}$ and
$P_e^{\mathsf{E}} (m)$ would then be convergent before $\bar{\alpha}_0$.
However if $P^{\mathsf{E}}_e (m)$ has not failed before $\beta_0$, it will do
so by $\beta_0$: the previous lemma shows that $\Lambda (e, m, \zeta_n) <
\Lambda (e, m, \zeta_{n + 1})$ holds in $M$. \ But these level facts are
absolute to $V$, as they are grounded just on the part of the computation tree
being built in $L_{\beta_0}$ as time goes towards $\beta_0$; so
$P_e^{\mathsf{E}} (m)$'s computation tree $\gott (e, m)$ will have an
illfounded branch at time $\beta_0$.\hfill {\qed} \\

Or, rephrased:
\begin{corollary}\label{Boundednesslengths}
If $\{e\}^\me(m)\da$ then the overall length of computation \linebreak $H(\me,e,m,\emp) < \bar\a_0$. Thus $\bar\a_0$ bounds the convergence times of all computations of the form $\{e\}^\me(m)$.
\end{corollary}

The above then shows that the initial segment $L_{\bar{\alpha}_0}$ of the
$L$-hierarchy contains all the information concerning looping or convergence
of computations of the form $P^{{\mathsf{E}}}_e (m)$. 
This leads to:

\begin{corollary}
$\kappa^\me = \a_0^\me \leq {\bar\alpha}_0.$
\end{corollary}

In the next section we shall see that $\a_0^\me = {\bar\alpha}_0$.\\

A computation may
then continue through the wellfounded part of the computation tree for all the
times $\beta < \beta_0$ but if so, it will be divergent. Relativisations to
real inputs $\tmmathbf{x}$ are then straightforward by defining $\beta_0
(\tmmathbf{x})$ as the least such that there is an infinite nesting based at
that ordinal in the $L_{\nosymbol} [\tmmathbf{x}]$ hierarchy {\etc}

Having ascertained ${\alpha}_0$ as an upper bound for convergent
computations of the form $\ptwoe{e} (m) \nocomma$, we now look to show that
this is best possible. We have seen that we can compute codes for ordinals for
increasing levels in the $E^{\alpha}$-hierarchy, indeed computation of the
next element of $E^k$ can be effected by a computation using a tree
$\mathfrak{T}^{\mathsf{E}}$ of rank $k$. We also have the natural finite depth
of nesting of convergent computations. We need to have computations that
approach $\beta_0$ in length. A convergent computation that somehow seemingly
required infinite depth nesting would seem to be impossible. This suggests the
following formulation.

We attempt to define a function $t$ with domain $\omega$ by a
$\Pi_1$-recursion, through finite approximation functions defined on initial
segments of $\omega$. Thus to some function $t$, but {\via}approximations $t
\rest k + 1$ where, for $i < k$, $\{ t_{i + 1} \}$ is $\Pi_1$ definable in
$\pa{t_0, \ldots, t_i}$. We shall implement the process as defining a
candidate for $t_k $ as being $\Pi_1$-definable over some $L_{\xi}$ where $\xi$
is admissible. The process will calculate larger and larger admissibles $\xi'$
and check that $t_k$, and indeed all the previous $t_i$ for $i < k$ still
satisfy the same $\Pi_1$ clauses over $L_{\xi'}$ as for the earlier $L_{\xi}$.
If so then an attempt to define a non-empty $\{ t_{k + 1} \}$ is made over
$L_{\xi'}$. On the other hand if for some least $i_0$, \ $t_{i_0}$ fails to
satisfy the $\Pi_1$ clause of its definition over $L_{\xi'}$, then all the
current approximating values $t_j$ for $i_0 \leq j \leq k$ are abandoned, and
the process continues with a search for a new approximating value $\{ t_{i_0}'
\}$ ({\tmem{if it exists}}), and advances to the next admissible level
$L_{\xi''}$. We emphasise existence here, since {\tmem{prima facie}} there is
no guarantee at any point of there being a non-trivial candidate in $L_{\xi'}$
for $t_{i_0}'$, and one may simply have to move to a higher $L_{\xi''}$ to
look again. All in all, we are looking for an $L_\xi$ over which we can find a $\Delta_2$
definable function $t$ whose range is an increasing $\omega$ sequence cofinal in $\xi$.

As just stated then, there is, in general, no guarantee that any instances of
a function $t$ being defined in this quasi-recursive fashion will eventually
be definable at some level $L_{\xi}$. \ Nor should we expect that
approximations will grow in a monotone fashion: $\Pi_1$ properties fail to be
upwards persistent.

Nevertheless we shall be able to find a $\Pi_1$-definable process which will
have the following essential property: suppose $\zeta_0 < \zeta_1 < \alpha <
\sigma_1 < \sigma_0$ is a nesting of $\alpha$ as above and $t \rest 1 = t_0
\in L_{\zeta_0}$ results from applying our $\Pi_1$ definition over
$L_{\zeta_0}$. Then there is to be guaranteed a candidate $t_1 \in
L_{\zeta_1}$ to extend $t \rest 1$ to $t \rest 2$. In fact we shall have
slightly more than that. Let $\varphi (v_0, v_1)$ be our intended
$\Pi_1$-formula for defining approximations on $\omega$ to such a $t$, and the
$L$-least $x$ such that $\varphi \left( x, t \rest k \right)^{L_{\xi}} $ is to
be set as $t (k)$ to extend the approximation to $t \rest k + 1$. Let
$\tilde{\varphi} (u)^{L_{\xi}}$ hold if $u$ is a finite sequence and for all
$k < \tmop{dom} (u)$ $u_k $ is $L$-least s.t. $\varphi \left( u_k, u \rest k
\right)^{L_{\xi}}$. This makes ``$y= t\rest k+1$'' a $\Delta_2^{L_\xi}$ definable predicate.

\

$\Upsilon(\vp)$: {\em 
 suppose $\zeta_0 < \zeta_1 < \cdots < \zeta_k < \alpha
< \sigma_k < \cdots < \sigma_1 < \sigma_0$ is any $k + 1$-nesting of $\alpha$; then
 $\exists t \upharpoonright k + 1 = \la t_0, \ldots, t_ k \ra$ so that. }
$$ \forall i\leq
k[t_i \in
L_{\zeta_i} \, \wedge \,{\varphi} (t_i, t \restriction i 
)^{L_{\zeta_i}}].$$

{\nod}Note that $\tilde{\varphi} \left( t \rest i \right)^{L_{\zeta_{i}}}$ and  $\tilde{\varphi} \left( t \rest i \right)^{L_{\zeta_{i}}}
\equi \tilde{\varphi} \left( t \rest i \right)^{L_{\sigma_{i }}}$ by $\Sigma_2$-reflection,   in any
case.  We can thus say that any $ t \rest i $ satisfying $\tilde{\varphi} \left( t \rest i \right)^{L_{\zeta_{i}}}$ is in fact `good up to $L_{\sigma_{i}}$.'

For any $p + 1$-nesting
we can define a candidate function $t \rest p + 1$ over $L_{\zeta_p}$. All
that would be needed to define a complete $\Delta_2$ function $t$ with domain $\omega$ at $L_{\alpha}$ is an
$\omega$-nesting! We have yet to justify that there is any $\varphi$ so that
$\Upsilon (\varphi)$ holds.

\subsection{Instantiating $\Upsilon$}

The first proof (in the early 2010's) concerning $\beta_0$ as the least upper bound for
wellfoundedness of computation trees of ittm-recursions in $\mathsf{E}$ (on
integer input) went {\via}a proof $\Sigma^0_3$-Determinacy. The principle
$\Upsilon$ was validated for a ${\varphi}$ that looked for successive
non-losing quasi-strategies for games $G (B_n)$ where a $\Sigma_3^0$ set $A =
\bigcup_n B_n$ for $B_n \in \Pi^0_2$, which Player \tmtextit{II} is to win. It
was known that an overall strategy for Player \tmtextit{II} in $G (A)$ was to
be found definably over $L_{\beta_0}$ ({\cf}{\cite{W2011}}). \ Hence having a
generalized ittm recursion that sought for such quasi-strategies brought one
ever closer to this level of the $L$-hierarchy.

\

Hachtman in {\cite{Ha18}} then showed that the reals of any admissible model
$L_{\beta} \models V = \tmop{HC}$ on which could be based an $\omega$-nesting,
formed a model of $\Pi^1_2$-monotone induction. He also showed ({\tmem{inter
alia}}) that $L_{\beta_0}$ yielded the least $\beta$-model of \ $\Pi^1_2$- MI
. That is, he showed that the reals of the least level of $L$ that satisfied \
$\Pi^1_2$- MI, were the reals of $L_{\beta_0}$. However this also went
indirectly {\via}the result on $\Sigma^0_3$-Determinacy first holding
definably over $L_{\beta_0}$.

\

It was thus desirable to have an ittm recursive in $\mathsf{E}$ procedure
which was more natural in the context of such recursions in order to, {\eg},
characterise the halting problem $H$ for such recursions/machines, rather than
searching for quasi-strategies in a game theoretic fashion.

What we effect below is a search for admissible levels of $L$ for finitely
many ordinals, where it is consistent for those ordinals to be the left hand
end of a nested sequence in an $\omega$-model which is an end extension of
that level. This sequence $(t_0, \ldots, t_k)$ will be of approximations to an
eventual $\omega$-sequence witnessing an infinite depth nesting in some
$\omega$-model which as we have seen yields an upper bound for the ranks of
wellfounded computation trees of ittm-recursions in $\mathsf{E}$ (on integer
input) - see the Boundedness lemma above Lemma \ref{Boundedness}.

We use the notion of an infinitary logic: for $\alpha$ an admissible, we set
$\call_{\alpha} \dfs L_{\omega_1, \omega} \cap L_{\alpha}$ and is thus an
`admissible fragment' in the sense of {\cite{Bar}}, {\cite{Ke71}}. More
specifically we let $\Gamma_{\alpha}$ be the axioms of $\call_{\alpha}$
comprising of: (i) $\tmop{KP}$; (ii) the $\call_{\alpha}$-infinitary atomic
diagram of $L_{\alpha}$; (iii) for any $y \in L_{\alpha}$: $\all x ( x
\in c_y \imp \bigvee_{u \in y} x = c_u )$ (we assume there are constants
$c_u$ in $\mathcal{\call_{\alpha}}$ for each $u \in L_{\alpha}$). Any model of
$\Gamma_{\alpha}$ has then (a transitive copy of) $(L_{\alpha}, \in)$
contained in its wellfounded part (by (ii)) and is indeed an end-extension of
$(L_{\alpha}, \in)$ by (iii). We set, with $\tmop{WFD} (A)$ denoting the
wellfounded part of a structure $A$:

\

{\tmem{Let}} $k \geq 0 : \Phi (\alpha, (t_0, \ldots, t_k)) \equiv$

\noindent {\tmem{``It is consistent in $\call_{\alpha}$-logic that there exists an
end-extension}} $(N, E) \supset (L_{\alpha},\in)$

$\left[ \tmop{WFD} (N) \supseteq L_{\alpha} \wedge \ex s_k, \ldots, s_0 \in
\tmop{On}^N (( \right.${\tmem{$ t_0, \ldots,
t_k, s_k, \ldots, s_0) $ forms a $k + 1$-nesting with}} $\nobracket t_k <
\alpha)$].''

{\nod}

Then $\Phi (\alpha, (t_0, \ldots, t_k))$ is a $\Pi_1^{L_{\alpha}}$-expressible
assertion about the sequence $(t_0, \ldots, t_k)$ since a {\tmem{consistency
property,}} containing the theory $\Gamma_{\alpha}$ together with the
assertion about the existence of a $k + 1$ nesting \etc\!, if it exists, is so
definable over $L_{\alpha}$, and thus a model $(N, E)$ of the
$\Gamma_{\alpha}$ axioms with the nesting property satisfying the property
$\Phi$ exists in $L_{\alpha^+}$, the next admissible set. Recall that we are low down in the $L$-hierarchy: every level $L_\gamma$ we are concerned with is a model of ``$V=HC$''. Thus models of consistent theories are available in the next admissible set.

\begin{lemma}
  \label{L4.7} If $(\xi, \sigma)$ is an extendible pair, there is $t_0 < \xi$
  so that
  $$L_{\sigma} \models\mbox{``}\ex
  t_0 \all \alpha \in (t_0, \sigma) \left( \alpha \in \tmop{ADM} \imp \Phi
  (\alpha, (t_0)) \right)\mbox{''}.\hspace{3em} (\ast) $$
\end{lemma}

{\pf} Note first that $ \all \alpha \in (\xi, \sigma) \left( \alpha \in
\tmop{ADM} \imp \Phi (\alpha, (\xi)) \right)$. This is because the
$\call_{\alpha}$-theory of $L_{\alpha}$ is consistent with there being such an
extension in which $\alpha$ is 1-nested, since $(N, E) = (L_{\sigma}, \in)$ witnesses
such. Thus $(\ast)$ holds
 (namely take $t_0
= \xi$). This is a $\Sigma_2$ statement and so by $\Sigma_2$ reflection goes
down to $L_{\xi}$; there is such a $t_0 < \xi$, which `survives' playing this
role, meaning that it satisfies $\Phi (\alpha, (t_0))$ for all admissibles
$\alpha < \sigma$ all the way through to $\sigma$. \ {\qed}

\

The moral of the last proof, is that as admissibles $\alpha \imp \xi$ then
such a $t_0$ can be found in $L_{\alpha}$, and eventually it will settle down
on some value $t$ which will shall say ``is good up till $\sigma$''.

\

\begin{remark}
  \label{R4.7} \ Note however in the above if $(\xi, \sigma)$ was the minimal
  extendible pair $(\zeta, \Sigma)$, that if $\alpha \in \tmop{ADM}$ is the least
  greater than $\Sigma$ then a $t < \zeta$ cannot satisfy $\Phi (\alpha,
  (t))$: if there were some model as specified with $(N, E) \supset
  (L_{\alpha}, \in)$ we should have a 2-nesting in $N$, with $t < \zeta <
  \Sigma < s$. By $\Sigma_1$-reflection in $N$ there would have to be an
  extendible pair $(\zeta', \Sigma')$ in $L_t$ which is impossible. Thus such
  a $t$ does not survive (playing its role) beyond $\Sigma$, {\ie}is not good
  beyond $\Sigma .$ On the other hand for other extendible pairs $(\bar{\xi},
  \bar{\sigma})$, for example those that are themselves contained within
  another extendible pair $(\xi, \sigma \nobracket$), then such a $t$ can and
  does survive well past $\bar{\sigma}$; indeed it may do so till $\sigma$.
\end{remark}

\begin{remark}
  From $\Phi$ being in $\Pi_1^{L_{\alpha}}$, one can write down a $\Pi_1$
  $\varphi$ so that $\Upsilon (\varphi)$. Namely, let 
  $\varphi \left( t_n, t \rest n \right)$ be $$\all \alpha > t_n \left( \alpha
  \in \tmop{ADM} \imp \Phi (\alpha, (t_0, \ldots, t_n)) \right) .$$
\end{remark}
\begin{lemma} \label{instupsilon}
  $\Upsilon (\varphi)$ holds. 
  \end{lemma}

{\pf} We verify $\Upsilon (\varphi)$ with $n$ replacing the $k$ there.
 As an example with $n=2$ let $\zeta_0 < \zeta_1 < \z_2 <\Sigma_2 < 
\Sigma_1 < \Sigma_0$ be any three nesting  of admissibles $\beta$ in the interval $(\zeta_2 ,
\Sigma_2)$. 
Then $\zeta_i, \Sigma_i$ are of the form \ $\zx{n - i}, \sx{n - i}$\quad and
$(\zeta_0, \Sigma_0)$ is an $n = 2$-extendible pair.  

The following $\S_2$ statement holds in $L_{\S_0}$:\\

\nod$L_{\S_0}\models$\\

\nod $ \mbox{ `` } \ex t_0 \left[ {\varphi}(t_0 , \emp)\wedge \mbox{\em For any $1$-extendible pairs  } 
 (\bar \zeta_1,\bar \z_2,\bar\S_2,\bar \Sigma_1) \mbox{\em with } t_0 < \bar\zeta_1:\right.$\\
 
 $\quad L_{\bar\S_1}\models \mbox{ `` } \ex t_1 \left[ \mbox{\em For any $0$-extendible pair  } 
 (\bar \z_2,\bar\S_2) \mbox{\em with } t_1 < \bar\zeta_2:\right.$\\

\indent\quad \quad$\left.\left.  L_{\bar\Sigma_2} \models \mbox{ `` } \ex t_2 \all \alpha >t_2( \alpha \in \tmop{ADM}\imp 
\Phi(\alpha, (t_0,t_1,t_2))\hspace{0.1em})\mbox{"}\hspace{0.1em} \right] \hspace{0.1em}\mbox{"}\hspace{0.1em}  \right] \hspace{0.1em}\mbox{"}. $\\

That this holds 
follows by taking $\zeta_0,\bar\zeta_1,\bar\zeta_2$ as witnesses for the existentially quantified $t_0, t_1,t_2$ respectively in the quoted formulae. Then clearly  ${\varphi}(\z_0 , \emp)^{L_{\S_0}}$  holds, for any admissible $\alpha\in (\z_0,\S_0)$. By $\S_2$-reflection the statement holds of some least $t_0<\zeta_0$ which we fix. Then, as $\vp$ is $\Pi_1$, we shall also have here ${\varphi} \left(t_0 , \emp \right)^{L_{\zeta_0}}$.  By taking our choice of $(\z_1,\S_1)$ in place of  $(\bar\z_1,\bar\S_1)$ we have:\\

 \nod$L_{\S_1}\models \mbox{ `` } \ex t_1 \left[\varphi(t_1 ,t_0 )\wedge \mbox{\em For any $0$-extendible pair  } 
 (\bar \z_2,\bar\S_2) \mbox{\em  with } t_1 < \bar\zeta_2:\right.$\\
 
\indent\quad $\left. L_{\bar\Sigma_2} \models \mbox{ `` } \ex t_2 \all \alpha >t_2( \alpha \in \tmop{ADM}\imp 
\Phi(\alpha, (t_0,t_1,t_2))\hspace{0.1em})\mbox{"}\hspace{0.1em} \right] \hspace{0.1em}\mbox{"}\hspace{0.1em}  . $\\

\nod That this holds 
in $L_{\Sigma_1}$ 
follows by now taking $\zeta_1,\bar\zeta_2$ as witnesses for the existentially quantified $ t_1,t_2$ respectively.\\

\nod Again by $\S_2$-reflection the least such $t_1$ is $<\z_1$ which we may fix. Thus 
${\varphi} \left(t_1 ,t \rest 1 \right)^{L_{\zeta_1}}$ is also verified. Now taking $(\z_2,\S_2)$ in place of  $(\bar\z_2,\bar\S_2)$ we have:

\indent\quad $  L_{\Sigma_2} \models \mbox{ `` } \ex t_2 \all \alpha >t_2( \alpha \in \tmop{ADM}\imp 
\Phi(\alpha, (t_0,t_1,t_2))\hspace{0.1em})\mbox{"}.$\\

\nod \nod Again by $\S_2$-reflection the least such $t_2$ is $<\z_2$. As our choice of 2-nesting ordinals was arbitrary this completes the verification of $\Upsilon(\vp)$ in this case. \\

 By expanding the
argument  for each nested
$$\zeta_0 < \zeta_1 < \cdots < \zeta_{n - 1} < \alpha < \Sigma_{n - 1} < \cdots
< \Sigma_1 < \Sigma_0$$
\noindent   working from the outermost interval inwards, just as above, there is some $(t_0, \ldots \nocomma, t_{n - 1}) $,  with
${\varphi} \left(t_j ,t \rest j \right)^{L_{\zeta_j}}$ holding for each $j< n$,
so that for each extendible  interval
$\left( \zx{}_n, \sx{}_n \right)$ contained in $\left( \zx{}_{n - 1}, \sx{}_{n -
1} \right)$ we can associate a $t_n < \zx{}_n$ for which
$$\left[ \all \alpha > t_n \left( \alpha \in \tmop{ADM} \imp \Phi
(\alpha, (t_0,  \ldots, t_n)) \right) \,  \right]^{L_{\Sigma_n}}.$${\qed}

\section{A Least Upper Bound to $\mathsf{^2 E}$ Computation Lengths}

We thus take the $\Pi_1$ formula $\varphi$ from the last section, with the
property $\Upsilon (\varphi)$ \ above, and with $\tilde{\varphi} (u) \equiv
\all k < \dom (u) (u_k \nobracket$ {\tmem{is $L$-least so that}} $\left.
\varphi \left( u_k, u \rest k \right) \right)$ the associated
$\Delta_2$-formula. We shall then see how to define a recursion in
$\mathsf{E}$ that fails, that is, has an illfounded computation tree, but it
only becomes illfounded at stage $\beta_0$. It thus outlasts all the
convergent computations $\ptwoe{e} (m)$, which by the Boundedness Lemma
\ref{Boundedness} all converge before ${\alpha}_0 < \beta_0$. One might be
curious about the interval $({\alpha}_0, \beta_0)$. At risk of whimsy
(something always to be avoided) - ${\alpha}_0$ appears as some kind of
event horizon: any computation \ $\ptwoe{e} (m)$ that enters this region is
destined to fail with an illfounded computation tree; the computation
continues, remaining wellfounded for all $\alpha < \beta_0,$ but eventually
falls into the black hole of $\beta_0$.

\

We give a scheme for a process that will be implementable as a recursion in
$\mathsf{E}$. In order for its description to be manageable we shall use some
abbreviations. We have seen above that recursions in a suitable $\mathsf{K}$
which for queries $Q^{\mathsf{K}} (e, m, y)$ formally return $\mathsf{K} (z)$
when $\{ e \}^{\mathsf{K}} (m, y) \da z$, can usefully be organised as
sequences of such queries which have (a) the effect of returning simply a
check as to whether $\{ e \}^{\mathsf{K}} (m, y) \da {\back \ua}$, (provided $\mathfrak{T}^\me(m,y)$ is wellfounded),\ or
indeed (b) returning the infinite sequence of values $z (0), \cdots \nocomma$ if $\{ e
\}^{\mathsf{K}} (m, y) \da {z}$.

We thus can ask for the whole result $z \in \bai$ where \ $\{ e
\}^{\mathsf{K}} (m, y) \da {z}$ during a computation call. We can also
ask for the {\tmem{length}} $\sigma$ of the looping computation $\{ e
\}^{\mathsf{K}} (m, y) \da z$ or its overall length \ $H (\mathsf{E}, e, m, y)
$(see Lemma \ref{H-lemma}), and get either in the form of a code $W_{\sigma}
\in \tmop{WO}$; we can also ask for $L_{\sigma} [z]$ {\etc\nolinebreak\!,} or for the next
admissible ordinal $\sigma^+$ {\etc} We shall just describe ``general
queries'' as those of this form that can be realised by (perhaps infinite)
sequences of official queries. We regard such general queries as passing
infinite amounts of information to lower levels of the computation tree as
data for that subcomputation; and for the results, again as infinite sequences
or amounts of information to be passed up to the controlling program at the
next level up. This is all just for the sake of brevity and comprehensibility.
The queries we shall freely formulate in English to describe a process to be
transacted at the next level down without involving us in all of the formal
{\tmem{minutiae}}. We shall do this without further comment in the confidence
that the enthusiastic reader (if there are any) could if they wished, with
effort, convert these general queries into the officialese of our formalism.

\

We now outline the algorithm at the various levels of computation in the
oracle calls of a {\em (master) computation} which runs at node $\nu_0$
(the only node at level $\Lambda = 0$). We proceed by describing the actions
of the programmes being called, which the reader, as we have just said, may
reformulate as official queries to the $\mathsf{E}$-functional as oracle.

\

At the end of the description we justify the claim that this is a bona fide
${\mathsf{E}}$-recursion of the form $\{e_0 \}^{{\mathsf{E}}}
(\nobracket k, t, W \nobracket)$. In the sequel $W$ is intended to be a subset
of $\omega$ coding a well ordering; if the well ordering is of length $\delta$
then it is intended to be the $L$-least well ordering of this length, and we
indicate this by writing $W_{\delta}$.  (In the region of $L$ under discussion
(i) $L_{\delta} \models${\em ``Every set is countable''}; (ii) there is always a
$\Sigma_2^{L_{\delta}}$ map (possibly requiring a parameter) of $\omega$ onto
$L_{\delta}$ (or a $\Sigma_3^{L_{\delta}}$ such map parameter free); (iii)
uniformly definable $\Sigma_n$-skolem functions for the $J_{\delta}$ (or
$L_{\delta}$ if $\tmop{Lim} (\delta)$); \ $W_{\delta}$ is thus always an
element of $L_{\delta + 1}$.)  During the proof \ $^0 M$ is a variable for a transitive
structure $L_{\tau} \models \tmop{KP}$.
\

The programme $\{e_0 \}^{{\mathsf{E}}} (\la 0, t_0, W_{\sigma} \ra)$
operates as follows, starting with $t_0 = \emp$ and $W_{\sigma} = W_{\omega}$.
The reader might like to bear in mind that the intention below is that a
positive solution to query $Q^k $ is a successful discovery of a next element
$t_k$ to add to a chain $(t_0,\, \ldots\, ,t_{k-1})$, and which takes place at level $\Lambda = k$. At a level $\Lambda = k$
there will be a register $R_k$ for purposes of bookkeeping theories
$T^2_{\alpha}$. The general maxim is that if during the program a node at
some level $k$ is passed an ordinal in the form $W_{\sigma}$, then $R_k$ is
immediately updated by writing out successively the theories $T^2_{\alpha}$
for $\delta < \alpha \leq \delta + \sigma$, where $T^2_{\delta}$ is the theory
currently held in $R_k$. The point of the bookkeeping is that a code for the ordinal $\d+\sigma$ is then, at worst, recursively enumerable in the liminf theory $\widehat T_{\d + \s}$ in $R_k$.

At $\nu_0$ (so at $\Lambda = 0$):

The master computation $\ptwoe{e_0} (0)$ is run at this level. \ As described
in the maxim above, it first writes out the sequence of theories
$T^2_{\alpha}$ for successive $\alpha$ to a register $\calr_0$ for
book-keeping purposes, for $\alpha$ any ordinal passed to it: here $\sigma$ as
given by $W_{\sigma}$. That done, it computes successive $\Sigma_1$-ADM
%[XXXXCheck $1$-extendible] 
structures of the form $^0 M$ of increasing
length starting from the least admissible $\geq \sigma$, looking for a
candidate $t_0$ for which $\varphi (t_0)$ holds, and which it may write to its
OT. Whilst doing this, it bookkeeps by continuing to add to the theories in
$R_0$, for each ordinal step in computation it takes. When it finds a $t_0$ in
some structure $^0 M$ so that ($\nobracket \varphi (t_0) )^{^0 M}$,
 it launches
a query to a subroutine one level down which essentially asks ?Q: {\tmem{Is
there a stable candidate for $t_1$ so that}} $\tilde{\varphi} \left( \la
\left. t_0, t_1 \ra \right. \right)$? \ Slightly more formally:

\

$Q^1 :  $
{\em Does $\ptwoe{e_0} \left( \pa{1, \nocomma \pa{t_0}, W_{\sigma}}
 \right)$ find, starting from the ordinal $\sigma = \tmop{On} \cap\,^0
M$, a candidate for $t_1$ which is stable in  admissible structures $^0M=L_\b$ of increasing length $\b>\sigma_0$, using $t_0$ so that $\tilde{\varphi} (
\la  t_0, t_1 \ra)
  $?
  
  }

\

Thus the information $x = \pa{1, \nocomma \pa{t_0}, \sigma}$ passed down in
the query contains the `current attempt' at $t_0$ and the `current ordinal',
that is the height of the `current admissible structure' $^0 M$. A `stable
candidate' is one which eventually settles down to a fixed value in
successively longer admissible structures $\,^0 M$, and could be written to
the local OT.

\

At $\Lambda = 1$:

$P_{e_0} (\nobracket \left. \pa{1, \nocomma \pa{t_0}, W_{\sigma}} \right)$
computes further successive admissible structures of the form $\,^0 M$
starting from $\sigma$, again writing theories $T^2_{\alpha} $ to a register
$\calr_1$ for $\alpha \geq \sigma$. Various alternatives may happen:

{\bu} (1) {\tmem{En route}} $\varphi (t_0)$ may become invalid in some $\,^0 M
$structure (as it is a $\Pi_1$ property, and thus $t_0$ failed to survive and
`stably' satisfy $\varphi$ in some longer structure). In which case, if this
happens in some least such structure, and a new candidate $t_0'$ is
present, it and the current ordinal $\sigma'$ (the ordinal height of the
current $\,^0 M$) is passed back up to $\nu_0$ at $\Lambda_0$ as
$W_{\sigma'}$. (If there is no $t_0'$ yet then the `empty candidate' $\emp$ is
passed up instead.)

{\bu} (2) However if $\varphi (t_0) $ continues to hold in successive
admissible structures, it may be that no potential $t_1$
is found. This means the process at this level continues until it loops. The
answer to $Q^1$ is then ``No''. Setting the new $\sigma'$ to be the length of
this loop, then ``No'' and $t_0$, and $W_{\sigma'}$ are passed back up to
$\Lambda = 0$, {\ie}the node $\nu_0$. (Strictly speaking, the looping length
ordinal $\sigma'$ is calculated{\tmem{ at}} the higher level node, see here the method of Lemma \ref{lengths}, thus here
at $\nu_0$, rather than being `passed up' to the node. But we say this in
interest of brevity: we shall use this circumlocution without further comment
below.) The theory, $T^2_{\delta}$ say, residing in $R_0$ at $\nu_0$, is then
extended by single steps from its previous length to that passed up: hence to
$T^2_{\delta + \sigma'} .$

{\bu} (3) Lastly if a $t_1$ is found so that $\tilde{\varphi} \left(  \la
t_0, t_1 \ra \right)$ holds in some admissible structure $^0 M$ of
height $\sigma'$, then a subcomputation query $Q^2$ is passed down to a new
node at level $\Lambda = 2$:

\

$Q^2 :$ \ {\em Does $\ptwoe{e_0} \left( \pa{2, \nocomma \pa{t_0, t_1},
W_{\sigma'}}  \right)$ find a stable candidate for $t_2$ starting from
$\pa{t_0, t_1}$
and the ordinal $\sigma' $  so that $\tilde{\varphi} \left( \la \left. t_0, t_1, t_2 \ra
\right. \right) $ ?}

\

At $\Lambda = 2$:

$P_{e_0} (\nobracket \left. \pa{2, \nocomma \pa{t_0, t_1}, W_{\sigma'}}
\right)$ updates $R_2$ using $W_{\sigma'}$, and computes further successive \
$\,^0 M$ admissible structures starting from $\sigma'$. Various alternatives again may
happen:

{\bu}(1) {\tmem{En route}} $\varphi (t_{i_0})$ may become invalid in some
$\,^0 M$ structure for a least $i_0 < 2$. In which case the computation at
this level HALTS, \ and if a new candidate $t_{i_0}'$ is present, it and the
current ordinal $\sigma'$ in the form of $W_{\sigma'}$ is passed back up to
$\Lambda = 1$, all other $t_j$ for $j > i_0$ (if any) being discarded. (If
there is no $t_{i_0}'$ yet then the `empty candidate' $\emp$ is passed up
instead.) If $i_0 = 0$, then it is arranged that this information then is
passed yet further up to level $\Lambda = 0$. Otherwise control is passed up
to $\Lambda = 1$ where it now remains.

{\bu} (2) However if $\tilde{\varphi} \left( \left. \la t_0, t_1 \right. \ra
\right)$ continues to hold in successive $\,^0 M$ structures, it may be that
no potential $t_2$ is found. This means the process at this level continues
until it loops. The answer to $Q^2$ is then ``No''. Setting the new $\sigma'$
to be the length of this loop, then ``No'' and $ \la t_0, t_1 
\ra$, and $W_{\sigma'}$ are passed back up to $\Lambda = 1.$

{\bu} (3) Lastly if a $t_2$ is found so that $\tilde{\varphi} \left( 
\la t_0, t_1 , t_2 \ra \right)$ holds in some structure $^0 M$ of
height $\sigma'$, then on the first such occasion a query $Q^3$ is passed down
to a new node at level $\Lambda = 3$:

\

$Q^3 :$ \ {\tmem{Does $\ptwoe{e_0} \left( \pa{3, \nocomma \pa{t_0, t_1,
t_2}, W_{\sigma'}} \right)$ find a stable candidate for $t_3$ starting
from $\pa{t_0, t_1, t_2}$ and the ordinal $\sigma' $ so that $\tilde{\varphi}
\left( \la t_0, t_1, t_2, t_3 \ra  \right)? $}}

\

And so forth. We trust the reader can see the successive formulations of
$Q^k$ {\etc} and the actions of $P^{\mathsf{E}}_{e_0}$ so intended. One
property of the construction is immediate:

\begin{lemma}
  The above construction has that for limit $\alpha$, $\Lambda (e_0, \alpha) =
  k > 0$ implies that in an interval $(\delta, \alpha)$ $t \rest k$ is stable.
\end{lemma}

{\pf}For otherwise, unboundedly in $\alpha$ control would have been passed up
to a level $\Lambda \leq k - 1$, and then by the liminf rule $\Lambda (e_0,
\alpha) < k$.{\qed}

\begin{lemma}
  \label{Lem4.2} If $(\delta, \sigma)$ is a $\Sigma_2$-extendible pair, with
  $\sigma$ still a defined stage in the computation of $\{ e_0 \}^{\mathsf{E}}
  (0)$, then $\Lambda (e_0, \delta) = \Lambda (e_0, \sigma) > 0$. 
\end{lemma}

{\rem} Being `still a defined stage' here means that $\{ e_0
\}^{\mathsf{E}} (0)$ has not crashed at any stage $\sigma' \leq \sigma$, that
is, $\mathfrak{{T^{\mathsf{E}}} } (e_0, 0)$ is still wellfounded below the node at which control of the computation is located, at stage $\sigma$.

\

{\pf} This is Lemma \ref{L4.7}: for some admissible $\tau <$ $\delta$, a
$t_0$ satisfying $\varphi \left( t_0, \emp \right)$ is defined in $L_{\tau}$,
and is stably so defined up to $\sigma$. Supose  first that $(\delta, \sigma)$
is a simply isolated extendible pair, {\ie} we have $\delta \in E^0$. By our
program, then, $\Lambda (e_0, \gamma) > 0$ for $\gamma \in (\tau, \sigma)$. (In
other words there is a query to $\Lambda = 1$ to find a $t_1$ underway during
these ordinal stages.) Thus we must have $\delta \in E^1$. Suppose for a
contradiction that $\Lambda (e_0, \delta) = \Lambda (e_0, \sigma) = 0$. Then
as $\Lambda (e_0, \delta) = \tmop{Liminf}_{\gamma \imp \delta} \Lambda (e_0,
\gamma)$, unboundedly in $\delta$ there are stages where control has passed to
$\nu_0$ from after a finished loop at level 1.

However, this does not happen unboundedly in $\delta$: suppose $(\delta',
\sigma')$ is an extendible pair with $\tau < \delta' < \sigma' < \delta$, and
$\delta' \in E^0 \back E^1$.\quad As $t_0$ is defined at $\tau$ control then
passes (or has already passed) to $\Lambda = 1$; repeatedly admissibles greater than $\tau$ are
checked for candidates $t_1$ for an extension to $t_0$. Note: as $t_0$ s stably defined up to $\sigma$,   control does not pass up to $\Lambda =0$ in $[\tau,\s]$ because $t_0$ fails its $\Pi_1$ check. But $\Upsilon (\vp)$
guarantees that for some $\tau' < \delta$, for larger ordinals, \ and \
$(\delta', \sigma')$ \ any extendible pair with $\tau' < \delta' < \sigma' <
\delta$, and $\delta' \in E^0 \back E^1$, we have some $t_1$ that is ``good up
to $\delta'$'', that is $\varphi (t_1, (t_0))^{L_{\delta'}}$ (and so  ``good up
to $\sigma'$'', that is $\varphi (t_1, (t_0))^{L_{\sigma'}}$ too). So at some
admissible at which such a suitable $t_1$ appears, which is less than
$\delta'$, control passes down to $\Lambda = 2$, for a search for a suitable
$t_2$ etc. This search continues at $\Lambda = 2$, at least until $\sigma'$.
After the latter, control then passes up to $\Lambda = 1$. This situation
holds for all extendible pairs $(\delta', \sigma')$ such as above, and also
for such nested inside $(\delta, \sigma)$. In other words, between $\tau$ and
$\sigma$, $\Upsilon(\vp)$ ensures that at any extendible pair $(\delta',\sigma')\subset (\tau,\sigma)$ control passes {\em downwards} to  lower levels  $\Lambda > 1$  rather than upwards to $\Lambda =0$, because $(t_0,t_1)$ is the current stable choice
in this interval $(\tau',\s')$. {\qed}

\begin{note}
  \label{exhaust}{ We can generalize the argument at the end of the last lemma to characterise the behaviour more broadly. In defining our process $\{e_0\}$ we also analysed the case that a loop occurred in
  which no appropriate extension of the current approximation was found going from one
  $\,^0 M$-structure to the next and the response to the calling query was thus negative.
  The alternative was phrased as ({\cf}the final case (3) when $\Lambda = 1 =
  k$ above) as $t_k$ being defined in a structure. One might have asked:
  `Where is the analysis that a loop occurred at level $k$ in which all of $t
  \rest k\smallfrown (t_k)$ was stable throughout for some $t_k$, and thus the answer to $Q^k$ would have been
  ``Yes'', (as always at the end of a loop, this would be followed by control
  passing up to $\Lambda = k - 1$)?' The following lemma shows that this
  behaviour does not occur.}
\end{note}

\begin{lemma}
  \label{L4.2}For $0 < k$, if $t \rest k$ is stably defined up to
  $\bar{\Sigma}$ (meaning $\tilde{\varphi} \left( t \rest k \right)^{L_{\tau}}
  $ for all sufficiently large $\tau < \bar{\Sigma}$), and $(\bar{\zeta},
  \bar{\Sigma})$ is a loop at $\Lambda = k$ resulting from a call from
  $\Lambda = k - 1$ to answer the query ?$Q^k$ {\tmem{Does
  $\ptwoe{e_0} \left( \pa{k,  t \rest k, W_{\sigma'}} \right)$
  find a stable candidate for $t_k$ starting from $t \rest k$ and the
  ordinal}} $\sigma' $?, then the response to $Q^k$ is ``No''.
\end{lemma}

{\pf} The short answer is simply that if the response were ``Yes'', with $t_k$
a stable candidate, and so that $\tilde{\varphi} \left( t \rest k \smallfrown
(t_k) \right)^{L_{\tau}} $ held for all sufficiently large $\tau <
\bar{\Sigma}$, then for all sufficiently large $\tau < \bar{\Sigma}$ (and then
by $\Sigma_2$-reflection, also for all sufficiently large $\tau <
\bar{\zeta}$) control would be at a level $\Lambda (\tau) = m > k$, looking
for some extension of $t \rest k \smallfrown (t_k)$, and thus $(\bar{\zeta},
\bar{\Sigma})$ is a loop at some level $m \neq k$, contrary to hypothesis. \qed

\begin{note}
  As the program runs there will eventually be subcomputation calls to
  arbitrary levels, as it uses approximations for as long as they survive
  fulfilling their role. But only after $\bar{\alpha}_0$ stages shall we be
  certain that $t_0$ really does stabilize to its final value. Thereafter we
  shall always have $\Lambda (e_0, (0, \emp, W_{\omega}), \alpha)
  \geq 1$. But only at $\beta_0$ shall we first have $$\mathrm{Liminf}_{\alpha
  \rightarrow \beta_0} \Lambda (e_0, (0, \emp, W_{\omega}), \alpha)=
  \omega$$ and so an illfounded computation.
\end{note}

We make some further observations on the flow of control during the recursion.

\begin{lemma}
  \label{extendibility-levels}  Let
  $(\zx{k}, \sx{k}) < \beta_0$ be a $k$-extendible pair. At times during
  $(\zx{k}, \sx{k})$ control will pass to depth at least $\Lambda = k + 1$,
  with queries $Q^{k + 1}$ asking that $t \rest k + 1$ be stable and seeking a
  candidate for $t_{k + 1}$. 
 \end{lemma}

{\pf}Formally a proof by induction on $k$, but really this is just a statement
on the template of the construction. 
Note first that for $k = 0$, if $(\xi,
\sigma) = \left( \zx{0}, \sx{0} \right)$ is $0$-extendible, but not
$\gamma$-extendible for any $\gamma > 0$, then by Lemma \ref{L4.7} there is $\delta
< \xi$ so that already $\Lambda (e_0, \alpha) \geq 1$ for $\alpha \in (\delta,
\sigma)$.

Consider the case $k=1$ and  that $\left( \zx{1}, \sx{1} \right) \supset (\zx{0}, \sx{0})$.
By $\Upsilon (\varphi)$ there is a $t_0 \in
L_{\zx{1}}$ satisfying the $\Pi_1$ condition $(\varphi (t_0,t\rest 0))^{L_{\zx{1}}}$.
As $L_{\zx{1}} \prec_{\Sigma_1} L_{\zx{0}}$, $(\varphi (t_0))^{L_{\zx{0}}}$
also. But $\Upsilon (\varphi)$ implies that there is $t_1 \in
L_{\zx{0}}$ such that ${\varphi} (t_1,  t\rest 1)^{L_{\zx{0}}}$. By
our description of $P^{\mathsf{E}}_{e_0}$ then a query $Q^2$ to $\Lambda = 2$ is
launched asking if $t_0 \smallfrown t_1$ is stable and seeking a $t_2$. (And
actually by $\Sigma_2$-reflection such a query or queries are being acted out
at level $\Lambda = 2$ unboundedly in $\zx{1}$ also.) {\qed}

\begin{corollary}
  \label{C4.6}If $\zeta_0 < \zeta_1 < \cdots < \zeta_k < \alpha < \sigma_k <
  \cdots < \sigma_1 < \sigma_0 < \beta_0$ is any $k + 1$-nesting then for some
  $\delta < \zeta_k$, $\Lambda (e_0, \gamma) \geq k +1$ for $\gamma \in (\delta,
  \sigma_k)$. 
\end{corollary}

\begin{lemma}
  \label{finalloop}   
   Our overall computation can never reach a final looping pattern.
\end{lemma}

{\pf} Suppose otherwise and that $(\nu, \sigma)$ is the first final looping period of
the master program. Then $\Lambda (e_0, \nu) = \Lambda (e_0, \sigma) = 0$.
Then $(\nu, \sigma)$ is a $\Sigma_2$-extendible pair, (as we have bookkept all
the theories $T^2_{\alpha}$ for $\alpha \leq \sigma$). However,  then by Lemma \ref{Lem4.2} $\Lambda (e_0,
\nu) = \Lambda (e_0, \sigma) > 0$. Contradiction!{\hspace*{\fill}} {\qed}

\begin{theorem}
  \label{leastnonrecursive}   
  $\alpha_0$ is the least ordinal which is not ITTM generalised recursive.
\end{theorem}

\pf  We just show that the supremum of such ITTM generalised recursive  ordinals is at least $\alpha_0$ (and let the reader fill in the other gaps). Let $\psi \in T^1_{\a_0}$. Suppose $\tau_\psi$ is the least ordinal so that $(\psi)^{L_{\tau_\psi}}$. Then we may run an adjusted version of the computation above which additionally looks for a level of the $L$-hierarchy in which $\psi$ holds, namely this $\tau_\psi$; and then halts when it locates it. This is a process which will converge in $\geq \tau_\psi$ steps. And such $\tau_\psi$ are unbounded in $\a_0$. However the arguments above show that no computation of the form $\{e\}^\me(m)$ converges at any stage in the interval $(\a_0,\b_0)$, and will crash at time $\b_0$ if it has not converged by stage $\a_0$.  Hence $\alpha_0$ is the exact supremum here. \qed

\begin{theorem}\label{H=Psi} $H^\me \equiv_1 \Psi$.
\end{theorem}
\pf This is a corollary to the argument of the last theorem: $L_{\a_0}$, and so $\Psi$, contains all the convergence information about computations of the form $\{e\}^\me(m)$, and so $H^\me \leq_1 \Psi$. However, as in the argument above, for any $\psi\in \Psi$ we have a computation that searches for $\tau_\psi$ and then halts. If we know $H^\me$ then we shall know for which $\psi \in \S_1$ this holds. \qed

%%%%%%%%%%%%%%
\section{$\Sigma^0_3$-Determinacy and $\omega$-nestings}

We consider here a classical application of the previous analysis of generalised ITTMs and infinite nestings to determinacy at the $G_{\d\s}$, or $\Sigma^0_3$, level. The games played are the usual two person perfect information infinite Gale-Stewart games with moves as integers creating sequences in a specified {\em game tree} which is simply a tree of finite sequences $T\sset ^{<\w}\w$ closed under initial segments growing from the empty sequence root. See, for example, \cite{Mosch4}.
 
\begin{definition}
  Let $\delta_0$ be the least ordinal so that for any game $G (A, T)$ with
  $A \in \Sigma^{0}_{3}$, $T \in L_{\delta_0}$ a game tree, then there is a
  winning strategy for a player definable over $L_{\delta_0}$.
\end{definition}

\nod It is our task to characterise $\delta_0$ in terms of the existence of infinitely nested ordinals.

In the next section we show that with $\beta_0$ as the first infinitely nested ordinal, that $L_{\beta_0}$ is not a model of $\Sigma^0_3$-Determinacy. Hence that $\d_0\geq\b_0$. We do this by showing that there are $\Sigma^0_3$ games which player \po wins, but winning strategies for such games must appear cofinally in $\a_0$. Were $L_{\b_0}\models $``$\Sigma^0_3$-Determinacy'' to hold, this would be a $\Sigma_1$ statement that would reflect down to $L_{\a_0}$, which would be impossible.

The subsequent section shows that $\d_0\leq\b_0$, by demonstrating how strategies for \po in such games may lie in $L_{\a_0}$ but those for \pt are at worst $\Delta_2(L_{\b_0})$-definable.

\subsection{$\Sigma^0_3$-Determinacy fails  below the least  $\omega$-nested ordinal}

Recall the following definition.

\begin{definition}
  Let $\Gamma$ be a pointclass. A set $Z \subseteq \mathbbm{N}$ is said to be
  in $\Game \Gamma$ if there is a set $X \subseteq \mathbbm{N}
   \times
  \mathbbm{N}^{\mathbbm{N}}$ in $\Gamma$ so that
  $$Z =\{n\mid \tmtextit{I} \mbox{ has a winning strategy in } G (X_n ; \omega^{<
  \omega}) \}$$ $ \mbox{ where }X_n =\{y\mid \langle n, y \rangle \in X.$\\
\end{definition}
%%%%%%%%%%%%%
\begin{theorem}\label{Friedman}
Let $\Psi =\{\psi | \psi \in \Sigma_1\cap Sent \wedge L_{\beta_0} \models \psi\}$ be
$  T^1_{\beta_0}$, the $\Sigma_1$-theory of
  $L_{\beta_0}$. 
  Let $G_3$ be a complete 
$\Game
  \Sigma_3^0 $ set of integers.  Then   $\Psi \leq_1 G_3 $ and is $\Game
  \Sigma_3^0 $.

\end{theorem}

 \pf  We do this, using a technique that goes back to H. Friedman, by defining
  certain games $G_{\psi}$ so that codes for initial segments of the
  $L$-hierarchy are recursive in any winning strategy for the game.

   For any $\psi
  \in \Psi$ we define: $\alpha_{\psi} =$ the least $\beta$ so that $L_{\beta}
  \models \tmop{KP} + \psi$. 

 \nod Note:  \label{skolemterms} The minimality of $\alpha_\psi$ ensures that every $x \in
    L_{\alpha_{\psi}}$ is $\Sigma_1$-definable by some parameter free
    $\Sigma_1$ term $t_x$. (In other words the $\Sigma_1$-Skolem hull inside
    $\langle L_{\alpha_{\psi}}, \in \rangle$ of $\varnothing$ is all of
    $L_{\alpha_{\psi}}$ itself.)\\

  The following is straightforward:\\

  (2) Let $\overline{\alpha} = \sup \{\alpha_{\psi} | \psi \in \Psi\}.$ Let
  $\alpha_0$ = the least $\beta ( L_{\beta} \prec_{\Sigma_1} L_{\beta_0}) .$
  Then $\alpha_0 = \bar{\alpha} .$\\

  We shall show for every $\psi \in \Sigma_1$ there is a game $G_{\psi}$ with a
  {\tmem{$\Sigma^0_3 ${\tmem{}}}} payoff set, for which \po has a winning strategy in $L_{\alpha_0}$
   if and only if $\psi\in \Psi$. This will show  that  $\Psi \leq_1 G_3$, and hence, by Tarski, that $G_3\notin L_{\alpha_0}$.   This will use the `smallness' assumption that no $\bar \beta < \beta_0 $ supports an $\omega$-nesting.\\

  For the rest of the argument fix a $\psi \in \Sigma_1$. Let {\tmem{$\alpha$}}
  denote $\alpha_{\psi}$. We consider the following game $G = G_{\psi} .$

  $I$ \ \ \ plays \ $m_{0,} m_{1,} \ldots, m_i$ \ \ \ \ \ \ \ \ \ \ \ \ \ \
  \ \ \ \ \ \ \ \ $x = (m_{0,} m_{1,} \ldots, m_i, \ldots)$

  {\tmstrong{}}{\tmem{II}}  \ \ plays \ \ \ \ \ \ $n_{0,} n_{1,} \ldots, n_i$ \ \  \ \ \ \
  \ \ \ \ \ \ \ \ \ \ \ \ \ \  \ $y = (n_{0,} n_{1,} \ldots, n_i,
  \ldots)$

  {\noindent}in the usual way, playing in the $i$'th round integers $(m_i,
  n_i) .$ Let $z = x \oplus y.$\\

  {\tmem{Rules for \em {\pt}}}

  Let $\text{$T$}$ be the theory ``$\tmop{KP }\, + \,{ V = L} \,+ \,\ex \bar\beta[ \bar \b \mbox{ \em  supports an $\w$-nesting} \wedge (\neg \psi)^{L_{\bar\beta}}]$". \pt's play $y$
  must be a set of G\"odel numbers for the complete $\Sigma_1$-theory of an
  $\omega$-model of \ $T +$``{\tmem{there is no set model of}} $T$''. Note  then if $T$ has a transitive model it is isomorphic to $L_{\b_0^+}$, the least admissible level of $L$ beyond $\beta_0$.

  The Note above on Skolem terms applies also to $\la L_{\b_0^+},\in\ra$: $\b_0^+$ is $\S_1$-definable, and every set is countable; hence every $x\in  L_{\b_0^+}$ is of the form $h^1_ {{\b_0^+}}(m,0)$ for some $m$. We denote by $\langle M, E \rangle$ the
  model \pt essentially constructs if he obeys this rule. We may regard also
  as part of the rule that $y$ as given by \pt should be specified simply by
  \pt stating ``$k \in T^1_M$'' or `$\text{`$k {\notin} T^1_M$''}$ where $T_M
  = T^1_M$ is the standard $\Sigma_1$-code or `truth-set' for
  his model.
  Then \red, just as in $L_{\alpha_{\psi}}$, \ered in $\langle M, E \rangle$,
  every set is likewise given by a $\Sigma_1$ parameter free Skolem term.
Consequently, as noted above, if $\langle M, E \rangle$ is wellfounded then it is isomorphic to $\langle
    L_{\b_0^+}, \in \rangle$.

  Amongst the codes for sentences that \pt plays are those of the form
  $$\ulcorner t_m \in \tmop{On} \wedge\ \  t_n \in \tmop{On} \wedge \, \ t_m < t_n
  \urcorner$$
  These we shall use to formulate rules for player \po. So far the
  Rules for \pt amount to a $\Pi^0_2$ condition on $y$ and so on $z$. (We may
  take a recursive listing of $\Sigma_1$-sentences $\langle \psi_k |k \in
  \omega \rangle$ and we then require $\forall k \exists k' (m_{k'} =
  \ulcorner \psi_k \urcorner \vee m_{k'} = \ulcorner \neg \psi_k \urcorner)$,
  thus the theory \pt constructs will be $\Sigma_1$-complete; we
  obtain that $M$ has at least the integers as standard, this is also by a $\Pi^0_2$
  condition.) \ Let $r : \omega \rightarrow \omega \times \omega$ be a
  recursive enumeration of $\omega^2$ in which each $(i, j)$ appears
  infinitely often.\\

  {\tmem{ Rules for }\po.}

  At round $k$:
if $(i, j) = r (k)$ and $n_k \neq 0,$ then we
  shall say that ``\po {\tmem{ makes the entry $n_k$ on list}} $L_{i,
  j}$.''. These `Listing' {\tmem{Rules}} here require her to list terms in a
  correct order. She may make an entry on list $L_{i, j}$ in round $k$ if:

  {\tmem{Either}} $L_{i, j}$ is empty at the current round, in which case
  $n_k$ can be any term $t_s$ {\tmem{as long as}} {\po} has asserted at
  an earlier round $\ulcorner t_s \in \tmop{On} \urcorner \in T_M$;

  {\tmem{or}} $L_{i, j} \neq \varnothing$, and if $t_s$ was the last entry
  {\po}  made on this list, then $n_k$ can be any term $t_r$, again
  {\tmem{provided that}} \pt has at an earlier round $k' <
  k$ asserted $m_{k'} = \text{ } \text{$\ulcorner t_r \in \tmop{On} \wedge\,  t_r
  < t_s \urcorner$} \in T_M$.\\

  {\tmem{The winning conditions.}} \pt wins immediately at a finite
  round if $\po$ breaks one of her  {\tmem{Listing Rules}} just
  enumerated.  $\po$ wins if \pt fails to obey his conditions on $y$,
  or both players obey their respective rules and additionally

 \quad\quad $\exists (i, j) [$\po {\tmem{makes infinitely many entries on
  list}} $L_{i, j} ]$.

  \nod This is a $\Sigma^0_3$ winning condition for \po on $z$. Hence
  $G_{\psi}$ has a $\S^0_3$ payoff set.\\

%  {\bf Remark:}
    In other words, if \pt obeys his rules, $\po$ can win if for some
    $ (i, j)$, $r^{- 1}$``$ (i,
    j)$ in effect picks out an infinite descending chain through the
    ordinals of $\tmop{the}$ model ${M}$ that \pt reveals {\tmem{via}}
    the g\"odel numbers of the $\Sigma_1 $sentences true in ${M}$.\\

 % {\bf Remark:}
    \po is not allowed to make an entry indicating that $t_s < t_r$
    until \pt has asserted this at some earlier stage. $\po$ is thus
    not predicting what the model will look like below $t_r$; by making an
    entry on a list she is merely adverting to the fact that \pt has already revealed
    that $\text{} t_s < t_r$.\\

(3)  Suppose $\psi\notin\Psi$. Then  \pt has a winning strategy.

  {\noindent}{\tmstrong{Proof:}} \pt plays out all ``$k \in T_M$'' for all $k
  \in T^1_{\beta_0^+}$, and ``$k {\notin} T_M$'' for all $k {\notin}
  T^1_{\b_0^+} $. Obviously then, $\langle M, E \rangle \simeq
  \langle L_{\b_0^+}, \in \rangle$. As $(\neg\psi)^{L_{\a_0}}$, by $\S_1$-elementarity,
  $(\neg\psi)^{L_{\b_0}}$. As $E$ is wellfounded
  \po has no
  chance to pick out any infinite descending chains. \hfill Q.E.D. (3)\\

  The point is the following:\\

   (4) \label{mainlemma}     Suppose $\psi\in\Psi$. Then  $\po$ has a winning strategy in $G_\psi$.

\nod  \tmtextit{}From this it then follows that $\tau \notin
  L_{\alpha_{0}}$, as otherwise this would imply that the latter's $\Sigma_1$-truth set is a member of itself, which would be a contradiction.\\

  {\noindent}{\tmstrong{Proof of (4)}
Since $\psi\in \Psi$ again by upwards persistence of $\S_1$ we have that $\psi$ holds in $L_{\b_0^+}$, hence \pt cannot play out a wellfounded model isomorphic to this level of the $L$-hierarchy. To do this would be to break the basic $\Pi^0_2$ rule of \pt, and so \po would win. So we may assume that \pt plays an illfounded $\pa{M,E}$.

By the same considerations, we cannot have that $WFD(M,E)\supseteq L_{\a_\psi}$ since $(\psi)^{L_{\a_\psi}}$. 
 Hence, by the Truncation Lemma, $WFD(M,E)= L_{\b_M}$ for some admissible $\b_M<\a_\psi$.
 
    \po assumes that \pt's
  eventual model $\langle M, E \rangle$ will be illfounded, and so she must
  act to discover a descending chain. 
  
   However she will not
  yet know, and in fact will not at any move know, where $\beta_M$ lies. All she will know is that
as  $(\tmop{KP})^M $ (if \pt plays correctly) that $\beta_M \in \tmop{ADM}$. 
   By our requirements on the theory
  $T^1_{\alpha_{\psi}}$, and upwards persistence of $\Sigma_1 $ formulae, we
  must have $\beta_M < \alpha_{\psi}$.

  \begin{definition}
    Let $F : \omega \twoheadrightarrow \tmop{ADM} \cap \ \alpha_{\psi} + 1$ be
    some fixed surjection.
  \end{definition}

  The idea is that at rounds $k$ where $r (k) = \text{$\text{$\text{}$} (i,
  j)$}$ \po will be making the working assumption that the ordinal
  height of the wellfounded part of $M$, $\beta_M,$ is precisely $F (i)$, and
  will be trying to find an illfounded chain through $\tmop{On}^M$ above
  $\beta_M .$ She will be working simultaneously on all such possible
  $\beta_M$.

  However we concentrate our description on an arbitrary but fixed $i$
  and hence on a fixed $\beta = \beta_M = F (i),$ and describe how \po can
  move in rounds $k$ with $r (k) = \text{$\text{$\text{}$} (i, j)$}$.\\

  (5) {\tmem{Claim}} $\exists \bar{a} \notin \tmop{WFP} (M) \forall b <
  \bar{a} (b \notin \tmop{WFP} (M) \rightarrow T^2_b \not\subset
  \widetilde{T} =_{df} T^2_{\beta}$.\\

  {\noindent}{\tmstrong{Proof.}} Supposed this failed, then $\forall \bar{a}
  \notin \tmop{WFP} (M) \exists b < \bar{a} (b \notin \tmop{WFP} (M) \wedge
  T^{2_{}}_b \subset \widetilde{T}).$

  Let $\eta =\eta_b =_{\tmop{df}} \sup \{c < b \mid \exists f \in \Sigma^{J_b}_2 , f
  : \omega \twoheadrightarrow c, f$ partial, onto$\}$.
  Note also, for use in a moment, that if $\{c\}$ is any $\Sigma_2^{J_b}$ definable ordinal
(that is, defined without parameters) then $c<\eta$, as the $L$-least onto map $f:\omega \rightarrow c$
is then also $\Sigma_2^{J_b}$ definable (it lies in $J_{c+1}$).
   We first claim that
  $\eta < \beta$. Clearly equality fails, as otherwise that would make $\beta$
  definable inside $M$ from $b$. If however $c \notin \tmop{WFP} (M)$, with $f
  \in \Sigma^{J_b}_2$, $f$ partial, but onto $c$, then the
  sentences ``$f (n) \downarrow, f (m) \downarrow \wedge\, f (n) < f (m) \in
  \tmop{On}$'' are all in $T^2_b$ and so in $\widetilde{T}$. This is absurd as
  $\beta$ is wellfounded! Hence $\eta<\beta$. Note that this somewhat trivially
  implies that $b$ is
  an $<_M$-limit ordinal: were $b=b_0 +1$ then $b_0$ itself is $\Sigma_2^{J_b}$
  and by the above reasoning we'd have the absurdity $b_0\in \tmop{WFP(M)}$!

It is not hard to see that $\eta$ is closed under the G\"odel pairing function
and this implies that there is a parameter free $\Delta_1^{J_\eta}$ bijection $\eta \leftrightarrow
J_\eta$ (\tmem{cf} \cite{De}).
Suppose
  $J_b \models \exists u \psi (u, \xi )$ where ${\xi} < \eta$ and
  $\psi \in \Pi_1$. (It suffices to verify $\Sigma_2$-elementarity just on
  formulae with single ordinal
  parameters $\xi$ by the above remarks.)

Let $\delta$ be the least ordinal such that $J_b\models$``$\forall
  \delta' > \delta J_{\delta'} \models \exists u \psi (u, \xi )$.'' Then $\{\delta\} \in
  \Pi_1^{J_b}(\{\xi\})$. There is thus a $\Sigma_2^{J_b}(\{\xi\})$ partial map $f_\delta:\omega
  \twoheadrightarrow \delta$ given by some formula: $f_\delta(m)=\tau \leftrightarrow
  \exists w \chi(w,m,\tau,\xi)$ for a $\Pi_1$ $\chi$.

  As ${\xi} = f_0 (n)$ for some $\Sigma^{J_b}_2$ $f_0$,  we have
  $$f_\delta(m)=\tau \leftrightarrow
  \exists x[x=f_0(n)\wedge \exists w \chi(w,m,\tau,x)].$$  Replacing ``$x=f_0(n)$''
  with its $\Sigma_2$ definition, this yields a parameter free $\Sigma_2^{J_b}$ definition
  of $f_\delta$. Hence $\delta <\eta$. By the definition of $\delta$ we shall have
  $J_\eta\models \exists u \psi (u, \xi )$ as required.

  Hence for such a $b$ we have $( J_{\eta_b} \prec_{\Sigma_2} J_b)^M$.
  However the supposition implies there is an infinite descending chain of
  such $b$ in the illfounded part of $M$. This implies that we have an
  infinite nested sequence of $\Sigma_2 $ reflecting intervals: there exists $\langle b_n
  |n < \omega \rangle$, $\langle \eta_n |n < \omega \rangle$ with $(\eta_n
  \leq \eta_{n + 1} \leq \ldots \beta \leq\ldots < b_{n + 1} < b_n)$, and with $(J_{\eta_n}
  \prec_{\Sigma_2} J_{b_n})^M$, for $n < \omega$.  However this would mean that $\beta_{M}$ supported an infinite nesting. But $\beta_{M}<\alpha_{\psi}<\beta_{0}$ where the latter is least supporting an infinite nesting. Contradiction! \hfill \qed\, (5)\\

   Let $\langle t_k |k \in \omega \rangle$ be our priorly fixed recursive
  enumeration of the $\Sigma_1$-Skolem terms (using the standard
  $\Sigma_1$-Skolem function, this could simply be an enumeration of  $\langle
  h^1 (i, n) \mid i, n < \omega \rangle$). \po makes the additional
  working assumption, or guess if you will, that $t_j^M = a_0,$ where $a_0$ is
  a witness for $\bar{a}$ to the truth of the last Claim. (Again the point is
  that \po does not know in advance which term in $M$ will denote
  such $a_0$.) {As} \pt reveals more and more facts about his model, he
  must, if $M$ is not going to be isomorphic to $L_{\b_0^+}$ at some
  point reveal a $\Sigma_1$-fact which is true in $M$ but false in
  $L_{\a_\psi}$. There really is then such an $M$-ordinal $a_0$.
  \po will, in effect, place her `guess'  $a_0 = t^M_j $ \ at the
  head of her putative descending chain, on list $L_{i, j}$, when round $k$ first satisfies $r(k)=(i,j)$. In order to
  choose the next element of the chain on this list, \po considers the set
  $\widetilde{T} = T^2_{\beta}$. Set $T_0 = (T^2_{t_j})^M$.

 \po now waits until \pt asserts that some $\sigma_0 $ is in $T_0$,
  (this itself being one of the $\Sigma_1$ facts about $M$ she must enumerate)
  but {\po} sees is not in $\tilde{T}$. (If \po is
  wrong in her guess about $t_j$ of course, then she may fruitlessly wait for
  ever.)\\

  (6) {\em Suppose $M \models$``$a_1 < a_0$ is least so that $\forall b \leq a_0 (b
  \geq a_1 \rightarrow (\sigma_0)_{L_b})$.'' Then } $a_1 \notin \tmop{WFP} (M)$.\\

  {\noindent}{\tmstrong{Proof}}: Were $a_1 \in L_{\beta}$ then we should have
  $\sigma_0 \in \widetilde{T}$. \hfill Q.E.D. (6)\\

  \po may thus wait until \pt asserts that some such $\sigma_0 \in T_0
  \backslash \widetilde{T}$ and additionally, perhaps later, the $\Sigma_1$
  fact that some term $t_{j_1}$ names the ordinal $a_1$ defined in (6) above.
  At some round $l$, \pt must then play the number 
  $$m_l =
  { \ulcorner t_{j_1} \in \tmop{On} \wedge\, t_j \in \tmop{On} \wedge
  t_{j_1 } < t_j \urcorner ;}$$ once all these facts have been gathered
  together, \po may at the next appropriate round $k$ with $r (k) = (i,
  j)$, set $n_k = t_{j_1}$.

  \po now has two elements of a descending chain in the illfounded
  part of $M$. Now she watches out for assertions that \pt makes about $T_1 =
  (T^2_{t_{j_1}})^M$, waiting for some $\sigma_1$ asserted by $\tmop{him}$ to
  be in $T_1$ but which does not lie in $\widetilde{T}$. By exactly the same
  considerations that held at (6) some $a_2, t_{j_2},$ are definable, and so
  she can continue. {\tmem{If}} this working
  assumption about {\tmem{{\tmem{{\tmem{$\beta_M $}}}}}}and $t_j$ was the
  correct one, by the end of the game the chain so defined by continuation of this process will be
  infinite, and she will have won.

  If \pt deviates from playing the correct $\Sigma_1$ truth set for $L_{\b_0^+}$, then at least
  one of \po's assumptions will turn out to be a correct one, that particular list will be infinite and
  hence she will be assured of winning.
  
  Notice finally, that if $\po$ wins $G_\psi$, then a winning strategy for \po 
 was definable just from knowledge of $T^1 _{\a_\psi}$, and so such a strategy can be found in $L_{\a_\psi^+}$.
   \mbox{ }\qed
  ({\tmstrong{(4)\& Theorem \ref{Friedman}}})\\

 Let $\sigma_3$ denote the least $\sigma$ so that every $\Sigma_3^0$ game that is a win for $\po$, has a winning strategy in $L_{\sigma}$.
\begin{corollary}
  \label{sig} (i) Each $T^1_{\alpha_{\psi}}$
  is in $\Game \Sigma_3^0 \cap \Game \Pi_3^0$ as a set of integers.\\
  (ii) ${\alpha}_0 \leq \sigma_3$.
\end{corollary}

{\noindent}{\tmstrong{Proof:}} The arguments are just variants of the above. Let $\alpha_{\psi}$ {\tmem{etc.}} be defined as
above.

Fix a $\psi\in \Psi$. We first show that $T^1_{\alpha_{\psi}}\in \Game \Sigma_3^0$. For $\varphi \in \Sigma_1$ let $G_{\psi,
\varphi}$ be the game described in the last theorem, except that \pt
must now play a code $y$ for a model of \ $T +$``there is no set model of
$T$''$+(\neg \varphi)^{L_{\a_\psi}}$. Everything else remains the same {\tmem{mutatis
mutandis}}: \po's task is still to find an infinite descending chain
through the ordinals of \pt's model. Note that if $\varphi \in
T^1_{\alpha_{\psi}}$ \po now has a winning strategy: for if
\pt obeys his rules, and $y$ codes an $\omega$-model $M$ of this
theory, then $M$ is not wellfounded, and has WFP$(M) \cap \tmop{On} <
\alpha_{\varphi}$, where the latter is the least admissible $\alpha$ where
$\varphi$ is true in $L_{\alpha}$. However \po playing can find a descending chain and win. Now, on
the other hand if $\varphi \notin T^1_{\alpha_{\psi}}$ \pt may just
play a code for the true wellfounded $L_{\b_0^+}$ and so win. This
shows that $T^1_{\alpha_{\psi}}$ is a $\text{$\Game \S_3^0$}$ set of
integers.\\

\nod $ T^1_{\alpha_{\psi}}\in \Game \Pi_3^0$: just replace ``$(\neg \varphi)^{L_{\a_\psi}}$'' by ``$(\varphi)^{L_{\a_\psi}}$" in the model $\pt$ must play, and reason similarly with $\pt$ winning iff $\vp\in T^1_{\a_\psi}$.\\

Suppose now that
% \label{sig}
%\marginpar{Why was there a label sig here?}
${\alpha}_0 > \sigma_3$. Let $\psi$ be such that
$\alpha_{\psi}$ is the second least admissible greater than $\sigma_3$. There
is thus a set $H \in L_{\alpha_{\psi}}$ (definable over the first admissible
$\gamma > \sigma_3)$ containing winning strategies for all
$\Sigma^{0_{}}_3$-games that are a win for player \po and in
particular a set $H_0$, definable at the same level,
 of those winning strategies for \po in games
of the form \ $G_{\psi, \varphi}$. Hence membership of $\varphi$ in
$T^1_{\alpha_{\psi}}$ is determined by searching through $H_0$ for a winning
strategy for \po; this is a bounded search. Hence
$T^1_{\alpha_{\psi}} \in \Delta_1^{L_{\alpha_{\psi}}} (\{H_0 \})$. Hence
$T^1_{\alpha_{\psi}} \in L_{\alpha_{\psi}}$ which contradicts Tarski. 
\hfill \qed

\begin{corollary}
   Each real $r\sset \w$ in $L_{\alpha_0}$ is 
  $\Game \Sigma_3^0 \cap \Game \Pi_3^0$.
\end{corollary}

\begin{corollary} \label{sig2} Assume $\S^0_3$-$Det$. Then 
   each real $r\sset \w$ in $L_{\alpha_0}$ is
  $\Delta(\Game \Sigma_3^0)$.
 \end{corollary}
\pf $\S^0_3$-$Det$ implies that $(\widecheck{\Game \S_3^0}) = \Game \Pi_3^0$. \qed
%%%%%%%%%%%%

\subsection{$\S^0_3$-Determinacy requires infinite nestings}

\begin{theorem}\label{strategies}
Let $A\in \Sigma^0_3$. Then if {\em\po} has a winning strategy $\sigma$ for $G(A)$ then there is such a strategy in $L_{\beta_0}$. If {\em \pt}has such a strategy $\tau$, then there is such a strategy $\Delta_2$ definable over $L_{\beta_0}$.
\end{theorem}
\pf We outline a proof.
We look at the
construction of the proof of Theorem 5 \ of {\cite{W2011}} in particular that
of Lemma 3. There we used an assumption that there is a triple of ordinals
$\gamma_0 < \gamma_{1} < \gamma_{2}$ with (a) $L_{\gamma_0}
\prec_{\Sigma_{2}} L_{\gamma_{1}}$ and (b) $L_{\gamma_0} \prec_{\Sigma_{1}}
L_{\gamma_{2}}$ and (c) $\gamma_{2}$ was the second admissible ordinal beyond
$\gamma_{1}$. \ One assumed that $I$ did not have a winning strategy in $L_{\g_0}$ for $G
(A;T)$.  Lemma 3 there ran as follows:

\begin{lemma}
  \label{2.12} Let $B \subseteq A \subseteq \lceil T \rceil$ with $B \in
  \Pi^{0}_{2}$. If $(G (A;T)$ is not a win for $I)^{L_{\gamma_0}}$, then
  there is a quasi-strategy $T^{\ast} \in L_{\gamma_0}$ for {\em \ptwo} with the
  following properties:
  
  (i) $\lceil T^{\ast} \rceil \cap B =  \varnothing \text{ }$
  
  (ii) $ ( G (A;T^{\ast} )$ is not a win for $I )^{L_{\gamma_0}}$.
\end{lemma}

For $p\in T$ we let $T_p=\{u\mid u\in T\wedge (p\sset u \vee u \sset p\}$. 
Here a {\em quasi-strategy} for a player, here player \pt, in a game $G(A;T)$ is a subtree $T^\ast\sset T$ which restricts only Player \pt's moves, compared to the original game $G(A;T)$. As is usual we let $\lceil T \rceil =\{x\in \bai\mid \all k (x\rest k \in T)\} \cup \{ u\in T\mid u \mbox{ is maximal sequence of } T\}$ denote the set of branches through $T$.

The format of the lemma's proof involved showing that the
$\Sigma^{L_{\gamma_0}}_{2}$ notion of `goodness' embodied in (i) and (ii)
held for the starting position $\emp$. \ To do this involved defining goodness for positions in general. We first
define $T'$ as {\ptwo}'s {{\em non-losing quasi-strategy\/}} for $G (A;T)$
(the set of positions $p \in T$ so that {\pone} does not have a winning
strategy in $G ( A;T_{p} )$) (where $T_p$ is  that part of the tree below $p$, \ie where all segments extend $p$); this is $\Pi_{1}$ definable over
$L_{\gamma_0}$ as the latter is a model KPI. Then `` $p \in T'$ '' is $\Pi^{L_{\zeta_{0}}}_{1}   $, where
$\zeta_{0} \dfs \min  S^{1}_{\gamma_0} \back \rho_{L} (T)$. More generally
we define: a position $ p  \in  T' $ is {{\em good\/}} if there is a quasi-strategy $T^{\ast}  $
for {{\em II\/}} contained in $T'_{p}  $ so that the following hold:

(i) $\lceil T^{\ast} \rceil   \cap  B =  \varnothing$; 

(ii) $G (A;T^{\ast} )$ is not a win for $I$.

\nod Here $T'_{p}$ is the subtree of $T'$ below the node $p$. \ The point of
requiring that the pair $( \gamma_0 , \gamma_{1} )$ have the
$\Sigma_{2}$-reflecting property of (a) above, is that the class $H$ of good
$p$'s of $L_{\gamma_{1}}$ is the same as that of $L_{\gamma_0}$ and so is a
set in $L_{\gamma_{1}}$ as it is thus definable over $L_{\gamma_0}$ by a
$\Sigma_{2} ( \{ T' \} )$ definition. \ The overall argument is a proof by
contradiction, where we assume that $\emp$ is in fact not good, and proceeds
to construct a strategy $\sigma$ for Player $I$ in the game $G (A;T' )$, which
is definable over $L_{\gamma_{1}}$, and is apparently winning in
$L_{\gamma_{2}}$. (The requirement (c) that $\gamma_{2}$ be a couple of
admissibles beyond $\gamma_{1}$ was only to allow for the strategy $\sigma$ to
be seen to be truly winning by going to the next admissible set, and verifying
that there are no winning runs of play for {\ptwo}.) The contradiction arises
since $T'$ - which was defined as the subtree of $ T$ of {\ptwo}'s non-losing
positions - is concluded still to be the same subtree of non-losing positions
in $L_{\gamma_{2}}$. Being a non-losing position, $p$ say, for {\ptwo} is a
$\Pi_{1}$ property of $p$. This carries up from $L_{\gamma_0}$ to
$L_{\gamma_{2}}$ as $L_{\gamma_0} \prec_{\Sigma_{1}} L_{\gamma_{2}}$, and
this is the reason for the requirement (b): we want $T'$ to survive beyond
$L_{\gamma_{1}}$ for our argument to work. There is then no winning strategy for $I$ in $G (A;T' )$
definable over $L_{\gamma_{1}}$, contradicting the reasoning that $\sigma$ is
such.

This proves the Lemma: $L_{\gamma_{1}}$ sees there is $T^{\ast}$ a subtree of
$T'$ witnessing that $\emp$ is good. The existence of such a subtree is a
$\Sigma_{2} ( \{ T' \} )$-sentence, and then again this reflects down to
$L_{\gamma_0}$. We thus have such a $T^{\ast}$ in $L_{\gamma_0}$.\\

The Theorem is proven by repeated applications of the Lemma, by using the
argument for each $\Pi^{0}_{2}$ set $B_{n}$ in turn where $A= \bigcup_{n}
B_{n}$ and refining the trees using this procession from a tree to a subtree
$T^{\ast}$. We thus repeat the argument with $T^{\ast}$ replacing $ T$. \
Because $T^{\ast} \in L_{\gamma_0}$ we have the same constellation of this
triple of ordinals $\gamma_{i}$ above the constructible rank of $T^{\ast}$,
and can do this.

However we can get away with less. An $\w$-nesting of an ordinal $\b$ will be  just what is needed.  Moreover the ordinal $\b_0<\g_0$ for the $\g_0$ in the triple constellation described above. 

The definition of the subtree of non-losing
positions of {\ptwo} now this time in the new $T^{\ast}$ can be considered as
taking place $\Pi_{1}$ over $L_{\delta_{0}}$ where $\eta_{0}$ is the least
element of $S^{1}_{\gamma_0}$ with $T^{\ast} \in L_{\eta_{0}}$. \ To get
our contradiction we actually use that $L_{\eta_{0}} \prec_{\Sigma_{1}}
L_{\gamma_{2}}$ ; we do not need that $L_{\gamma_0} \prec_{\Sigma_{1}}
L_{\gamma_{2}}$. \ Notice that our argument that $T^{\ast}  $ exists is
non-constructive: we simply say that the $\Sigma_{2}$-sentence of its
existence reflects to $L_{\gamma_0}$: we do not have any control over its
constructible rank below $\gamma_0$. \ Moreover any sufficiently large
$\gamma'$ greater than $\gamma_{1}$ would do for the upper ordinal, as long as
it is a couple of admissibles larger than $\gamma_{1}$. Thus we could apply
the Lemma repeatedly for different $B_{n}$ if we have a guarantee that
whenever a $T_{n}^{\ast}$-like subtree is defined there exists a $\zeta_{n}
\in S^{1}_{\gamma_0}$ and a suitable upper ordinal $\gamma_{n} > \gamma_{1}$
with $T_{n}^{\ast} \in L_{\zeta_{n}} \prec_{\Sigma_{1}} L_{\gamma_{n}}$. Of
course if there are arbitrarily large $\zeta_{n}$ below $\gamma_0$ with this
extendability property, then this is tantamount to $L_{\gamma_0}
\prec_{\Sigma_{1}} L_{\gamma'}$ for some suitable $\gamma'$, and this shows
why our original constellation of the triplet of $\gamma_{i}$ provides a sufficient
condition.

Actually as the final paragraph of the Theorem 5 there, {\em op.cit.} shows, we are doing
slightly more than this: we are, each time, applying the Lemma infinitely
often to each possible subtree of $T^{\ast}$ below some node $p_{2}$ of it
which is of length $2$, to define our strategy $\tau$ applied to moves of
length $3$. \ We then move on to the next $\Pi^{0}_{2}$ set. Although we are
applying the Lemma infinitely many times for each such $p_{2}$, and thus
infinitely many new $\Sigma_{2}$-sentences, or trees, have to be instantiated,
we had that $L_{\gamma_0}$ is a $\Sigma_{2}$-admissible set, and as the
class of such $p_{2}$ is just a set of $L_{\gamma_0}$, \
$\Sigma_{2}$-admissibility works for us to find a bound for the ranks of the
witnessing trees, as some $\delta < \gamma_0$. \ We thus can claim that our
final $\tau$ is an element of $L_{\gamma_0}$ even after $\omega$-many
iterations of this process.

$( \beta_{0} \geq \delta_0 )$ We argue for this. Let $(M,E)$ be a
non-standard model of $\ensuremath{\operatorname{KP}}$ with an infinite
nesting $( \zeta_{n} ,s_{n} )$ about $\beta_{0}$ as described. Note that
$S^{1}_{\beta_{0}}$ must be unbounded in $\beta_{0}$ (so that $L_{\beta_{0}}
\models \Sigma_{1}$-Separation), and each $\zeta_{n}$ is a limit point of
$S^{1}_{\beta_{0}}$. We do not assume that $\beta_{0}$ is
$\Sigma_{2}$-admissible (which in fact it is not as the proof shows). Let $T
\in L_{\beta_{0}}$ be a game tree. By omitting finitely much of the outer
nesting we assume $T \in L_{\zeta_{0}}$. \ We assume that Player $I$ has no
winning strategy for $G (A;T)$ in $L_{\beta_{0}}$ (for otherwise we are done).
Note that in $M$ we have that $L_{s_{0}}$ also has no winning strategy for
this game (otherwise the existence of such would reflect into $ L_{\beta_{0}}$).
We show that \pt has a winning strategy definable over $L_{\beta_{0}}$. Let
$A= \bigcup B_{n}$ with each $B_{n} \in \Pi^{0}_{2}$. For $n=0$ we apply the
argument of the Lemma using the pair $( \zeta_{1} ,s_{1} )$ in the role of $(
\gamma_0 , \gamma_{1} )$ from before, with $( \zeta_{0} ,s_{0} )$ in the
role of $( \eta_{0} , \gamma_{2} )$ described above, {\ie} we use only that
$T_{} \in L_{\zeta_{0}}$ and that $L_{\zeta_{0}} \prec_{\Sigma_{1}}
L_{s_{0}}$.

The Lemma then asserts the existence of a quasi-strategy for {\itshape{II}}
definable using the pair $( \zeta_{1} ,s_{1} )$: $T^{*} ( \varnothing )$. By
$\Sigma_{2}$-reflection the $L$-least such lies in $L_{\zeta_{1}}$, and we
shall assume that $T^{*} ( \varnothing )$ refers to it.\\

{{\em Claim: For any pair $( \zeta_{n}, s_{n} ) \text{ for } n \geq 1$ the
same tree $T^{\asterisk} ( \varnothing )$ would have resulted using this
pair.\/}}

Proof: \ Note that we can define such a tree like $T^{*} ( \varnothing )$
using such pairs, since for all of them we have that $( \zeta_{0} ,s_{0} )
\supset ( \zeta_{1} ,s_{1} ) \supset ( \zeta_{m} ,s_{m} )$ for $m>1$. As
$T^{*} ( \varnothing ) \in L_{\zeta_{1}}$ and satisfies a $\Sigma_{2}$
defining condition there, and since we also have $\zeta_{1} \in
S^{1}_{\zeta_{n}}$ for $n\geq 1$, it thus satisfies the same $\Sigma_{2}$ condition in
$L_{\zeta_{n}}$. \mbox{ }\qed\, {{\em Claim}}\\

For any position $p_{1} \in T$ with $\ensuremath{\operatorname{lh}} (p_{1} )
=1$, let $\tau (p_{1} )  $  be some arbitrary but fixed move in $T' (
\varnothing )$, this now {\itshape{II}}'s non-losing quasi-strategy for the
game $G (A,T^{*} ( \varnothing ))$ as defined in $L_{\zeta_{2}}$. \ The
relation ``$p \in T' ( \varnothing )$'' is $\Pi_{1}^{L_{\zeta_{2}}} (\{T^{*} (
\varnothing )\})$ or equivalently $\Pi_{1}^{L_{\zeta_{1}}} (\{T^{*} (
\varnothing )\})$, or indeed $\Pi_{1}^{L_{\delta}} (\{T^{*} ( \varnothing
)\})$ where $\delta$ is least in $S^{1}_{\zeta_{1}}$ above $\rho_{L} (T^{*} (
\varnothing ))$. \ Hence ``$y=T' ( \varnothing )$''  $\in
\Delta^{L_{\delta}}_{2} (\{$\/$T^{*} ( \varnothing ) \})$ and thus $T' (
\varnothing )$ also lies in $L_{\zeta_{1}}$. For definiteness we let $\tau
(p_{1} )$ be the numerically least move.

For any play, $p_{2}$ say, of length 2 consistent with the above definition
of $\tau$ so far, we apply the lemma again with $B=A_{1}$ replacing $B=A_{0}$
and with $(T^{*} ( \varnothing ))_{p_{2}}$ replacing $T$. We use the nested
pair $( \zeta_{2} ,s_{2} )$ to define quasi-strategies for {\itshape{II,}}
call them $T^{*} (p_{2} )$, one for each of the countably many $p_{2}$. These
are each definable in a $\Sigma_{2}$ way over $L_{\zeta_{2}}$, in the
parameter $(T^{*} ( \varnothing ))_{p_{2}}$. This argument uses that $(T^{*} (
\varnothing ))_{p_{2}} \in L_{\zeta_{1}} \prec_{\Sigma_{1}} L_{s_{1}}$. Let
$T' (p_{2} ) \in L_{\zeta_{2}}$ be {\itshape{II}}'s non-losing quasi-strategy
for $G (A,T^{*} (p_{2} ))$, this time with ``$y=T' (p_{2} )$''$\in
\Delta_{2}^{L_{\zeta_{2}}} (\{T^{*} (p_{2} )\})$. \ (Again these will satisfy
the same definitions as over $L_{\zeta_{m}}$ for any $m \geq 2$.) Note that we
may assume that the countably many trees $T' (p_{2} )$ appear boundedly below
$\zeta_{2}$ (using the $\Sigma_{2}$-admissibility of $\zeta_{2}$). Again for
$p_{3} \in T^{*} (p_{2} )  $ any position of length 3, let $\tau (p_{3} )  $ be
some arbitrary but fixed move in $T' (p_{2} )$. Now we consider appropriate
moves $p_{4}$ of length 4, and reapply the lemma with $B=A_{2}$ and $(T^{*}
(p_{2} ))_{p_{4}}$. \ Continuing in this way we obtain a strategy $\tau$ for
{\itshape{II}}, so that $\tau \upharpoonright\!\!\!\mbox{ }^{[1,2k+2)} \omega , $ for
$k< \omega$, is defined by a length $k$ recursion that is
$\Sigma_{2}^{L_{\zeta_{k}}} (\{T\})$.

As the argument continues, more and more of the strategy $\tau$ is defined
using successive $( \zeta_{m} ,s_{m} )$ to justify the existence of the
relevant trees in $L_{\zeta_{m}}$. {{\em Knowing\/}} that the trees are there
for the asking, we see that $\tau$ can actually be defined by a
$\Sigma_{2}$-recursion over $L_{\beta_{0}}$ in the parameter $T$ in precisely
the manner given above (the $\Sigma_{2}$-inadmissibility of $\beta_{0}$
notwithstanding).\\

If $x$ is any play consistent with $\tau$, then for every $n$, by the
defining properties of $T^{*} (p_{2n} )  $ given by the relevant application of
the lemma, $x \in \lceil T^{*} (x \upharpoonright 2n) \rceil \subseteq \neg
A_{n}$. Hence $x \notin A$, and $\tau$ is a winning strategy for{\itshape{
II}} as required. \ Thus $\beta_{0} \geq \delta_0$ is demonstrated.

$( \beta_{0} \leq \delta_0 )  $: suppose $\beta_{0} > \delta_0$. Then,
since the existence of a winning strategy for a player in any particular
$\Sigma^{0}_{3}$ game would be part of the theory $T^{1}_{\beta_{0}}
=T^{1}_{\alpha_{0}}$ where $\alpha_{0}$ is least with $L_{\alpha_{0}}
\prec_{\Sigma_{1}} L_{\beta_{0}}$, and since moreover that the existence of a
stage $\delta_0$ over which {{\em all\/}} such games have strategies,
amounts also to an existential statement, we have that $\delta_0 <
\alpha_{0}$. But this is an immediate contradiction: find a $\psi \in
T^{1}_{\alpha_{0}}$ with $\delta_0 < \alpha_{\psi} < \alpha_{0}$. But as
in the Friedman-like game of Theorem \ref{Friedman}, {\ptwo} has as winning strategy $\sigma$ to play a code for
$L_{\alpha_{\psi}}$. Hence as $\delta_0 < \alpha_{\psi}$ such a strategy and
so such a code can be found in $L_{\alpha_{\psi}}$; but  this
contradicts Tarski.  Hence $\beta_{0} \leq \delta_0$.
 {\hspace*{\fill}} {\qed} Theorem \ref{strategies}\\

As a corollary:
\begin{theorem}\label{Psicomplete}
$\Psi \equiv_1 G_3 $ and so is a $\Game
  \Sigma_3^0$-complete set.
\end{theorem}
\pf By Theorem \ref{Friedman} we just need to show that $  G_3\leq_1 \Psi $.  But we have just showed that any $\S^0_3$-game $G(A; ^{<\omega}\w)$ that is win for\po must have a winning strategy $\sigma$ for the game which is an element of $L_{\b_0}$, (otherwise the game would be a win for $\pt$ by the last theorem). The existence of such a  strategy for $\po$ is a $\Sigma_1$ sentence, which then reflects down to $L_{\a_0}$, and is an element of $\Psi$.  
\qed

\begin{theorem}\label{Thm1.2}
   If $A$ is a $\Sigma^0_3 (x)$ set so that the game $G (A)$ is won by Player
  $\mathit{I}$, then there is a generalised-ittm-recursively computable (in
  $x$) winning strategy $\sigma$ for $\mathit{I}$. That is, for some index $e$
  dependent on the definition of $A$, but not $x$, $\{ e \}^{\mathsf{E}} (x)
  \da \sigma$.
  \end{theorem}
\pf We take $x = \emp$. By the previous theorem, we have that the existence of a winning
 strategy $\s$ for $\po$ is a $\S_1$-statement, which if true, is true in $L_{\beta_0}$ and hence by $\S_1$-elementarily, also in $L_{\a_0}$. However then we can adapt the 
computation $\{e_0\}^\me$ of Section 5 above to search for a level $L_\g$ containing such a strategy $\s$, and which
 passes the check in $L_{\g^{++}}$ that it is truly winning; such a computation can then halt with $\s$ written to the OT. 
\qed
%%%%%%%%%%

\section{Discussion and some questions}

We conclude by discussing the results here, with questions concerning this kind of ittm generalised recursion, but also possible extensions of the model. There are of course, plenty of questions about the model expounded here analogous to what was discovered for the Kleeneian theory. For example:\\

\nod{\em Q Characterise the {\tmem{superjump}} for generalised ittm recursion theory.}\\

The superjump was invented by Gandy. Let $\mi$ be a general type 2 functional. \begin{definition}[Gandy \cite{Ga67a}]\mbox{ }\\
$$\begin{array}{rcl}
& & 0, 
\mbox{ if }\{e\}^\mi(\vec m, \vec x)\da\\
\mi^{\mathbbm{sJ}}((e,\vec m),\vec x) & \simeq & \\
& & 1, \mbox{ otherwise}.
\end{array}
$$
\end{definition}

This definition can be taken over verbatim to the current ittm context. Then $\mathbbm{sJ}$, regarded as a functional, is itself a type-3 functional.

\nod{\em Q Does the class of ittm semi-recursive in $\mathsf{I}$ sets have the Scale
Property (\cf \cite{KeMo72})? }\\

\nod It has the Prewellordering property by the above Lemma \ref{PWO}.
The pointclass $\Sigma^0_3$ has the scale property (\cf \cite{Kechris1973}) and, as we have $\S^0_3$-Det,  by the Third Periodicity Theorem $\Game\S^0_3$ does also (\cf \cite{Mosch4}). So for $\mi = \me$ the answer is affirmative.

\

We can summarise, and extend, the results discussed.

\begin{theorem}\label{Summary} The following characterisations of the ordinal $\beta_0$ are equivalent: it is the least ordinal that:\\
(i) supports an infinite nesting;\\
(ii)  $\power(\omega)\cap L_\beta$ form a model of {\boldmath{$\Pi^1_2$}}-Monotone Induction;\\
(iii)  any strategy for a $\Sigma^0_3$ game is definable over $L_\beta$;\\
(iv) no ITTM generalised computation $\{e\}^\mathsf{E}(m)$ computes for more than $\beta$ stages.
\end{theorem}

The new idea here is due to Hachtman \cite{Ha18} who inserted $(ii)$ and showed:  $(i) \Imp (ii) \Imp (iii)$ (and more, for example that for {\em any} $\beta$ that supported an infinite nesting, and so that $L_\beta\models V=HC$, we have that $(ii)$ holds). One should be aware that this does not claim that $L_{\beta_0}$ is $\Pi^1_2(Z)$-correct, for real parameters $Z\in L_\b$, but only that the reals of the model internally form a model of  {\boldmath{$\Pi^1_2$}}-Monotone Induction. Thus we make the definitions as follows:

\begin{definition}
$ L_\beta\models $ {\boldmath$ \Pi^1_{n+1}$}-$\sf{MI}\Equi
L_\beta\models$ 

\

\nod{ ``If $X \imp \Phi(X)$ is monotone $\wedge \,$ ``$m \in \Phi(X)$'' is $\Pi_n(m,X,Z)$ then for any 
$X\sset \omega$ and repeated applications of } $\Phi$:
 $$\all \gamma \Phi_{\gamma+1}(X)= \Phi(\bigcup_{\alpha < \gamma}\Phi_{\alpha}(X)) \imp \ex \gamma_\infty 
\Phi_{\gamma_\infty+1}(X)= \Phi_{\gamma_\infty}(X)\mbox{''}.$$
\nod  Write $\Phi_\infty(X)$ for $\Phi_{\gamma_\infty}(X)$.
\end{definition}

We also have shown (Theorem \ref{H=Psi}, Theorem\ref{Psicomplete}):

\begin{theorem}\label{summary}
$H^\mathsf{E} \equiv_1 \Sigma_1$-$Th(L_{\beta_0}) \equiv_1  G_3 \mbox{ the latter a complete } \Game\Sigma^0_3$ set.
\end{theorem}

It is natural to seek generalisations of the above. For example, consider the definition of an ordinal $\b$ supporting  $\w$-nestings of $\Sigma_{n+1}$-elementary extendibles for any fixed finite $n\geq 1$ rather than just $n=1$. Rather than write out this definition separately, we give that for a slightly stronger notion that appears to be needed in several arguments, that of   {\em strong} $\Sigma_n$ nestings.

\begin{definition}[\cite{AgWe2026} def.40]
Let $\beta$ be an ordinal. A \emph{strong $\Sigma_m$-(infinite)-nesting on $\beta$} (for $m\geq 2$) consists of a nonstandard model $M$ of $\mathsf{KP}$ and sequences $\{\zeta_i:i\in\w\}$, $\{s_i:i\in\w\}$ of $M$-ordinals such that the following hold:
\begin{enumerate}
\item $M$ extends $L_\beta$;
\item \label{StrongNestingCond2} for all $i$, $\zeta_i \leq \zeta_{i+1} < \beta$ and $s_{i+1} <^M s_i$;
\item for all $i$, $L_{\zeta_i} \prec_{\Sigma_{m}} L_{s_i}^M$;
\item for all $i$, $L_{s_{i+1}}^M \prec_{\Sigma_{m-2}} L_{s_i}^M$. 
\end{enumerate}
\end{definition}

\begin{definition} Let $\beta_n$ be the least level of $L$ on which can be based an
$\omega$-model supporting an $\omega$-depth nesting of $\Sigma_n$-extendibles.
\end{definition}

\nod Notice that for $m=2$ the notion of $\S_2$-nesting we were working with above at Def. \ref{Def3.11}, and strong $\Sigma_2$-nesting, coincide: requirement {\em 4} is superfluous when $m=2$. For $m>2$ one can show that they differ: it can be shown that the ordinal $\b_m$ is strictly greater than the least $\b_m^w$ supporting a plain (or {\em `weak'})$\S_m$-nesting. \\

One can look at the various parts that are the characterisations in Theorem \ref{Summary} and see how the linking implications can be generalised. 

\begin{definition}
A set $X$ is $1$-$\Pi^0_3$ if it is $\Pi^0_3$.  It is 
$2$-$\Pi^0_3$ if there are $\Pi^0_3$ sets $A_0,A_1$ and $X = A_0\backslash A_1$.\\
A set $X$ is $m+1$-$\Pi^0_3$ if there are  $A_1 \in m$-$\Pi^0_3$, $A_0\in \Pi^0_3$ and $X = A_0\backslash A_1$.
\end{definition}
 
We then have that $\Sigma_{n + 2}$-$\mathsf{KP}$ (but not $\Sigma_{n + 1}$-$\mathsf{KP}$, nor even $\Sigma_{n + 1}$-$\mathsf{KP}+ \Sigma_{n+1}$-Separation) proves $\tmop{Det}
(n$-$\Pi^0_3)$ \
 For $n = 1$ see {\cite{W2011}}; for $1 < n < \omega$ this is Mont\'alban-Shore
{\cite{MoSh12}}. These can be rephrased as below.
\begin{theorem}[Welch \cite{W2011}, $m=1$; Mont\'alban-Shore \cite{MoSh12}, $m >1$]
(i) {\boldmath$\Pi^1_{m+2}$}-${\mathsf{CA_0}}\vdash m$-{\boldmath$\Pi^0_m$}-$Det.$\\
(ii) {\boldmath$\Delta^1_{m+2}$}-${\mathsf{CA_0}}\not\vdash m$-{\boldmath$\Pi^0_m$}-$Det.$ 
\end{theorem}

This was then improved to: 
\begin{theorem}[Welch \cite{W2011}, $m=1$; Mont\'alban-Shore {\cite{MoSh14}}, $m > 1$]
\quad $m$-{\boldmath$\Pi^0_2$}-$Det. \vdash $ ``There exists a $\beta$-model of  
 {\boldmath$\Delta^1_{m+2}$}-${\mathsf{CA_0}}$''.

\end{theorem}

From Determinacy in the difference hierarchy we get {\em strong} nestings:
\begin{theorem}[Aguilera-Welch \cite{AgWe2026}, $\Pi^1_1$-$\mathsf{CA_0}$] \label{TheoremReversalStrong}
Let $2 \leq m$.\\ 
$m$-$\Pi^0_3$-$Det \Imp  \ex \beta( \beta $ admits   a strong $\Sigma_{m+1}$-nesting$)$ (and is a limit of such.)
\end{theorem}

{\em Q $\Sigma^0_3$-determinacy and $\game \Sigma^0_3$-monotone
inductive definitions; {\cf}{\cite{We12}}.}

\

As we have seen the architecture of an ittm, and so of generalised ittm recursions we have been discussing, is closely connected to the ideas of $\Sigma_2$ extendibility, and the fact that we have a $\Sigma_2$ definable liminf rule for cell, head position and instruction  updates in the machine.  In {\cite{FrWe11}} is developed a notion of machine using a $\S_3$, or $\S_n$ for larger $n$, definable limit rule. One can perhaps without much effort see that the generalised ittm recursion theory here, can be extended to a generalised $\S_n$-ittm recursion theory.
What is not clear is that concerning the resulting halting problem, whether this will utilise the notion of infinite $\S_n$-strong nesting, or just the generalisation of the weaker variant?\\

\nod {\em Q Does the generalised type-2 $\Sigma_n$-ittm recursion theory, require a characterisation through infinite $\S_n$-strong nestings, or just the obvious generalisation of the weaker variant?}\\

For Kleene degrees, the classes of semi-decidable sets of reals are the co-analytic sets. (There is a boldface component to the notion of reducibility here.) By work of Harrington and Steel, Det( {\boldmath$\Pi^1_1$}) is equivalent to there being only two such degrees: the  {\boldmath$\Delta^1_1$}-sets and the   {\boldmath$\Pi^1_1$}. Whilst in $L$ and set generic extensions thereof there are many incomparable semi-decidable Kleene degrees. What is the corresponding phenomenon here? \\

\nod{\em Q Which sharp of which inner model $M$ is equivalent to all ITTM-generalised semi-decidable sets of reals falling into just two equivalence classes?}\\

 We are well within {\boldmath$\Delta^1_2$} here, and so are below an inner model with a Woodin cardinal. By results of \cite{we11a} the inner model $M$ must contain a proper class of strong cardinals, and an upper bound is given by any Type 2 mouse of Feng and Jensen \cite{FeJe04}. \\

\nod{\em Q Characterise $\beta_n ?$ Is it the least level of the $L$-hierarchy over
which strategies for $n$-Boolean combinations of $\Sigma^0_3$ sets \ are
definable?}\\

%\cite{ReSl2022}

\small

\bibliographystyle{plain}
%\bibliography{settheory10v}

\ed

\end{document}